\documentclass[10pt,a4paper]{amsart}

\usepackage{mathrsfs, amsmath, amssymb, stmaryrd}%
\usepackage{hyperref}%
\usepackage[all]{xy}
\xyoption{2cell}
\UseAllTwocells

\DeclareMathOperator{\id}{id}
\DeclareMathOperator{\el}{el}
\DeclareMathOperator{\ev}{ev}
\DeclareMathOperator{\op}{op}
\DeclareMathOperator{\co}{co}
\DeclareMathOperator{\Nat}{Nat}
\DeclareMathOperator{\Lan}{Lan}
\DeclareMathOperator{\Mor}{Mor}
\DeclareMathOperator{\Mod}{\mathbf{Mod}}
\DeclareMathOperator{\Rep}{Rep}
\DeclareMathOperator{\Vect}{\mathbf{Vect}}
\DeclareMathOperator{\Coalg}{\mathbf{Coalg}}
\DeclareMathOperator{\fgp}{fgp}

\DeclareMathOperator{\len}{length}
\DeclareMathOperator{\proj}{proj}
\DeclareMathOperator{\Fil}{Fil}
\DeclareMathOperator{\MF}{MF}
\DeclareMathOperator{\colim}{colim}

\DeclareMathOperator{\Cocts}{\mathbf{Cocts}}
\DeclareMathOperator{\cat}{\mathbf{cat}}
\DeclareMathOperator{\CAT}{\mathbf{CAT}}
\DeclareMathOperator{\Set}{\mathbf{Set}}
\DeclareMathOperator{\CGTop}{\mathbf{CGTop}}

\DeclareMathOperator{\Comon}{\mathbf{Comon}}

\newcommand{\ca}[1]{\mathscr{#1}}
\newcommand{\VNat}{\ca{V}\mbox{-}\Nat}
\newcommand{\Vcat}{\ca{V}\mbox{-}\cat}
\newcommand{\VCAT}{\ca{V}\mbox{-}\CAT}
\newcommand{\Prs}[1]{\mathcal{P}\ca{#1}}
\newcommand{\Bimod}[1]{{_\ca{#1}}{\mathcal{M}}{_\ca{#1}}}
\newcommand{\modules}[2]{{_\ca{#1}}{\mathcal{M}}{_\ca{#2}}}
\newcommand{\CC}[1]{\mathbf{Comon}\left(\Bimod{#1}\right)}

\newcommand{\defl}{\mathrel{\mathop:}=}

\theoremstyle{plain}
\newtheorem{thm}[subsection]{Theorem}
\newtheorem{prop}[subsection]{Proposition}
\newtheorem{lemma}[subsection]{Lemma}
\newtheorem{cor}[subsection]{Corollary}

\theoremstyle{definition}
\newtheorem{dfn}[subsection]{Definition}

\newtheoremstyle{citing}{}{}{\itshape}{}{\bfseries}{.}{ }{\thmnote{#3}}
\theoremstyle{citing}
\newtheorem{cit}{}

\newtheoremstyle{citingdfn}{}{}{}{}{\bfseries}{.}{ }{\thmnote{#3}}
\theoremstyle{citingdfn}
\newtheorem{citdfn}{}

\title{Tannaka duality for comonoids in cosmoi}
\author{Daniel Sch\"appi}
\address{University of Chicago\\Department of Mathematics\\5734 S University Avenue\\60637 Chicago}
\email{schaeppi@math.uchicago.edu}
\keywords{Tannaka duality, comonads, structure-semantics duality}
\subjclass[2000]{18D20, 16T15}
\begin{document}

\begin{abstract}
 A classical result of Tannaka duality is the fact that a coalgebra over a field can be reconstructed from its category of finite dimensional representations by using the forgetful functor which sends a representation to its underlying vector space. There is also a corresponding recognition result, which characterizes those categories equipped with a functor to finite dimensional vector spaces which are equivalent to the category of finite dimensional representations of a coalgebra.

 In this paper we study a generalization of these questions to an arbitrary \emph{cosmos}, that is, a complete and cocomplete symmetric monoidal closed category. Instead of representations on finite dimensional vector spaces we look at representations on objects of the cosmos which have a dual. We give a necessary and sufficient condition that ensures that a comonoid can be reconstructed from its representations, and we characterize categories of representations of certain comonoids. We apply this result to certain categories of filtered modules which are used to study $p$-adic Galois representations.
\end{abstract}

\maketitle
\setcounter{tocdepth}{1}
\tableofcontents

\section{Introduction}

\subsection[GENERAL_OVERVIEW_SECTION]{}\label{GENERAL_OVERVIEW_SECTION}
 Let $k$ be a field. For any $k$-coalgebra $C$ we have an associated $k$-linear category $\Rep(C)$ of finite dimensional $C$-comodules. This category comes equipped with a forgetful functor, i.e., a $k$-linear functor $V \colon \Rep(C) \rightarrow \Vect_f$ to the category of finite dimensional $k$-vector spaces. Classical Tannakian duality for coalgebras was developed by N.\ Saavedra Rivano in \cite[Chapter~2]{SAAVEDRA}. It concerns the study of the relationship between a coalgebra $C$ and its category of representations. There are two basic questions one would like to answer:
\begin{enumerate}
 \item[(1)] The reconstruction problem: can a coalgebra be reconstructed from its category of representations?
 \item[(2)] The recognition problem: which $k$-linear functors $\omega \colon \ca{A} \rightarrow \Vect_f$ are equivalent to a forgetful functor $V \colon \Rep(C) \rightarrow \Vect_f$ for a coalgebra $C$?
\end{enumerate}

 The reconstruction problem can be solved with the following construction which was introduced independently in several places, see \cite{ULBRICH, DELIGNE, JOYAL_STREET}. For any small $k$-linear category $\ca{A}$ and any $k$-linear functor $\omega \colon \ca{A} \rightarrow \Vect_f$, the coend
\[
L(\omega) \defl \int^{A\in \ca{A}} \omega(A)\otimes \omega(A)^{\vee}
\]
 is a $k$-coalgebra. In \cite{DELIGNE, JOYAL_STREET} it is shown that for $\omega=V \colon \Rep(C) \rightarrow \Vect_f$, $L(\omega)$ is isomorphic to the coalgebra $C$. The recognition problem was solved in \cite[Proposition~2.6.3]{SAAVEDRA}: a $k$-linear functor $\omega \colon \ca{A} \rightarrow \Vect_f$ is equivalent to $V \colon \Rep\bigl(L(\omega)\bigr) \rightarrow \Vect_f$ if and only if $\ca{A}$ is a small abelian $k$-linear category and $\omega$ is faithful and exact.

 If we rephrase the universal property of coends we find that the underlying vector space of $L(\omega)$ represents the functor $\Vect \rightarrow \Vect$ which sends a vector space $X$ to the vector space $\Nat(\omega,X\otimes \omega)$ of natural transformations $\omega \Rightarrow X\otimes \omega$, where $X\otimes \omega$ stands for the functor $\ca{A}\rightarrow \Vect$ which sends $A\in \ca{A}$ to $X\otimes \omega(A)$ and a morphism $f$ to $\id_X \otimes \omega(f)$. Note that in \cite{ULBRICH} it was shown directly that the functor $\Nat(\omega,-\otimes \omega)$ is representable, without giving the explicit coend formula for the representing object. For $X=C$ a coalgebra, the natural bijection $\Vect\bigl(L(\omega), C\bigr)\cong \Nat(\omega, C\otimes \omega)$ restricts to a bijection between morphisms of coalgebras $L(\omega)\rightarrow C$ and \emph{$C$-coactions} on $\omega$, that is, natural transformations $\varrho \colon \omega \Rightarrow C\otimes \omega$ for which $\bigl(\omega(A),\varrho_A\bigr)$ is a $C$-comodule for all $A\in \ca{A}$. Such a coaction is precisely what is needed to lift $\omega$ along $V$ to a functor $\overline{\omega} \colon \ca{A} \rightarrow \Rep(C)$. We write $k\mbox{-}\cat$ for the category of small $k$-linear categories and $k$-linear functors, and we write $k\mbox{-}\cat\slash\Vect_f$ for the category of small $k$-linear categories over $\Vect_f$. The above considerations show that we have a natural bijection
\[
 \Coalg\bigl(L(\omega),C\bigr)\cong k\mbox{-}\cat\slash \Vect_f\Bigl((\ca{A},\omega),\bigl(\Rep(C),V\bigl)\Bigr)
\]
 between morphisms of coalgebras $L(\omega) \rightarrow C$ and $k$-linear functors $\overline{\omega} \colon \ca{A} \rightarrow \Rep(C)$ with $V\overline{\omega}=\omega$. In other words, the assignment which sends $(\ca{A},\omega)$ to $L(\omega)$ gives a left adjoint to the functor $\Coalg \rightarrow k\mbox{-}\cat/\Vect_f$ which sends a coalgebra $C$ to $\bigl(\Rep(C),V \bigr)$. This observation is due to R.\ Street; see \cite[\S~16]{STREET_QUANTUM_GROUPS} for a more extensive discussion. Using the terminology of adjunctions, the reconstruction and recognition results give answers to the following questions:
\begin{enumerate}
 \item[(1)]
 When is the counit $L(V) \rightarrow C$ an isomorphism?
 \item[(2)]
 When is the unit $\ca{A} \rightarrow \Rep\bigl(L(\omega)\bigr)$ an equivalence of categories?
\end{enumerate}

 This is only a special case of a class of Tannakian adjunctions which we want to generalize to other cosmoi: for any $k$-algebra $B$, the category of $B$-$B$-bimodules whose left and right $k$-actions agree is a monoidal category, so it makes sense to speak about comonoids in this category. We call such comonoids \emph{$B$-$B$-coalgebras} (they are called `$k$-cog{\`e}bro{\"i}des sur $B$' in \cite{DELIGNE}). A (left) comodule of a $B$-$B$-coalgebra $C$ is a left $B$-module $M$ together with a morphism $\varrho \colon M\rightarrow C\otimes_B M$ satisfying the evident generalizations of the usual axioms for a comodule. We get a category $\Rep(C)$ with objects the $C$-comodules whose underlying $B$-module is finitely generated and projective. This category is equipped with a $k$-linear forgetful functor $V \colon \Rep(C) \rightarrow \Mod^{\fgp}_B$. The assignment which sends a $B$-$B$-coalgebra $C$ to $V\colon \Rep(C) \rightarrow \Mod^{\fgp}_B$ again has a left adjoint
\[
 \xymatrix{
**[l]B\mbox{-}B\mbox{-}\Coalg  \rtwocell<5>^{L(-)}_{*!<0pt,2pt>+{\Rep(-)}}{`\perp} & **[r]k\mbox{-}\cat\slash\Mod^{\fgp}_{B}
}
\]
 where $L(\omega)$ is given by the coend
\[
L(\omega) \defl \int^{A\in \ca{A}} \omega(A)\otimes \omega(A)^{\vee}\rlap{,}
\]
 with left and right $B$-action induced by the left $B$-action on $\omega(A)$ and the right $B$-action on $\omega(A)^{\vee}$ respectively. We call this adjunction the \emph{Tannakian adjunction}. We recover the original Tannakian adjunction if we let $B=k$. We will show that an analogous adjunction still exists if we replace the category $\Vect$ of $k$-vector spaces by any complete and cocomplete symmetric monoidal closed category (a similar generalized adjunction was already considered in \cite[\S~16]{STREET_QUANTUM_GROUPS}). Following B\'enabou and Kelly, we call such a category a \emph{cosmos} (cf.\ \cite[\S~4]{BENABOU} and \cite{KELLY_COSMOS}).

\subsection[F_MODULES_INTRO_SECTION]{}\label{F_MODULES_INTRO_SECTION}
 We now describe an example which was the motivation for studying the reconstruction problem in a cosmos other than $\Vect$. Let $k$ be a perfect field of characteristic $p>0$, and let $W$ be the ring of Witt-vectors with coefficients in $k$. The ring $W$ is a complete discrete valuation ring with uniformizer $p$ and residue field $k$ (see \cite[\S\S~II.5-II.6]{SERRE}). We write $W_n$ for the quotient ring $W/p^n W$. In \cite{FONTAINE_LAFFAILLE}, J.-M.\ Fontaine and G.\ Laffaille defined the category $\MF_{fl}$, which consists of filtered $W$-modules of finite length with some additional structure (see Section~\ref{F_MODULES_SECTION} for a precise definition). In \cite{WINTENBERGER}, the objects of $\MF_{fl}$ are called `F-modules filtr\'es sur $W$'. The category $\MF_{fl}$ is abelian and $\mathbb{Z}_p$-linear, and the forgetful functor $\omega \colon \MF_{fl}\rightarrow \Mod_W$ is a $\mathbb{Z}_p$-linear functor. We let $\MF_{fl}^n$ be the full subcategory of $\MF_{fl}$ consisting of those filtered $F$-modules whose underlying $W$-module is annihilated by $p^n$. The category $\MF^n_{fl}$ is $\mathbb{Z}/p^n \mathbb{Z}$-linear, and the forgetful functor gives a $\mathbb{Z}/p^n \mathbb{Z}$-linear functor $\omega \colon \MF^n_{fl} \rightarrow \Mod_{W_n}$. By tensoring the hom-modules with the $\mathbb{Z}/p^n \mathbb{Z}$-algebra $W_n$ we obtain a $W_n$-linear category $\MF^n_{fl}\otimes W_n$. The $\mathbb{Z}/p^n \mathbb{Z}$-linear functor $\omega$ induces a $W_n$-linear functor $\MF^n_{fl} \otimes W_n \rightarrow \Mod_{W_n}$. R.\ Pink asked the following question: is it possible to apply constructions similar to the ones found in the theory of Tannakian categories developed in \cite{DELIGNE} to the functor $\MF^n_{fl} \otimes W_n \rightarrow \Mod_{W_n}$ to obtain an affine group scheme over $W_n$?

 In Section~\ref{F_MODULES_SECTION} we study a different but related question. Instead of looking at the $W_n$-linear functor $\MF^n_{fl} \otimes W_n \rightarrow \Mod_{W_n}$ obtained by base change, we study the $\mathbb{Z}/p^n \mathbb{Z}$-linear forgetful functor $\omega \colon \MF^n_{fl} \rightarrow \Mod_{W_n}$ itself. We write $\MF^n_{\proj}$ for the full subcategory of $\MF^n_{fl}$ consisting of objects whose underlying $W_n$-module is finitely generated and projective. In Section~\ref{F_MODULES_SECTION} we will prove the following theorem.

\begin{cit}[Theorem~\ref{FILTERED_MODULE_THEOREM}]
 The category $\MF^n_{\proj}$ is equivalent to the category of those left comodules of the $W_n$-$W_n$-coalgebra
\[
 L = \int^{\MF^n_{\proj}} \omega(M) \otimes_{\mathbb{Z}/p^n \mathbb{Z}} \omega(M)^\vee
\]
 whose underlying $W_n$-module is finitely generated projective. The right action on $L$ is induced by the $W_n$-actions on $\omega(M)^\vee$, and the left action is induced by the $W_n$-actions on $\omega(M)$. The $W_n$-$W_n$-coalgebra $L$ is flat as right $W_n$-module.
\end{cit}

 We prove this result by showing that the forgetful functor $\omega$ satisfies certain criteria which ensure that the unit of the generalized Tannakian adjunction is an equivalence.

\subsection[COSMOS_SECTION]{}\label{COSMOS_SECTION}
 In order to define the Tannakian adjunction for arbitrary cosmoi, we first have to find the appropriate generalization of `finite dimensional vector space'. For our purposes, the appropriate notion is `object with a dual'. The finite dimensional vector spaces are precisely the objects with duals in $\Vect$, so this is a natural generalization of the classical case. Another reason for our choice is the following: if our comonoid is the underlying comonoid of a Hopf monoid, we want our representations to have duals, and a comodule of a Hopf monoid has a dual if and only if its underlying object has a dual (see \cite[Proposition~15.1]{STREET_QUANTUM_GROUPS}). We will eventually be interested in the reconstruction of Hopf monoids instead of mere comonoids, and the existence of duals in the category of representations is crucial for the reconstruction of the antipode map of the Hopf monoid (see \cite[\S~16]{STREET_QUANTUM_GROUPS}). In order to apply reconstruction results for Hopf monoids similar to those found in \cite{STREET_QUANTUM_GROUPS}, we should therefore ask the following question: is it possible to reconstruct a comonoid $L$ from the category of those $L$-comodules for which the underlying objects have duals?

 For certain classes of cosmoi this question has already been studied: T.\ Wedhorn studied the reconstruction problem for Hopf algebras over Dedekind rings, and the recognition problem for valuation rings (see \cite{WEDHORN}). B.\ Day solved both problems for finitely presentable cosmoi for which the full subcategory of objects with duals is closed under finite limits and colimits (see \cite{DAY}). P.\ McCrudden used a result of B.\ Pareigis to solve the reconstruction problem for \emph{Maschkean} categories, which are certain abelian monoidal categories in which all monomorphisms split (see \cite{PAREIGIS}, \cite{MCCRUDDEN_MASCHKE}). All these approaches make the assumption that the category of objects with duals is closed under finite limits. But an $R$-module has a dual if and only if it is finitely generated and projective, and a kernel of a morphism between projective modules is in general not projective; therefore, the above results cannot be applied to the case where $\ca{V}$ is the cosmos $\Mod_R$ of $R$-modules for a general commutative ring $R$, such as the case of the example described in Section~\ref{F_MODULES_INTRO_SECTION}. In the present paper we study the reconstruction and the recognition problem without assuming that the category of objects with duals is closed under finite limits. We succeed in giving a necessary and sufficient condition for solving the reconstruction problem (see Theorems~\ref{RECONSTRUCTION_THM} and \ref{DENSE_RECONSTRUCTION_THM}), and we provide a partial solution of the recognition problem (see Theorems~\ref{RECOGNITION_THM} and \ref{DENSE_RECOGNITION_THM}). We now turn to a discussion of these results.

\subsection[DESCRIPTION_OF_RESULTS_SECTION]{}\label{DESCRIPTION_OF_RESULTS_SECTION}
 We fix a cosmos $\ca{V}$, and we let $\ca{V}^c$ be the full subcategory of objects which have a dual. In order to construct the Tannakian adjunction, we have to choose a cosmos for which $\ca{V}^c$ is an (essentially) small category (i.e., it has only a set of isomorphism classes). For a comonoid $T$, we let $\ca{V}_T^c$ be the category of $T$-comodules whose underlying objects have duals. Such comodules are called \emph{Cauchy comodules} (cf.\ \cite[Proposition~10.6]{STREET_QUANTUM_GROUPS}; see Section~\ref{CAUCHY_COMPLETION_SECTION} for a motivation of this terminology). It is well known that in the case where $\ca{V}$ is the cosmos of modules over some commutative ring $R$, $\ca{V}_T^c$ is an $R$-linear category. For general cosmoi, it is still true that $\ca{V}_T^c$ is enriched in $\ca{V}$. We will use this additional structure on the category of comodules for our constructions. The reader who is unfamiliar with the theory of enriched categories should find enough background material in Section~\ref{PRELIMINARIES_SECTION} to follow the arguments in the special case where $\ca{V}$ is the cosmos $\Mod_R$ of $R$-modules for some commutative ring $R$. One advantage of working in full generality is that we also cover the cases where $\ca{V}$ is the cosmos of $A$-graded $R$-modules for any abelian group $A$, or of modules or comodules of an $R$-bialgebra.

 The \emph{comodule functor} $\ca{V}^c_{(-)} \colon \Comon(\ca{V}) \rightarrow \Vcat\slash \ca{V}^c$ sends a comonoid $T$ in $\ca{V}$ to the $\ca{V}$-category $\ca{V}_T^c$, equipped with the restriction $V^c_T \colon \ca{V}_T^c \rightarrow \ca{V}^c$ of the forgetful functor $V_T$. In Section~\ref{TANNAKIAN_ADJUNCTION_SECTION} we show that the comodule functor has a left adjoint. The resulting adjunction
\[
\xymatrix{
**[l]\Comon(\ca{V})  \rtwocell<5>^{L(-)}_{*!<0pt,2pt>+{\ca{V}_{(-)}^c}}{`\perp} & **[r]\Vcat\slash\ca{V}^c
}
\]
 is called the \emph{Tannakian adjunction}. The existence of such an adjunction is not new (cf.\ \cite[Section~16]{STREET_QUANTUM_GROUPS} and \cite{MCCRUDDEN_REPR_COALGEBROIDS}, where similar adjunctions are defined), but we give a new construction: in Section~\ref{TANNAKIAN_ADJUNCTION_SECTION} we show that the Tannakian adjunction can be written as a composite of two partial adjunctions. One of them is very basic; it exists in all categories with pullbacks. The other is a form of the so-called `semantics-structure adjunction', which we explain in Section~\ref{SEMANTICS_STRUCTURE_SECTION}. For a $\ca{V}$-functor $\omega \colon \ca{A} \rightarrow \ca{V}^c$ with small domain, we denote the $(\ca{A},\omega)$-component of the unit of the Tannakian adjunction by
\[
 \xymatrix{ \ca{A} \ar[rr]^N \ar[rd]_{\omega} && \ca{V}^c_{L(\omega)} \ar[ld]^{V^c_{L(\omega)}}\\
 & \ca{V}^c}
\]
 and for a comonoid $T$ we write $\nu \colon L(V^c_T) \rightarrow T$ for the $T$-component of the counit.

 Using this new construction, we study the counit and the unit of the Tannakian adjunction in Sections~\ref{RECONSTRUCTION_SECTION} and \ref{RECOGNITION_SECTION} respectively. The statements of our theorems become less technical if we make some additional assumptions about the cosmos $\ca{V}$. Let $\ca{C}$ be a category, and let $\ca{A} \subseteq \ca{C}$ be a subcategory. For every object $C\in \ca{C}$, we have the domain functor $D \colon \ca{A}\slash C \rightarrow \ca{C}$ which sends an object $\varphi \colon A \rightarrow C$ of $\ca{A}\slash C$ to its domain. There is a \emph{tautological natural transformation} $\tau \colon D \rightarrow \underline{C}$ from $D$ to the constant functor at $C$. The component of $\tau$ at $\varphi \colon A \rightarrow C$ is $\varphi$ itself.

\begin{dfn}\label{SET_DENSITY_DFN}
 A subcategory $\ca{A}\subseteq \ca{C}$ is called \emph{dense} if for all $C\in \ca{C}$, the tautological natural transformation $\tau \colon D \Rightarrow \underline{C}$ exhibits $C$ as colimit of the domain functor $D \colon \ca{A}\slash C \rightarrow \ca{C}$. There is a generalization of this notion to categories enriched in $\ca{V}$ (see Section~\ref{DENSE_FUNCTORS_SECTION} for a definition). Any ordinary category is a category enriched in the cosmos $\Set$ of sets, so we generally speak of $\Set$-dense subcategories to avoid confusion.
\end{dfn}

If we consider a small category $\ca{A}$ as a subcategory of the functor category $[\ca{A}^{\op}, \Set]$ via the Yoneda embedding, then $\ca{A} \subseteq [\ca{A}^{\op},\Set]$ is an example of a dense subcategory. For any commutative ring $R$, the full subcategory of $\Mod_R$ consisting of the single object $R\oplus R$ is another example of a dense subcategory.

\begin{citdfn}[Definition~\ref{DAG_DFN}]
 An essentially small full subcategory $\ca{X} \subseteq \ca{V}$ of a cosmos $\ca{V}$ is called a \emph{dense autonomous generator} if $\ca{X}$ consists of objects with duals, is closed under the tensor product and under the formation of duals, and is $\Set$-dense in $\ca{V}$.
\end{citdfn}

 We use this terminology because a dense autonomous generator $\ca{X}$ of $\ca{V}$ is an \emph{autonomous monoidal category}, i.e., all the objects of $\ca{X}$ have duals. For example, for $R$ a commutative ring, the finitely generated free $R$-modules form a dense autonomous generator of $\Mod_R$. If $\ca{V}$ is the cosmos of chain complexes of $R$-modules, then the full subcategory consisting of bounded chain complexes of finitely generated free $R$-modules is a dense autonomous generator. Note that if $\ca{V}$ has an essentially small $\Set$-dense subcategory $\ca{A}$ consisting of objects with duals, then the closure of $\ca{A}$ under duals and tensor products is a dense autonomous generator.

 The classical reconstruction theorems all rely on the fact that the comonoid in question is a union of Cauchy comodules (see e.g.\ \cite{JOYAL_STREET, MCCRUDDEN_MASCHKE, WEDHORN}). A union is of course just a special case of a colimit. In Theorem~\ref{DENSE_RECONSTRUCTION_THM}, we prove that we get a necessary and sufficient condition for the solution of the reconstruction problem when we consider colimits of diagrams consisting of arbitrary morphisms instead of monomorphisms. A special case of Theorem~\ref{DENSE_RECONSTRUCTION_THM} is the following result.

\begin{thm}\label{NEUTRAL_RECONSTRUCTION_PROBLEM}
 Let $\ca{V}$ be a cosmos which has a dense autonomous generator, and let $T$ be a comonoid in $\ca{V}$. Then the counit $\nu \colon L(V^c_T) \rightarrow T$ of the Tannakian adjunction is an isomorphism if and only if the tautological natural transformation exhibits $T$, considered as a comodule over itself, as the colimit of the diagram $D \colon \ca{V}^c_T \slash T \rightarrow \ca{V}_T$ of Cauchy comodules over $T$.
\end{thm}

 In Section~\ref{RECONSTRUCTION_REMARKS_SECTION} we will see that there are cosmoi for which the reconstruction theorem does not hold, i.e., for which the counit of the Tannakian adjunction is not an isomorphism. The question whether the reconstruction theorem holds in cosmoi of modules over an arbitrary commutative ring $R$ remains open.

 The classical recognition result concerns \emph{exact} $k$-linear functors. It turns out that we have to generalize left and right exactness separately. Right exactness concerns the preservation of finite colimits. In the general case, we no longer assume that $\ca{V}^c$ is closed under finite colimits, but there is a certain class of colimits which plays an important role in the recognition problem.

\begin{dfn}\label{OMEGA_RIGID_DFN}
 Let $\omega \colon \ca{A} \rightarrow \ca{V}^c$ be a $\ca{V}$-functor, and let $\ca{D}$ be a small ordinary category. A diagram $D \colon \ca{D} \rightarrow \ca{A}$ is called \emph{$\omega$-rigid} if the colimit of $\omega D \colon \ca{D} \rightarrow \ca{V}$ lies in $\ca{V}^c$.
\end{dfn}

 Note that $\omega$-rigid diagrams need not be finite. For any comonoid $T\in \ca{V}$, if we take $\omega=V^c_T \colon \ca{V}^c_T \rightarrow \ca{V}$, then $\ca{V}^c_T$ has all colimits of $\omega$-rigid diagrams, because colimits in $\ca{V}_T$ are computed in the same way as in $\ca{V}$. Any $(\ca{A},\omega)$ for which the unit of the Tannakian adjunction is an equivalence must therefore satisfy a similar property.

 We now turn to the generalization of left exactness. To do this, we need to introduce the \emph{category of elements} of a functor $F \colon \ca{A} \rightarrow \Set$ between ordinary categories. This category is denoted by $\el(F)$, and its objects are pairs $(A,x)$, where $A \in \ca{A}$ and $x\in FA$. The morphisms $(A,x) \rightarrow (A^\prime, x^\prime)$ are given by the morphisms $f\colon A \rightarrow A^\prime$ in $\ca{A}$ with $Ff(x)=x^\prime$. If $\ca{A}$ is an abelian $k$-linear category, then a functor $\omega \colon \ca{A} \rightarrow \Mod_k$ is left exact if and only if the category $\el(\omega)$ is \emph{cofiltered}, i.e., if and only if
\begin{itemize}
 \item the category $\el(\omega)$ is non-empty;
 \item for any two objects $A,A^\prime \in \el(\omega)$, there is an object $B \in \el(\omega)$ together with morphisms $B \rightarrow A$, $B \rightarrow A^\prime$; and
 \item for any two morphisms $f,g \colon A \rightarrow A^\prime$ in $\el(\omega)$, there is an object $B$ and a morphism $h\colon B \rightarrow A$ such that $fh=gh$.
\end{itemize}
 More generally, if $\ca{C}$ is an ordinary category with finite limits, then a functor $F \colon \ca{A} \rightarrow \Set$ preserves finite limits if and only if $\el(F)$ is cofiltered. Even if we no longer assume that $\ca{C}$ has finite limits, we can still talk about functors $F$ for which the category $\el(F)$ is cofiltered. Such functors are called \emph{flat} functors.

 In order to state the result, we need the \emph{forgetful functor} $V=\ca{V}(I,-) \colon \ca{V} \rightarrow \Set$, which sends an object $M$ of $\ca{V}$ to its `underlying set', i.e., to the set of maps $I \rightarrow M$, where $I$ denotes the unit of the monoidal category $\ca{V}$. In the case $\ca{V}=\Mod_R$, this functor coincides with the usual forgetful functor to the category of sets. If $\ca{V}$ is the category of graded $R$-modules, then $V \colon \ca{V} \rightarrow \Set$ sends a graded module $(M_i)_{i\in\mathbb{Z}}$ to the set of elements of $M_0$. Using the forgetful functor we can define the underlying ordinary category $\ca{A}_0$ of a $\ca{V}$-category $\ca{A}$, and the underlying ordinary functor $F_0$ of a $\ca{V}$-functor $F$. In the case $\ca{V}=\Mod_R$, this just means that we forget about the scalar multiplication and addition of morphisms.

Our results concerning the recognition problem require the additional assumption that the category $\ca{V}$ is locally finitely presentable. If $\ca{V}$ has a dense autonomous generator, this is equivalent to the fact that the forgetful functor $V \colon \ca{V} \rightarrow \Set$ preserves filtered colimits (see the proof of Theorem~\ref{DENSE_RECOGNITION_THM}). This holds for example in the cosmoi of (graded) $R$-modules, chain complexes of $R$-modules, etc. A very general result showing that `algebraic' categories are locally finitely presentable is \cite[Corollary~3.7]{ADAMEK_ROSICKY}.

\begin{thm}\label{NEUTRAL_RECOGNITION_THEOREM}
 Let $\ca{V}$ be a locally finitely presentable cosmos which has a dense autonomous generator $\ca{X}$. Let $\ca{A}$ be a $\ca{V}$-category which has tensor products\footnote{Tensor products are also known as \emph{copowers}; see Section~\ref{WEIGHTED_COLIMIT_SECTION} for a definition.} with objects in $\ca{X}$, and let $\omega \colon \ca{A} \rightarrow \ca{V}$ be a $\ca{V}$-functor such that $\omega(A)\in \ca{V}^c$ for all $A \in \ca{A}$. The unit $N \colon \ca{A} \rightarrow \ca{V}^c_{L(\omega)}$ of the Tannakian adjunction is an equivalence if
\begin{enumerate}
 \item[i)]
 the functor $\omega$ reflects isomorphisms;
 \item[ii)]
 the category $\el(V\omega_0)$ of elements of $V\omega_0 \colon \ca{A}_0 \rightarrow \Set$ is cofiltered; and
 \item[iii)]
 the underlying category $\ca{A}_0$ of $\ca{A}$ has colimits of $\omega$-rigid diagrams, and $\omega_0 \colon \ca{A}_0 \rightarrow \ca{V}_0$ preserves them.
\end{enumerate}
\end{thm}

 We prove a more general result in Theorem~\ref{DENSE_RECOGNITION_THM}. In Section~\ref{SEMANTICS_STRUCTURE_SECTION} we will see that conditions i) and iii) are necessary conditions: if the unit $ N \colon \ca{A} \rightarrow \ca{V}^c_{L(\omega)}$ of the Tannakian adjunction is an equivalence of categories, then $\omega$ automatically satisfies i) and iii). Condition ii), on the other hand, is rather strong. It implies for example that the comonoid $L(\omega)$ associated to $(\ca{A},\omega)$ is \emph{flat}, i.e., that $L(\omega)\otimes- \colon \ca{V} \rightarrow \ca{V}$ preserves finite limits (in the sense of \cite{KELLY_FINLIM}). It would be interesting to know if condition ii) is necessary in the case of flat comonoids. Is it true that the category of elements $\el(V\omega_0)$ is cofiltered if $\omega \colon \ca{V}^c_L\rightarrow \ca{V}$ is the forgetful functor for a flat comonoid $L$? This question is still open.

 As in the classical case outlined in Section~\ref{GENERAL_OVERVIEW_SECTION}, there is a more general form of the Tannakian adjunction, and the example mentioned in Section~\ref{F_MODULES_INTRO_SECTION} requires this additional level of generality. The generalized adjunction concerns comonoids in the category of \emph{modules} from $\ca{B}$ to $\ca{B}$, where $\ca{B}$ is a fixed $\ca{V}$-category (see Section~\ref{COMODULE_FUNCTOR_OVERVIEW_SECTION} for a definition). Taking $\ca{V}=\Vect$ and $\ca{B}$ a $k$-linear category with one object (i.e., a $k$-algebra $B$), our adjunction reduces to the adjunction from Section~\ref{GENERAL_OVERVIEW_SECTION}. There are other examples related to other generalizations of the Tannakian adjunction. For example, if the cosmos $\ca{V}$ is an additive category, and if we let $\ca{B}$ be the free $\ca{V}$-category on the finite set $X$, considered as a category with no non-identity arrows, then a comonoid in the category of modules from $X$ to $X$ is the same as a $\ca{V}^{\op}$-category whose set of objects is $X$. A Tannakian adjunction for $\ca{V}^{\op}$-categories (with varying set of objects) was defined in \cite{DAY_STREET} and \cite{MCCRUDDEN_REPR_COALGEBROIDS}. Here the assumptions that $X$ is finite and that $\ca{V}$ is additive are crucial; the above characterization of comonoids relies on the fact that coproducts indexed by $X$ are isomorphic to the corresponding products.

\subsection[OUTLINE_SECTION]{}\label{PAPER_OUTLINE_SECTION}
 We now give a short outline of this paper. In Section~\ref{PRELIMINARIES_SECTION}, we provide some background material from enriched category theory. The purpose of this section is twofold: it is intended to serve as a brief introduction to enriched categories for people who are less familiar with them and to fix some notations and give names to certain maps which we will use in later parts of the paper. In Section~\ref{DENSITY_SECTION}, we also use some concepts from the theory of locally finitely presentable enriched categories, which was developed in \cite{KELLY_FINLIM}. Throughout the paper we will often use the formalism of pasting composites of $\ca{V}$-natural transformations which was introduced in \cite{KELLY_STREET}. We briefly explain this in Appendix~\ref{PASTING_COMPOSITES_APPENDIX}.

 In Section~\ref{CATEGORIES_OF_COMODULES_SECTION}, we define enriched categories of comodules for comonads. These are then used to define the \emph{comodule functor} in Section~\ref{COMODULE_FUNCTOR_SECTION}, which is the right adjoint of the Tannakian adjunction. To show that the comodule functor has a left adjoint, we exhibit it as a composite of the right adjoints of two partial adjunctions. We make use of the following observation: a comonoid $T$ in $\ca{V}$ gives rise to a comonad on $\ca{V}$, given by the $\ca{V}$-functor $T\otimes - \colon \ca{V} \rightarrow \ca{V}$. This construction gives an equivalence between the categories of comonoids in $\ca{V}$ and the category of cocontinuous comonads on $\ca{V}$.

 The first partial adjunction is the \emph{semantics-structure adjunction}, which we describe in Section~\ref{SEMANTICS_STRUCTURE_SECTION}; it is dual to the one found in in \cite{DUBUC}, and describes a relationship between comonads on $\ca{V}$ and $\ca{V}$-categories over $\ca{V}$. Since we are interested in finding descriptions of the unit and the counit, we have to be quite explicit about the involved natural isomorphisms. The second partial adjunction relates $\ca{V}$-categories over $\ca{V}$ to $\ca{V}$-categories over $\ca{V}^c$. We introduce it in Section~\ref{TANNAKIAN_ADJUNCTION_SECTION}, where we also combine the two to get a new construction of the Tannakian adjunction.

 In Section~\ref{RECONSTRUCTION_SECTION}, we use this new construction of the Tannakian adjunction to prove our reconstruction result. Basically, we have to find a convenient description of the counit of the Tannakian adjunction, and this can be done in a straightforward way using the pasting composites from \cite{KELLY_STREET}.

 The main result of the paper is proved in Section~\ref{RECOGNITION_SECTION}, where we give conditions which ensure that a category is equivalent to a category of Cauchy comodules for some comonoid. The proof of this result is more involved. The basic strategy behind the proof is the following. If $\omega \colon \ca{A} \rightarrow \ca{V}$ \emph{is} the forgetful functor of the $\ca{V}$-category of Cauchy $T$-comodules for some comonoid $T$, then there is a natural embedding $K \colon \ca{A} \rightarrow \ca{V}_T$ of $\ca{A}$ in the $\ca{V}$-category of \emph{all} $T$-comodules, and this embedding is compatible with the forgetful $\ca{V}$-functor $V_T \colon \ca{V}_T \rightarrow \ca{V}$. The $\ca{V}$-category $\ca{V}_T$ is cocomplete in the enriched sense (see Proposition~\ref{COLIMIT_CREATION_PROP}), so we get an induced cocontinuous $\ca{V}$-functor $L_K \colon \Prs{A} \rightarrow \ca{V}_T$ from the \emph{free cocompletion} $\Prs{A}$ of $\ca{A}$. Similarly, the forgetful functor $\omega \colon \ca{A} \rightarrow \ca{V}$ induces a cocontinuous $\ca{V}$-functor $L_\omega \colon \Prs{A} \rightarrow \ca{V}$. The situation is summarized by the diagram
\[
 \xymatrix{\ca{A} \ar[r]^{K} \ar[rd]_{\omega} & \ca{V}_T \ar[d]^{V_T} & \ar[l]_{L_K} \Prs{A} \ar[ld]^{L_\omega}\\
& \ca{V}}
\]
 which commutes up to natural isomorphism. For certain comonoids $T$, the category $\ca{V}_T$ is a `quotient' of $\Prs{A}$. More precisely, the $\ca{V}$-functor $L_K \colon \Prs{A} \rightarrow \ca{V}_T$ is sometimes a \emph{localization} of $\Prs{A}$; i.e., the category $\ca{V}_T$ is in some cases obtained from $\Prs{A}$ by formally inverting the morphisms in $\Prs{A}$ which get sent to isomorphisms by $L_\omega$. This happens for example if the underlying object of $T$ is itself a Cauchy object.

 If we start with an arbitrary $\ca{V}$-category $\ca{A}$ and a $\ca{V}$-functor $\omega \colon \ca{A} \rightarrow \ca{V}$, we still get an induced $\ca{V}$-functor $L_\omega \colon \Prs{A} \rightarrow \ca{V}$ from the free cocompletion of $\ca{A}$ to $\ca{V}$. If $\ca{V}$ is locally finitely presentable, we can formally invert the class $\Sigma$ of morphisms in $\Prs{A}$ which get sent to isomorphisms by $L_\omega$. We prove our reconstruction result by giving conditions on the $\ca{V}$-functor $\omega \colon \ca{A} \rightarrow \ca{V}$ which ensure that the localized category $\Prs{A}[\Sigma^{-1}]$ is of the form $\ca{V}_T$ for some comonoid $T$, and that the composite
\[
 \xymatrix{\ca{A} \ar[r]^-Y & \Prs{A} \ar[r] & \Prs{A}[\Sigma^{-1}] \ar[r]^-{\simeq} & \ca{V}_T}
\]
induces an equivalence between the $\ca{V}$-category $\ca{A}$ and the $\ca{V}$-category of Cauchy comodules of the comonoid $T$. We use a form of Beck's monadicity theorem to prove that $\Prs{A}[\Sigma^{-1}]$ is equivalent to $\ca{V}_T$ for some comonoid $T$.

 In Section~\ref{DENSITY_SECTION}, we specialize both the reconstruction result and the recognition result to cosmoi with dense autonomous generator. We use this in Section~\ref{F_MODULES_SECTION} to prove our result about the categories $\MF^n_{\proj}$ of filtered F-modules mentioned in Section~\ref{F_MODULES_INTRO_SECTION}.

 In Section~\ref{TANNAKIAN_BIADJUNCTION_SECTION}, we give a short outline of how one could generalize the reconstruction and recognition results found in this paper to comonoids with additional structure. We approach this problem from a categorical perspective.

\section*{Acknowledgments}
 This paper contains generalizations of results of my master's thesis, which was written under the advice of Prof.\ Richard Pink at ETH Z\"urich, Switzerland. I thank Peter May, Richard Pink, Mike Shulman and Ross Street for kindly answering questions and for giving suggestions for improvement. I am especially grateful to Mike Shulman, who pointed out that my original construction of the Tannakian adjunction arises as a composite of two partial adjunctions. This viewpoint greatly clarifies certain proofs. I thank Emily Riehl and Claire Tomesch for their help with editing the paper.

\section{Preliminaries}\label{PRELIMINARIES_SECTION}

\subsection[ENRICHED_CATEGORIES_SECTION]{}\label{ENRICHED_CATEGORIES_SECTION}
 A \emph{cosmos} is a complete and cocomplete symmetric monoidal closed category. We fix a cosmos $\ca{V}$, with tensor product $-\otimes - \colon \ca{V}\times \ca{V} \rightarrow \ca{V} $, unit object $I\in \ca{V}$ and internal hom $[-,-]$. A category $\ca{A}$ \emph{enriched} in $\ca{V}$ has objects $A, A^\prime, \ldots$ and instead of hom-sets, it has \emph{hom-objects} $\ca{A}(A,A^\prime) \in \ca{V}$. The standard source for the theory of enriched categories is \cite{KELLY_BASIC}. The basic concepts of category theory can be generalized to this context. For example, for a small $\ca{V}$-category $\ca{A}$, there is a $\ca{V}$-category $\Prs{A}$ of enriched presheaves on $\ca{A}$ (i.e., $\ca{V}$-functors $\ca{A}^{\op} \rightarrow \ca{V}$), and a corresponding Yoneda embedding. We denote the category of small $\ca{V}$-categories and $\ca{V}$-functors by $\Vcat$, and we write $\VCAT$ for the (very large) $\ca{V}$-category of all large $\ca{V}$-categories and $\ca{V}$-functors. The reader who is unfamiliar with the general theory of enriched categories should keep in mind the case $\ca{V}=\Mod_R$, $R$ a commutative ring, where $\ca{V}$-category, $\ca{V}$-functor and $\ca{V}$-natural transformation correspond to the notions of $R$-linear category, $R$-linear functor and ordinary natural transformation respectively. Note that we do not require that an $R$-linear category has finite direct sums. Most of the general concepts are self-explanatory in this context. In the next few sections we outline some of those which are not.

\subsection[WEIGHTED_COLIMIT_SECTION]{}\label{WEIGHTED_COLIMIT_SECTION}
 When we enrich the notion of colimits, we naturally arrive at the concept of a \emph{weighted colimit}\footnote{Weighted colimits are called \emph{indexed colimits} in \cite{KELLY_BASIC}.}: an object $K$ of a $\ca{V}$-category $\ca{E}$ is said to be the colimit of $G\colon \ca{D} \rightarrow \ca{E}$ weighted by  $J \colon \ca{D}^{\op} \rightarrow \ca{V}$ if there is an isomorphism
\[
\xymatrix{\ca{E}(K,E) \ar[r]^-{\varphi_E} & \Prs{D}\bigl(J,\ca{E}(G-,E)\bigr)}
\]
 of $\ca{V}$-functors which is $\ca{V}$-natural in $E$. The object $K$ is usually denoted by $J \star G$. For any small $\ca{V}$-category $\ca{B}$, the category $\Prs{B}$ of enriched presheaves on $\ca{B}$ has all weighted colimits (see \cite[\S~3.3]{KELLY_BASIC}). The identity of $J \star G$ corresponds under $\varphi$ to the \emph{unit}
\[
 \xymatrix{J \ar[r]^-{\lambda} & \ca{E}(G-,J\star G)}
\]
 of $J \star G$, which has the property that for any $\ca{V}$-natural transformation $\alpha \colon J \Rightarrow \ca{E}(G-,E)$, there is a unique morphism $a \colon J\star G \rightarrow E$ such that $\alpha = \ca{E}(G-,a)$. For $\ca{V}=\Mod_R$, the existence of a $\lambda$ with this property is equivalent to the existence of the natural isomorphism $\varphi$ (see \cite[\S~3.1]{KELLY_BASIC}). In particular, if $L \colon \ca{E} \rightarrow \ca{E}^{\prime}$ is a $\ca{V}$-functor such that the colimit $J \star LG$ exists, there is a unique morphism $\widehat{L} \colon J \star LG \rightarrow L(J\star G)$ for which the diagram
\[
 \xymatrix{J  \ar[rr]^-{\lambda} \ar[d]_-{\lambda^{\prime}} && \ca{E}(G-,J\star G) \ar[d]^-{L} \\
 \ca{E}^{\prime}(LG-,J\star LG) \ar[rr]_-{\ca{E}(LG-,\widehat{L})} && \ca{E}^{\prime}\bigl(LG-,L(J\star G)\bigr)}
\]
 is commutative, where $\lambda^{\prime}$ denotes the unit of $J \star LG$. The morphism $\widehat{L} \colon J\star LG \rightarrow L(J\star G)$ is called the \emph{comparison morphism}, and we say that $L$ \emph{preserves} the colimit $J \star G$ if $J\star LG$ exists and $\widehat {L}$ is an isomorphism. A $\ca{V}$-functor is said to be \emph{cocontinuous} if it preserves all small weighted colimits that exist. A $\ca{V}$-functor $L\colon \ca{E} \rightarrow \ca{E}^{\prime}$ is called a \emph{left $\ca{V}$-adjoint} or simply \emph{left adjoint} if there is a $\ca{V}$-functor $R\colon \ca{E}^{\prime} \rightarrow \ca{E}$ and $\ca{V}$-natural transformations $\eta \colon \id\Rightarrow RL$ and $\varepsilon \colon LR \Rightarrow \id$ satisfying the usual triangular identities. Recall that we get underlying ordinary categories, functors and natural transformations if we apply the forgetful functor $V=\ca{V}(I,-)\colon \ca{V} \rightarrow \Set$ to the hom-objects of a $\ca{V}$-category. The condition that $L$ is a left $\ca{V}$-adjoint is in general stronger than saying that the underlying ordinary functor $L_0$ is a left adjoint, but if $\ca{V}=\Mod_R$, the two notions agree (see \cite[\S~1.11]{KELLY_BASIC}). As one would expect, if $L$ is a left adjoint, then it is cocontinuous (see \cite[\S~3.2]{KELLY_BASIC}). The category of cocontinuous $\ca{V}$-functors $\ca{A} \rightarrow \ca{B}$ will be denoted by $\Cocts[\ca{A},\ca{B}]$.

 Let $V\in \ca{V}$. If the $\ca{V}$-functor $[V,\ca{E}(E,-)] \colon \ca{E} \rightarrow \ca{V}$ is representable, we denote the representing object by $V \odot E$ and we call it the \emph{tensor product} of $V$ and $E$. This concept is a special case of a weighted colimit: for $\ca{D}=\ca{I}$, the unit $\ca{V}$-category with one object $0$ and $\ca{I}(0,0)=I$, giving a weight amounts to giving an object $V \in \ca{V}$, giving a $\ca{V}$-functor $\ca{I} \rightarrow \ca{E}$ amounts to giving an object $E \in \ca{E}$, and the colimit of $E$ weighted by $V$ is precisely the tensor product $V\odot E$. For a subcategory $\ca{X} \subseteq \ca{V}$ we say that $\ca{E}$ is \emph{$\ca{X}$-tensored} if the tensor product $V\odot E$ exists for all $E\in \ca{E}$ and all $V\in \ca{X}$. If $\ca{E}$ is $\ca{X}$-tensored for $\ca{X}=\ca{V}$ we simply say that $\ca{E}$ is \emph{tensored}.

If a $\ca{V}$-category has all small weighted colimits, then the colimit of $G\colon \ca{D} \rightarrow \ca{E}$ weighted by  $J \colon \ca{D}^{\op} \rightarrow \ca{V}$ is given by the coend
\[
J \star G = \int^{D\in \ca{D}} JD \odot GD
\]
 (see \cite[\S~3.10]{KELLY_BASIC}). There are also weights corresponding to ordinary diagrams in the underlying category. To distinguish them from general weights the corresponding colimits are called \emph{conical colimits}. For $\ca{V}=\Mod_R$, a conical colimit exists if and only if the corresponding ordinary colimit exists in the underlying category (see \cite[\S~3.8]{KELLY_BASIC}), so in this case there is no need to distinguish the two notions.

\subsection[CAUCHY_COMPLETION_SECTION]{}\label{CAUCHY_COMPLETION_SECTION}
 We can now define Cauchy objects and Cauchy completions of $\ca{V}$-categories. Let $\ca{B}$ be a small $\ca{V}$-category. An object $F\in \Prs{B}$ is called a \emph{Cauchy object} if the representable functor $\Prs{B}(F,-)$ is cocontinuous. The \emph{Cauchy completion} $\overline{\ca{B}}$ of $\ca{B}$ is the full subcategory of Cauchy objects in $\Prs{B}$. If $F$ is represented by the object $B\in \ca{B}$, then the Yoneda lemma implies that $\Prs{B}(F,-)$ is isomorphic to the functor which evaluates a presheaf in $\ca{B}$. It follows that $\overline{\ca{B}}$ contains all the representable functors, i.e., we have $\ca{B} \subseteq \overline{\ca{B}}$.
 
 Cauchy completions are best explained by giving a few examples. If $\ca{V}=\Mod_R$ is the cosmos of $R$-modules for some commutative ring $R$, then an $R$-algebra $B$ can be considered as a one object $\ca{V}$-category $\ca{B}$, where $\ca{B}(\ast,\ast)=B$. The presheaf category $\Prs{B}$ is isomorphic to the $R$-linear category of right $B$-modules. A $B$-module $M$ is a Cauchy object if and only if it is finitely generated and projective. This is why Cauchy objects in an arbitrary cosmos are sometimes called \emph{small projective} objects. The name `Cauchy completion' comes from a different example due to F.\ W.\ Lawvere. Let $\ca{V}$ be the cosmos $[0,\infty]$ of extended non-negative real numbers. For $x,y$ objects of $[0,\infty]$, there is a unique morphism $x \rightarrow y$ if and only if $x\geq y$, and the tensor product is given by addition of real numbers. A $\ca{V}$-category is a (generalized) metric space. Any ordinary metric space $X$ gives an example of a $[0,\infty]$-category, and the Cauchy completion as a $[0,\infty]$-category coincides with the usual Cauchy completion of $X$ as a metric space (see \cite{LAWVERE_METRIC}). If $\ca{V}=\Set$, then a small $\ca{V}$-category $\ca{B}$ is just a small ordinary category, and $\overline{\ca{B}}$ is the \emph{Karoubi envelope} of $\ca{B}$, which is the universal category containing $\ca{B}$ in which all idempotents split.

 An important example of a Cauchy completion which works in any cosmos is the following. For $\ca{B}=\ca{I}$, the unit $\ca{V}$-category, we have $\Prs{B}\simeq \ca{V}$ and the representable functors $\Prs{B}(X,-)$ correspond to $[X,-]$ under this equivalence. Since $X$ has a dual if and only if the internal hom-functor is cocontinuous we conclude that the Cauchy completion of $\ca{I}$ is equivalent to $\ca{V}^c$, the full subcategory of $\ca{V}$ consisting of objects with duals.

\subsection[DENSE_FUNCTORS_SECTION]{}\label{DENSE_FUNCTORS_SECTION}
 There is another concept from the theory of enriched categories which we are going to use frequently, namely the notion of a \emph{dense} $\ca{V}$-functor. This notion is motivated as follows (see \cite[Chapter~5]{KELLY_BASIC}). A continuous map $f\colon X\rightarrow Y$ between Hausdorff topological spaces has dense image if and only if a continuous map $g\colon Y\rightarrow Z$ into another Hausdorff space is uniquely determined by the composite $gf$. A dense functor has an analogous property, where `continuous map' gets replaced by `cocontinuous functor'. A $\ca{V}$-functor $K \colon \ca{A} \rightarrow \ca{B}$ is \emph{dense} if precomposing with $K$ induces a fully faithful functor $\ca{A}\mbox{-}\Cocts[\ca{B},\ca{C}] \rightarrow [\ca{A},\ca{C}]$ for every $\ca{V}$-category $\ca{C}$, where $\ca{A}\mbox{-}\Cocts[\ca{B},\ca{C}]$ stands for the full subcategory of those $\ca{V}$-functors which preserve those weighted colimits whose weights have domain $\ca{A}^{\op}$. Thus $\ca{A}\mbox{-}\Cocts[\ca{B},\ca{C}]$ contains $\Cocts[\ca{B},\ca{C}]$ if $\ca{A}$ is small. It follows that for any dense $\ca{V}$-functor $K \colon \ca{A} \rightarrow \ca{B}$ with small domain, for any two cocontinuous $\ca{V}$-functors $F,G \colon \ca{B} \rightarrow \ca{C}$ and for any $\ca{V}$-natural transformation $\alpha \colon FK \Rightarrow GK$ there is a unique $\ca{V}$-natural transformation $\beta \colon F \Rightarrow G$ such that for all $A\in \ca{A}$, $\beta_{KA}=\alpha_A$. One of the most important examples of a dense functor is the Yoneda embedding $Y\colon \ca{A} \rightarrow \Prs{A}$ of a small $\ca{V}$-category $\ca{A}$ (see \cite[Proposition~5.16]{KELLY_BASIC}). A $\ca{V}$-functor $K \colon \ca{A} \rightarrow \ca{C}$ is dense if and only if any object $C$ in $\ca{C}$ is a canonical weighted colimit of $K$ (see \cite[Theorem~5.1]{KELLY_BASIC}). This relates the notion of $\ca{V}$-dense functors to the $\Set$-dense subcategories introduced in Definition~\ref{SET_DENSITY_DFN}: a subcategory is $\Set$-dense if and only if the inclusion functor is $\Set$-dense.

\section{Enriched categories of comodules}\label{CATEGORIES_OF_COMODULES_SECTION}
 A comonoid $T$ in $\ca{V}$ gives rise to the comonad $T \otimes - \colon \ca{V} \rightarrow{\ca{V}}$, which in turn gives rise to an adjunction $V_T \colon \ca{V}_T \rightleftarrows \ca{V} \colon W_T$, and we want to use this to give our new construction of the Tannakian adjunction. More precisely, we will use the semantics-structure adjunction, which we define in Section~\ref{SEMANTICS_STRUCTURE_SECTION}. It relates arbitrary comonads on $\ca{V}$ to $\ca{V}$-categories over $\ca{V}$, and, more generally, comonads on a $\ca{V}$-category $\ca{C}$ to $\ca{V}$-categories over $\ca{C}$. In order to define the semantics-structure adjunction we need to introduce enriched categories of comodules of a comonad.

\subsection[COMONAD_DEFINITION_SECTION]{}\label{COMONAD_DEFINITION_SECTION}
 A \emph{comonad} on a $\ca{V}$-category $\ca{C}$ is a comonoid in the category of endo-$\ca{V}$-functors of $\ca{C}$ and $\ca{V}$-natural transformations. In other words, a comonad consists of a $\ca{V}$-functor $T\colon \ca{C}\rightarrow\ca{C}$, a \emph{comultiplication} $\delta\colon T\Rightarrow TT$ and a \emph{counit} $\varepsilon \colon T\Rightarrow \id$ such that the diagrams
\[
\vcenter{ \xymatrix{
T \ar@{=>}[r]^{\delta}\ar@{=>}[d]_{\delta} & TT \ar@{=>}[d]^{\delta T}\\
TT \ar@{=>}[r]_{T\delta} & TTT\\
}}\quad\textrm{and}\quad
\vcenter{\xymatrix{
& T \ar@{=>}[d]^{\delta}\\
T \ar@{=}[ru] & \ar@{=>}[l]^{\varepsilon T} TT \ar@{=>}[r]_{T \varepsilon} &
\ar@{=}[lu] T \\
}} 
\]
 are commutative. A morphism of comonads on $\ca{C}$ is a $\ca{V}$-natural transformation which is compatible with the counit and the comultiplication.

\subsection[COMODULE_DEFINITION_SECTION]{}\label{COMODULE_DEFINITION_SECTION}
 For any comonad $T$ on $\ca{C}$, there is an associated $\ca{V}$-category $\ca{C}_T$ with an adjunction $\eta^T, \varepsilon^T \colon V_T \colon \ca{C}_T \rightleftarrows \ca{C} \colon W_T$. The objects of the $\ca{V}$-category $\ca{C}_T$ are the \emph{comodules}\footnote{Comodules for a comonad are usually called coalgebras, but in our context this terminology would be misleading; if $\ca{V}$ is the category of modules for some commutative ring, the comonoids giving rise to the relevant comonads are themselves called coalgebras.} of $T$, i.e., of objects $M$ of $\ca{C}$ equipped with a morphism $\varrho \colon M \rightarrow TM$ such that the diagrams
\[
\vcenter{ \xymatrix{
M \ar[r]^{\varrho}\ar[d]_{\varrho} & TM \ar[d]^{\delta_M}\\
TM \ar[r]_{T\varrho} & TTM\\
}}\quad\textrm{and}\quad
\vcenter{\xymatrix{
M \ar[r]^{\varrho} \ar@{=}[rd] & TM \ar[d]^{\varepsilon_M}\\
& M
}}
\]
 are commutative. The hom-object $\ca{C}_T(\mathbf{M},\mathbf{M^{\prime}})$ between the comodules $\mathbf{M}=(M,\varrho)$ and $\mathbf{M^{\prime}} = (M^{\prime},\varrho^{\prime})$ is given by the equalizer
\[
\xymatrix{
\ca{C}_T(\mathbf{M},\mathbf{M^{\prime}}) \ar[rr]^-{V_T} & & 
\ca{C}(M,M^{\prime})
\ar@<0.5ex>[rr]^{\ca{C}(\varrho,TM^{\prime}) \circ T}
\ar@<-0.5ex>[rr]_{\ca{C}(M,\varrho^{\prime})} & & 
\ca{C}(M,TM^{\prime})
}
\]
 in $\ca{V}$, and the structure map of this equalizer gives the functor $V_T\colon \ca{C}_T \rightarrow \ca{C}$. The right adjoint $W_T\colon \ca{C} \rightarrow \ca{C}_T$ sends an object $M$ of $\ca{C}$ to $(TM,\delta_M)$, called the \emph{cofree comodule} on $M$, and $(W_T)_{M,M^{\prime}}\colon \ca{C}(M,M^{\prime}) \rightarrow \ca{C}_T\bigl((TM,\delta_M), (TM^{\prime}, \delta_{M^{\prime}}) \bigr)$ is the induced map
\[
\xymatrix@!C=35pt{
\ca{C}_T\bigl( (TM,\delta_M),(TM^{\prime},\delta_{M^{\prime}})\bigr)
\ar[rr]^-{V_T} & & 
\ca{C}(TM,TM^{\prime})
\ar@<0.5ex>[rr]^{\ca{C}(\delta_M, T^2 M^{\prime}) \circ T}
\ar@<-0.5ex>[rr]_{\ca{C}(TM, \delta_{M^{\prime}})} & & 
\ca{C}(TM,T^2 M^{\prime}) \\
\ca{C}(M,M^{\prime}) \ar[u]^{W_T} \ar[rru]_{T}
}
\]
 into the corresponding equalizer. The unit $\eta^T$ and counit $\varepsilon^T$ of this adjunction are given by $\xymatrix@1{(M,\varrho) \ar[r]^-{\varrho} & (TM,\delta_M)=W_T V_T (M,\varrho)}$ and $\xymatrix@1{V_T W_T (M)=TM \ar[r]^-{\varepsilon} & M}$ respectively.

 The next proposition shows that weighted colimits (see Section~\ref{WEIGHTED_COLIMIT_SECTION}) in the category $\ca{C}_T$ of $T$-comodules are computed the same way as in $\ca{C}$.

\begin{prop}\label{COLIMIT_CREATION_PROP}
 Let $T$ be a comonad on $\ca{C}$. The functor $V_T \colon \ca{C}_T \rightarrow \ca{C}$ creates weighted colimits, i.e., if $K \colon \ca{A} \rightarrow \ca{C}_T$ is a $\ca{V}$-functor and $J \colon \ca{A}^{\op} \rightarrow \ca{V}$ is a weight such that the colimit $J \star V_T K$ exists, then $J\star K$ exists and is preserved by $V_T$.
\end{prop}

\begin{proof}
 If $\ca{C}$ is cocomplete, then $\ca{C}_T$ is cocomplete by the dual of \cite[Proposition~III.2.6]{DUBUC}. The functor $V_T$ reflects isomorphisms and has a right adjoint, so by the dual of the remark near the beginning of \cite[\S~3.6]{KELLY_BASIC} it creates all small weighted colimits. One can prove the general result using the concept of a \emph{$T$-coaction}, defined in Section~\ref{COACTION_DEFINITION_SECTION}. We will not need the more general result because we are only interested in the case where $\ca{C}$ is a category of enriched presheaves.
\end{proof}

\section{The comodule functor}\label{COMODULE_FUNCTOR_SECTION}
\subsection[COMODULE_FUNCTOR_OVERVIEW_SECTION]{}\label{COMODULE_FUNCTOR_OVERVIEW_SECTION}
 The goal of this section is to give a construction of the right adjoint of the Tannakian adjunction, which sends comonoids to categories of representations. For $\ca{V}$-categories $\ca{A}$, $\ca{B}$, a \emph{module} $M$ from $\ca{A}$ to $\ca{B}$ (also called \emph{profunctor} or \emph{distributor}) consists of an object $M(b,a) \in \ca{V}$ for each $b \in \ca{B}$, $a\in \ca{A}$, together with maps
\[
M(b,a) \otimes \ca{B}(b^\prime,b) \rightarrow M( b^{\prime},a) \quad \mathrm{and} \quad \ca{A}(a,a^\prime) \otimes M(b,a) \rightarrow M( b,a^{\prime})
\]
 subject to certain natural coherence conditions. There are many equivalent ways to define the category of modules from $\ca{A}$ to $\ca{B}$. Since $\ca{V}$ is closed, a module is the same as a $\ca{V}$-functor $\ca{B}^{\op}\otimes \ca{A} \rightarrow \ca{V}$. By adjunction, this is the same as a $\ca{V}$-functor $\ca{A} \rightarrow \Prs{B}$. Since $\Prs{A}$ is the free cocompletion of $\ca{A}$ (cf.\ \cite[Theorem~4.51]{KELLY_BASIC} or Theorem~\ref{FREE_COCOMPLETION_THM}), the category of modules is also equivalent to the category $\Cocts[\Prs{A},\Prs{B}]$ of cocontinuous functors $\Prs{A} \rightarrow \Prs{B}$. We denote the category of modules from $\ca{A}$ to $\ca{B}$ by $\modules{A}{B}$. Composition of functors in $\Cocts[\Prs{B},\Prs{B}]$ induces a monoidal structure on $\Bimod{B} \simeq \Cocts[\Prs{B},\Prs{B}]$. The \emph{comodule functor} (see Definition~\ref{COMODULE_FUNCTOR_DFN}) is a functor from the category $\CC{B}$ of comonoids in $\Bimod{B}$ to the category $\Vcat\slash\overline{\ca{B}}$ of small $\ca{V}$-categories over the Cauchy completion $\overline{\ca{B}}$ of $\ca{B}$. An important special case is the \emph{neutral} case, where $\ca{B}=\ca{I}$, and the comodule functor sends comonoids in $\ca{V}\simeq\Bimod{I}$ to categories over $\overline{\ca{I}}\simeq \ca{V}^c$, the full subcategory of $\ca{V}$ consisting of objects with duals (see Section~\ref{CAUCHY_COMPLETION_SECTION}). In this case the above mentioned equivalence $\ca{V}\simeq \Bimod{I}\simeq \Cocts[\ca{V},\ca{V}]$ sends an object $M \in \ca{V}$ to the cocontinuous functor $M\otimes- \colon \ca{V} \rightarrow \ca{V}$.

 From now on, we usually work with $\Cocts[\Prs{B},\Prs{B}]$ instead of $\Bimod{B}$, because the former is a strict monoidal category. Using the equivalence from the previous paragraph one can easily transfer results about one of these categories to the other; we will not do this explicitly.

\subsection[CAUCHY_COMODULE_SECTION]{}\label{CAUCHY_COMODULE_SECTION}
 If $T \colon \Prs{B} \rightarrow \Prs{B}$ is a comonad on the category $\Prs{B}$ of presheaves on $\ca{B}$, then we say that a comodule $(M,\varrho)$ is a \emph{Cauchy comodule} if its underlying object $M$ lies in the Cauchy completion of $\ca{B}$. We denote the category of Cauchy comodules of $T\colon \Prs{B} \rightarrow \Prs{B}$ by $\Prs{B}_T^c$. This category constitutes the object part of the right adjoint of the Tannakian adjunction. In order to show functoriality, we need a nice description of functors with codomain $\Prs{B}_T$, which requires the introduction of the following concept.

\subsection[COACTION_DEFINITION_SECTION]{}\label{COACTION_DEFINITION_SECTION}
 A \emph{coaction} of a comonad $T\colon \ca{C}\rightarrow \ca{C}$ on a functor $S\colon\ca{A}\rightarrow \ca{C}$ is a natural transformation $\varrho\colon S\Rightarrow TS$ such that
\[
\vcenter{ \xymatrix{
S \ar@{=>}[r]^{\varrho}\ar@{=>}[d]_{\varrho} & TS \ar@{=>}[d]^{\delta S}\\
TS \ar@{=>}[r]_{T\varrho} & TTS\\
}}\quad\textrm{and}\quad
\vcenter{\xymatrix{
S \ar@{=>}[r]^{\varrho} \ar@{=}[rd] & TS \ar@{=>}[d]^{\varepsilon S}\\
& S
}}
\]
 are commutative. Note that if we identify functors $M \colon \ca{I} \rightarrow \ca{C}$ from the unit $\ca{V}$-category to $\ca{C}$ with objects in $\ca{C}$, then a coaction on $M$ is precisely a morphism $\varrho \colon M\rightarrow TM$ which turns $M$ into a comodule.

\begin{prop}\label{COACTION_LIFT_PROP}
 Let $\varrho\colon S\Rightarrow TS$ be a coaction. Then the assignment $\overline{S}(A)=(SA,\varrho_A)$ extends uniquely to a $\ca{V}$-functor $\overline{S} \colon \ca{A} \rightarrow \ca{C}_T$ such that
\[
\xymatrix{ \ca{A} \ar[rd]_S \ar[rr]^{\overline{S}} & & \ar[ld]^{V_T} \ca{C}_T\\
& \ca{C}\\
}
\]
 is commutative. This gives a bijection, natural in $T$, between coactions on $S$ and functors $\overline{S}$ which make the above triangle commutative. If $S^{\prime} \colon \ca{A} \rightarrow \ca{C}_T$ is another $\ca{V}$-functor, together with a $T$-coaction $\varrho^{\prime} \colon S^{\prime} \rightarrow TS^{\prime}$, then whiskering with $V_T \colon \ca{C}_T \rightarrow \ca{C}$ (see Appendix~\ref{PASTING_COMPOSITES_APPENDIX}) gives a bijection between $\ca{V}$-natural transformations $\overline{\varphi} \colon \overline{S} \Rightarrow \overline{S^{\prime}}$ and $\ca{V}$-natural transformations $\varphi \colon S \Rightarrow S^{\prime}$ which make the diagram
\[
\xymatrix{
S \ar@{=>}[r]^{\varphi} \ar@{=>}[d]_{\varrho} & S^{\prime}
\ar@{=>}[d]^{\varrho^{\prime}} \\
TS \ar@{=>}[r]_{T\varphi} & TS^{\prime}
}
\]
commutative.
\end{prop}

\begin{proof}
 The first statement is dual to \cite[Proposition~II.1.1]{DUBUC}, and the statement about $\ca{V}$-natural transformations between lifts follows from the fact that for any two $T$-comodules $\mathbf{M}$ and $\mathbf{M^{\prime}}$, the component $(V_T)_{\mathbf{M}, \mathbf{M^{\prime}}} \colon \ca{C}_T (\mathbf{M}, \mathbf{M^{\prime}}) \rightarrow \ca{C}(V_T \mathbf{M}, V_T \mathbf{M^{\prime}})$ of $V_T$ is by definition an equalizer (cf.\ \cite[p.~64]{DUBUC}).
\end{proof}

The definition of monads and actions of monads makes sense in any strict 2-category. This is known as `the formal theory of monads', see \cite{STREET_FTM}. With the terminology introduced there, the above result says that $\VCAT^{\co}$, the strict 2-category with the same objects and 1-cells as $\VCAT$ and with reversed 2-cells, admits the construction of algebras.

\subsection[FUNCTORIALITY_OF_COMODULES_SECTION]{}\label{FUNCTORIALITY_OF_COMODULES_SECTION}
 Proposition~\ref{COACTION_LIFT_PROP} enables us to show that $\ca{C}_T$ and $\Prs{B}_T^c$ are functorial in the comonad $T$. A morphism of comonads $\varphi \colon T \Rightarrow T^{\prime}$ induces a $T^{\prime}$-coaction on $V_{T}$, with $(M,\varrho)$-component given by
\[
\xymatrix{V_T(M,\varrho)=M \ar[r]^-{\varrho} & TM \ar[r]^-{\varphi_M} &
T^{\prime}M=T^{\prime}V_T(M,\varrho)\rlap{,}}
\]
 and by Proposition~\ref{COACTION_LIFT_PROP} this corresponds to a functor $\ca{C}_\varphi \colon \ca{C}_T \rightarrow \ca{C}_{T^{\prime}}$ making the diagram
\[
\xymatrix{\ca{C}_T \ar[rr]^-{\ca{C}_\varphi} \ar[rd]_{V_T} & &
\ca{C}_{T^{\prime}} \ar[ld]^{V_{T^{\prime}}}\\ & 
\ca{C}}
\]
 commutative. Functoriality is an immediate consequence of the uniqueness condition in Proposition~\ref{COACTION_LIFT_PROP}, and in the special case where $T$ is a comonad on $\ca{C}=\Prs{B}$, commutativity of the above diagram implies that $\Prs{B}_\varphi$ sends Cauchy comodules to Cauchy comodules.

\begin{dfn}\label{COMODULE_FUNCTOR_DFN}
 Let $\ca{B}$ be a $\ca{V}$-category whose Cauchy completion $\overline{\ca{B}}$ is small. Let $T$ be a cocontinuous comonad on $\Prs{B}$, which we identify with the corresponding module from $\ca{B}$ to $\ca{B}$ under the equivalence from Section~\ref{COMODULE_FUNCTOR_OVERVIEW_SECTION}. Let $\Prs{B}^c_T$ be the category of Cauchy comodules of $T$ (see Section~\ref{CAUCHY_COMODULE_SECTION}). We let $V^c_T \colon \Prs{B}^c_T \rightarrow \overline{\ca{B}}$ be the restriction of the forgetful functor $V_T\colon \Prs{B}_T \rightarrow \Prs{B}$ to the full subcategory $\Prs{B}^c_T$. The \emph{comodule functor} is the functor
\[
\Prs{B}_{(-)}^c \colon \CC{B} \longrightarrow \Vcat / \overline{\ca{B}}
\]
 from the category of cocontinuous comonads on $\Prs{B}$ to the category of small $\ca{V}$-categories over $\overline{\ca{B}}$ which sends a cocontinuous comonad $T$ to the pair $(\Prs{B}^c_T,V^c_T)$.
\end{dfn}

 Note that $\Prs{B}_T^c$ really is a small category: for any $M \in \overline{\ca{B}}$, there is only a set of morphisms $M\rightarrow TM$ which turn $M$ into a $T$-comodule. This concludes the construction of the right adjoint of the Tannakian adjunction.

\section{The semantics-structure adjunction}\label{SEMANTICS_STRUCTURE_SECTION}
\subsection[SEMANTICS_STRUCTURE_EXPLANATION_SECTION]{}\label{SEMANTICS_STRUCTURE_EXPLANATION_SECTION}
 The functor which sends a comonad to its category of comodules is traditionally known as the \emph{semantics} functor. It has a partial left adjoint, which is called the \emph{structure} functor (cf.\ \cite{DUBUC}, \cite{STREET_FTM}). The \emph{partial} here means that the left adjoint of the semantics functor is only defined on a certain subcategory. These names go back to Lawvere (see \cite[p.~77]{LAWVERE_SEMANTICS}). A monad can be viewed as a sort of logical theory, and from this viewpoint the semantics functor sends it to its category of models; the study of the models of a logical theory is generally called its semantics. On the other hand, a monad can also be seen as a `type of structure' with which objects of $\ca{C}$ can be equipped. From this point of view, the structure functor sends a category $\ca{A}$ over $\ca{C}$ to the universal structure with which the objects of $\ca{A}$ can be equipped. We give a precise statement of this result in Theorem~\ref{SEMANTICS_STRUCTURE_THM}, which also contains explicit descriptions of the natural bijections of the adjunction. Theorem~\ref{COMONADICITY_THM} characterizes those categories over $\ca{C}$ for which the unit of the semantics-structure adjunction is an equivalence. This result is an enriched version of the dual of Beck's famous monadicity theorem.

\subsection[KAN_EXTENSION_SECTION]{}\label{KAN_EXTENSION_SECTION}
 Given $\ca{V}$-functors $F\colon \ca{A} \rightarrow \ca{B}$, $K \colon \ca{A} \rightarrow \ca{C}$ and $L \colon \ca{B} \rightarrow \ca{C}$, and a $\ca{V}$-natural transformation
\[
\xymatrix{ \ca{A} \ar[r]^-{F} \ar[rd]_{K} \xtwocell[1,1]{}\omit{^<-2>\pi} & \ca{B} \ar[d]^{L}\\& \ca{C}}
\]
we say that $\pi$ exhibits $L$ as \emph{left extension\footnote{Left extensions are called \emph{weak left Kan extensions} in \cite{KELLY_BASIC}.} of $K$ along $F$} if for any $\ca{V}$-functor $S \colon \ca{B} \rightarrow \ca{C}$, the assignment which sends a $\ca{V}$-natural transformation $\varphi \colon L \Rightarrow S$ to the pasted composite
\[
\xymatrix@!C=25pt@R=30pt{ \ca{A} \ar[r]^-{F} \ar[rd]_{K} \xtwocell[1,1]{}\omit{^<-2.5>\pi} & \ca{B} \dtwocell_{L}^{S}{^\varphi} \\ & \ca{C}} 
\]
gives a bijection $\VNat(L,S) \rightarrow \VNat(K,SF)$. The $\ca{V}$-natural transformation $\pi$ is called the \emph{unit} of the left extension. A $\ca{V}$-functor $G \colon \ca{C} \rightarrow \ca{C}^{\prime}$ \emph{preserves} the left extension $L$ if $G\pi$ exhibits $GL$ as left extension of $GK$ along $F$. A left extension $L\colon \ca{B} \rightarrow \ca{C}$ is called a \emph{left Kan extension} if the $\ca{V}$-functors $\ca{C}(-,C) \colon \ca{C} \rightarrow \ca{V}^{\op}$ preserve the left extension $L$ for all $C \in \ca{C}$. The left Kan extension of $K$ along $F$ is usually denoted by $\Lan_F K \colon \ca{B} \rightarrow \ca{C}$.

For example, if $F \colon \ca{A} \rightarrow \ca{B}$ has a right adjoint $R \colon \ca{B} \rightarrow \ca{A}$, with unit $\eta \colon \id \Rightarrow RF$ and counit $\varepsilon \colon FR \Rightarrow \id$, then $K\eta \colon K \Rightarrow KR \circ F$ exhibits $KR$ as left extension of $K$ along $F$. Indeed, the triangular identities imply that the assignment $\VNat(K,SK) \rightarrow \VNat(KR,S)$, $\psi \mapsto  S\varepsilon \circ \psi R$ gives the desired inverse. Moreover, this left extension is preserved by \emph{any} $\ca{V}$-functor $G \colon \ca{C} \rightarrow \ca{C}^\prime$, for we can replace $K$ by $GK$ in the above argument. It follows in particular that $KR$ is the left Kan extension of $K$ along $F$.

In \cite[Proposition~4.33]{KELLY_BASIC} it was shown that the left Kan extension $\Lan_F K$ of $K \colon \ca{A} \rightarrow \ca{C}$ along $F \colon \ca{A} \rightarrow \ca{B}$ always exists if $\ca{A}$ is small and $\ca{C}$ is cocomplete.

\subsection[DENSITY_COMONAD_SECTION]{}\label{DENSITY_COMONAD_SECTION}
 Let $\omega \colon \ca{A} \rightarrow \ca{C}$ be a $\ca{V}$-functor such that the left extension of $\omega$ along itself exists (see Section~\ref{KAN_EXTENSION_SECTION}). We denote this left extension by $L(\omega) \colon \ca{C} \rightarrow \ca{C}$, and we let $\pi \colon \omega \Rightarrow L(\omega) \circ \omega$ be its unit. By definition of left extensions, there are unique $\ca{V}$-natural transformations $\varepsilon \colon L(\omega) \Rightarrow \id$ and $\delta \colon L(\omega) \Rightarrow L(\omega) \circ L(\omega)$ such that the equalities
\[
\vcenter{
\xymatrix@C=25pt@R=20pt{
\xtwocell[0,1]{+<13pt,0pt>{}}\omit{<1.75>*!<12pt,-2pt>{\pi}}
\ar@/^/[rr]^{\omega} \ar@/_/[rd]_{\omega}  && \\
 & \urtwocell^{*!<4pt,0pt>{L(\omega)}}_{\id}
{*!<-2pt,-2pt>{\varepsilon}}&
}}=\xymatrix@1{\rrtwocell^{\omega}_{\omega}{*!<-2pt,0pt>{\id}} & & }
\quad\text{and}\quad
\vcenter{
 \xymatrix@C=25pt@R=35pt{ \rrtwocell \omit{<2>*!<12pt,0pt>{\pi}}
\ar@/^/[rr]^{\omega} \ar@/_/[rd]_{\omega}  && \\
 & \uruppertwocell<4>^{*!<11pt,11pt>{L(\omega)}}{<0.5>\delta}
\urcompositemap<-3>_{*!<-4pt,1pt>{L(\omega)}}^{*!<-4pt,2pt>{L(\omega)}}{\omit}
}}
=
\vcenter{
 \xymatrix@C=25pt@R=18pt{ \rrtwocell \omit{<1>\pi} \ar@/^/[rr]^{\omega}
\ar[rd]^{\omega} \ar@/_/[dd]_{\omega} && \\
\xtwocell[0,1]{+<-5pt,0pt>} \omit{*!<12pt,0pt>{\pi}} & \ar[ru]_{L(\omega)} \\
\ar[ru]_{L(\omega)}
}}
\]
 hold. It follows easily from the definition of left extensions that $\bigl(L(\omega),\delta,\varepsilon \bigr)$ is a comonad on $\ca{C}$. This comonad is called the \emph{density comonad} of $\omega$. The name comes from the fact that if $\omega$ is dense, then the identity natural transformation $\omega \Rightarrow  \id \circ \omega$ exhibits $\id \colon \ca{C} \rightarrow \ca{C}$ as left Kan extension of $\omega$ along itself (see \cite[Theorem~5.1]{KELLY_BASIC}), so the density comonad measures the failure of a $\ca{V}$-functor to be dense. From the example of left extension along left adjoints given in Section~\ref{KAN_EXTENSION_SECTION} it follows that for an adjunction $\eta, \varepsilon \colon F \colon \ca{E} \rightleftarrows \ca{C} \colon G$, the density comonad of $F$ is $(FG,F\eta G,\varepsilon)$. In particular, the density comonad of the functor $V_T \colon \ca{C}_T \rightarrow \ca{C}$ is the comonad $T$ itself (see Section~\ref{COMODULE_DEFINITION_SECTION}).

 Note that there is a natural coaction of the density comonad $L(\omega)$ on $\omega$, given by the unit $\pi \colon \omega \Rightarrow L(\omega) \circ \omega$ of the left extension $L(\omega)$. The fact that $\pi$ is a coaction is an immediate consequence of the definition of the counit and comultiplication of the density comonad. It is in this sense that the objects of $\ca{A}$ can be equipped with an $L(\omega)$-structure, which justifies the name `structure functor' for the partial functor which sends $\omega$ to the density comonad $L(\omega)$.

\begin{thm}[Semantics-structure adjunction]\label{SEMANTICS_STRUCTURE_THM}
 Let $\omega \colon \ca{A} \rightarrow \ca{C}$ be a $\ca{V}$-functor, let $\pi \colon \omega \Rightarrow L(\omega) \circ \omega$ exhibit $L(\omega) \colon \ca{C} \rightarrow \ca{C}$ as left extension of $\omega$ along itself, and let $T\colon \ca{C}\rightarrow \ca{C}$ be a comonad. The assignment
\[
\vcenter{ \xymatrix@C=30pt{ \ca{C} \rtwocell^{L(\omega)}_{T}{\varphi} & \ca{C}} } \quad\mapsto\quad
\vcenter{
\xymatrix@C=20pt@R=18pt{ \ca{A}
\xtwocell[0,1]{+<13pt,0pt>{}}\omit{<2>*!<12pt,-2pt>{\pi}}
\ar@/^/[rr]^{\omega} \ar@/_/[rd]_{\omega}  && \ca{C} \\
 & \ca{C} \urtwocell^{*!<4pt,-1pt>{L(\omega)}}_{T}
{{\varphi}}&
}}
\]
 induces a bijection, natural in $T$, between morphisms of comonads $\varphi \colon L(\omega) \Rightarrow T$ and $T$-coactions $\varrho \colon \omega \Rightarrow T\omega$. Together with Proposition~\ref{COACTION_LIFT_PROP} this gives a bijection between $\ca{V}$-functors $\overline{\omega} \colon \ca{A} \rightarrow \ca{C}_T$ which make the diagram
\[
\xymatrix{\ca{A} \ar[rd]_\omega \ar[rr]^{\overline{\omega}} & & \ca{C}_T
\ar[ld]^{V_T}\\
& \ca{C} }
\]
 commutative, and morphisms of comonads $\varphi \colon L(\omega) \Rightarrow T$. In other words: the assignment which sends a $\ca{V}$-functor $\omega \colon \ca{A} \rightarrow \ca{C}$ to its density comonad $L(\omega)$ extends to a partial left adjoint
\[
\xymatrix{\VCAT\slash \ca{C} \ar@{-->}[rr]^-{L(-)} && \Comon([\ca{C},\ca{C}])}
\]
 of the functor which sends a comonad to the $\ca{V}$-functor $V_T \colon \ca{C}_T \rightarrow \ca{C}$. The domain of this partial left adjoint consists of $\ca{V}$-functors $\ca{A} \rightarrow \ca{C}$ which admit a density comonad. The $T$-component of its counit is given by $\id \colon L(V_T) \rightarrow T$.
\end{thm}

\begin{proof}
 The first statement is dual to \cite[Proposition~{II.1.4}]{DUBUC}, and the second assertion follows from the first by Proposition~\ref{COACTION_LIFT_PROP}. A more general version of this result is \cite[Theorem~6]{STREET_FTM}. We have $L(V_T)=V_T W_T$ (see Section~\ref{KAN_EXTENSION_SECTION}). If we apply the bijection to $\id_T$, then we get the lift of the coaction $V_T\eta^T$. By definition of $\eta^T$ (see Section~\ref{COMODULE_DEFINITION_SECTION}) this lift is the identity functor on $\ca{C}_T$. Thus the counit of the semantics structure adjunction is indeed given by the identity $\id \colon V_T W_T \rightarrow T$.
\end{proof}

\subsection[COSPLIT_EQUALIZER_SECTION]{}\label{COSPLIT_EQUALIZER_SECTION}
 In order to state the dual of Beck's Monadicity Theorem we have to introduce one more concept. A \emph{cosplit equalizer} in an (ordinary) category $\ca{C}_0$ is a diagram of the form
\[
\turnradius{5pt}
\xymatrix{
E \ar[r]^{s} & A \ar@<0.5ex>[r]^u \ar@<-0.5ex>[r]_v \ar `u[l] `[l]+<0.4pt,8pt>_p [l]+<0pt,5pt> &  B \ar `u[l] `[l]+<2.4pt,8pt>_{q} [l]+<2pt,5pt>\\
}
\]
 such that the equalities $us=vs$, $ps=\id_E$, $qu=\id_A$ and $qv=sp$ hold. These identities imply that $s$ exhibits $E$ as equalizer of $u$ and $v$. If $F \colon \ca{A}_0 \rightarrow \ca{C}_0$ is a functor, then we say that a pair $f,g\colon A \rightarrow A^\prime$ in $\ca{A}_0$ is \emph{$F$-cosplit} if there is a cosplit equalizer in $\ca{C}_0$ as above with $u=Ff$, $v=Fg$.

\begin{thm}\label{COMONADICITY_THM}
 Let $\eta, \varepsilon \colon F \colon \ca{E} \rightleftarrows \ca{C} \colon G$ be an adjunction, and let
\[
 \xymatrix{ \ca{E} \ar[rd]_{F} \ar[rr]^J && \ca{C}_{FG} \ar[ld]^{V_{FG}} \\ & \ca{C}}
\]
 be the unit of the semantics-structure adjunction (see Theorem~\ref{SEMANTICS_STRUCTURE_THM}). Then $J$ is an equivalence if and only if
\begin{enumerate}
\item[i)] 
 The functor $F \colon \ca{E} \rightarrow \ca{C}$ reflects isomorphisms.
\item[ii)]
 The category $\ca{E}$ has conical equalizers of $F_0$-cosplit pairs (see above), and $F$ preserves these conical equalizers.
\end{enumerate}
 In this case, $F$ is called \emph{comonadic}.
\end{thm}

\begin{proof}
Assume that i) and ii) hold. Then $F$ reflects equalizers of $F$-cosplit pairs, i.e., a morphism $f\colon A \rightarrow B$ in $\ca{E}$ is an equalizer of an $F$-cosplit pair $g,h \colon B \rightarrow C$ if and only if $Ff$ is an equalizer of $Fg$ and $Fh$. It follows by the dual of \cite[Theorem~2.II.1]{DUBUC} that $J$ is an equivalence.
Conversely, if $J$ is an equivalence, then i) holds because $V_T \colon \ca{C}_T \rightarrow \ca{C}$ clearly reflects isomorphisms. Condition ii) is again a consequence of \cite[Theorem~2.II.1]{DUBUC}.
\end{proof}

\section{The Tannakian adjunction}\label{TANNAKIAN_ADJUNCTION_SECTION}

\subsection[PARTIAL_ADJUNCTION_SECTION]{}\label{PARTIAL_ADJUNCTION_SECTION}
 Let $\ca{B}$ be a small $\ca{V}$-category whose Cauchy completion $\overline{\ca{B}}$ (see Section~\ref{CAUCHY_COMPLETION_SECTION}) is small. In Theorem~\ref{SEMANTICS_STRUCTURE_THM} we have seen that there is a partial adjunction
\[
 \xymatrix{ *!<33pt,0pt>+{\Comon([\Prs{B},\Prs{B}])} \rtwocell<7>~^{+{}~**\dir{--}}^{L(-)}_{*!<0pt,2pt>+{\Prs{B}_{(-)}}}{`\perp} &*!<-20pt,0pt>+{\VCAT\slash \Prs{B}} }
\]
 where $L(\omega)$ is defined whenever the left extension of $\omega \colon \ca{A} \rightarrow \Prs{B}$ along itself exists. The inclusion $i\colon \overline{\ca{B}} \rightarrow \Prs{B}$ induces a functor $i\circ - \colon \Vcat\slash \overline{\ca{B}} \rightarrow \VCAT\slash \Prs{B}$. For $F \colon \ca{C} \rightarrow \Prs{B}$, let $(\ca{C}^c,F^c)$ be the pullback
\[
 \xymatrix{\ca{C}^c \ar[r] \ar[d]_{F^c} & \ca{C} \ar[d]^F \\
 \overline{\ca{B}} \ar[r]_-i & \Prs{B}}
\]
 of $F$ along $i$, i.e., $\ca{C}^c$ is the full subcategory of $\ca{C}$ consisting of objects $C\in \ca{C}$ with $FC \in \overline{\ca{B}}$. This construction gives a partial right adjoint
\[
 \xymatrix@C=40pt{ *!<20pt,0pt>+{\VCAT\slash \Prs{B}} \rtwocell<7>~_{+{}~**\dir{--}}^{i \circ -}_{*!<0pt,2pt>+{(-)^c}}{`\perp} &*!<-10pt,0pt>+{\Vcat\slash \overline{\ca{B}}} }
\]
 of $i\circ -$. It is defined at $(\ca{C},F)$ if $\ca{C}^c$ is a small $\ca{V}$-category. For any $\omega \colon \ca{A} \rightarrow \overline{\ca{B}}$ with $\ca{A}$ small, the left extension $L(i\omega)$ of $i\omega \colon \ca{A} \rightarrow \Prs{B}$ along itself exists since $\Prs{B}$ is cocomplete (see \cite[Proposition~4.33]{KELLY_BASIC}). If $T$ is a comonad on $\Prs{B}$, then $\Prs{B}_T^c$ is a small $\ca{V}$-category: for any object $M\in \Prs{B}$, there is only a set of $T$-coactions on $M$. Thus the composite
\[
 \xymatrix{ *!<33pt,0pt>+{\Comon([\Prs{B},\Prs{B}])} \rtwocell<7>^{L(i\circ -)}_{*!<0pt,2pt>+{\Prs{B}^c_{(-)}}}{`\perp} &*!<-10pt,0pt>+{\Vcat\slash \overline{\ca{B}}} }
\]
of the two partial adjunctions is defined everywhere. If we restrict the right adjoint to the full subcategory $\CC{B} \subseteq \Comon([\Prs{B},\Prs{B}])$, we clearly get the comodule functor from Section~\ref{COMODULE_FUNCTOR_SECTION}. In order to show that this functor has a left adjoint, it suffices to prove that the image of $L(i\circ -)$ is contained in the full subcategory $\CC{B}$ of cocontinuous comonads. We prove this in Proposition~\ref{DENSITY_COMONAD_COCONT_PROP}. In order to do this we need some more background on left extensions.

\subsection[FREE_COCOMPLETION_SECTION]{}\label{FREE_COCOMPLETION_SECTION}
An important special case of left extensions are left Kan extensions along the Yoneda embedding $Y \colon \ca{A} \rightarrow \Prs{A}$ for some small $\ca{V}$-category $\ca{A}$. If the $\ca{V}$-category $\ca{C}$ is cocomplete, then $\Lan_Y K$ exists for any $\ca{V}$-functor $K \colon \ca{A} \rightarrow \ca{C}$, and it is naturally isomorphic to $-\star K \colon \Prs{A} \rightarrow \ca{C}$ (see \cite[Formula~4.31]{KELLY_BASIC}). The unit of this left Kan extension is an isomorphism because the Yoneda embedding is fully faithful (cf.\ \cite[Proposition~4.23]{KELLY_BASIC}). Since we frequently use this type of Kan extension, we use the abbreviation $L_K$ for $\Lan_Y K$, and we denote the unit by $\alpha_K \colon K \Rightarrow L_K Y$. The $\ca{V}$-functor $L_K$ has a right adjoint $\widetilde{K} \colon \ca{C} \rightarrow \Prs{A}$ with $\widetilde{K}C=\ca{C}(K-,C)$, by definition of weighted colimits (see Section~\ref{WEIGHTED_COLIMIT_SECTION}). We denote the unit and counit of this adjunction by $\eta^K \colon \id \Rightarrow \widetilde{K} L_K$ and $\varepsilon^K \colon L_K \widetilde{K} \Rightarrow \id$ respectively. We have the following theorem, which shows that $\Prs{A}$ is the free cocompletion of $\ca{A}$.

\begin{thm}\label{FREE_COCOMPLETION_THM}
 Let $Y \colon \ca{A} \rightarrow \Prs{A}$ be the Yoneda embedding of the small $\ca{V}$-category $\ca{A}$, and let $\ca{C}$ be a cocomplete $\ca{V}$-category. Then for any cocontinuous $\ca{V}$-functor $S \colon \Prs{A} \rightarrow \ca{C}$ we have
\[
S\cong L_K \cong -\star K \colon \Prs{A} \rightarrow \ca{C}
\]
 where $K=SY \colon \ca{A} \rightarrow \ca{C}$. The assignment $S \mapsto SY$ is an equivalence of $\ca{V}$-categories
\[
[Y,\ca{C}] \colon \Cocts[\Prs{A},\ca{C}] \rightarrow [\ca{A},\ca{C}]\rlap{.}
\]
 The inverse to this equivalence sends $K$ to (a choice of) $\Lan_Y K$.
\end{thm}

\begin{proof}
 This is (part of) \cite[Theorem~4.51]{KELLY_BASIC}.
\end{proof}

\subsection[KAN_EXTENSION_II_SECTION]{}\label{KAN_EXTENSION_II_SECTION}
If $\pi \colon K \Rightarrow LF$ exhibits $L$ as left extension of $K$ along $F$ and $\pi^\prime \colon L \rightarrow L^\prime F^\prime$ exhibits $L^\prime$ as left extension of $L$ along $F^\prime$, then
\[
 \xymatrix@!=30pt{\ca{A} \ar[r]^-{F} \ar[rrd]_{K} \xtwocell[1,2]{}\omit{^<-2>\pi} & \ca{B}  \ar[r]^-{F^\prime} \ar[rd]_(0.4){L} \xtwocell[1,1]{}\omit{^<-3>\pi^\prime} & \ca{B}^\prime \ar[d]^{L^\prime}\\ & & \ca{C}}
\]
exhibits $L^\prime$ as left extension of $K$ along $F^\prime F$, simply because a composite of two bijections is a bijection. Moreover, it is clear that if $G \colon \ca{C} \rightarrow \ca{C}^\prime$ preserves the left extensions $L$ of $K$ along $F$ and $L^\prime$ of $L$ along $F^\prime$, then $G$ preserves the left extension $L^\prime$ of $K$ along $F^\prime F$. In particular, if $L=\Lan_F K$ and $L^\prime =\Lan_{F^\prime} L$, then $\Lan_{F^\prime F} K$ exists and is equal to $\Lan_{F^\prime} \Lan_F K$. If $\lambda\colon F \Rightarrow G$ is a $\ca{V}$-natural isomorphism, then $L\lambda \circ \pi$ clearly exhibits $L$ as left Kan extension of $K$ along $G$. If we combine these two observations with the two examples of left Kan extensions we have, we get the following proposition.

\begin{prop}\label{KAN_EXTENSION_PROP}
 Let $\ca{C}$ be a cocomplete $\ca{V}$-category, $\ca{A}$ a small $\ca{V}$-category, and let $K\colon \ca{A} \rightarrow \ca{C}$, $F \colon \ca{A} \rightarrow \ca{B}$ be $\ca{V}$-functors. Then
\[
 \xymatrix@C=30pt@R=8pt{\ca{A} \xuppertwocell[0,2]{}<3>^{F}{^<1>*!<5pt,0pt>{\alpha_F^{-1}}}
\xtwocell[0,2]{}\omit{\omit *!<-7.5pt,5.5pt>{\cong}} \ar[rd]_-{Y} \xlowertwocell[3,2]{}<-8>_{K}{^<1.5>*!<3pt,-1pt>{\alpha_K}}
\xtwocell[3,2]{}\omit{\omit *!<-3pt,11pt>{\cong}} &   & \ca{B} \ar[ddd]^{L_K \widetilde{F}}\\& \Prs{A}  \ar[ru]^-{L_F} \xtwocell[2,1]{}\omit{^<-3>*!<3pt,-3pt>{L_K \eta^F}} \ar[rdd]_{L_K}\\\\ & & \ca{C}}
\]
exhibits $L_K \widetilde{F}$ as left Kan extension of $K$ along $F$.
\end{prop}

\begin{proof}
 We know that $\alpha_K$ exhibits $L_K$ as $\Lan_Y K$ (see Section~\ref{FREE_COCOMPLETION_SECTION}) and that $L_K \eta^F$ exhibits $L_K \widetilde{F}$ as $\Lan_{L_F} L_K$ (see Section~\ref{KAN_EXTENSION_SECTION}). Since we can compose left Kan extensions as described in Section~\ref{KAN_EXTENSION_II_SECTION}, we find that $L_K \eta^F Y \circ \alpha_K$ exhibits $L_K \widetilde{F}$ as $\Lan_{L_F Y} K$. The isomorphism $\alpha_F^{-1} \colon L_F Y \Rightarrow F$ shows that the pasting composite in question exhibits $L_K \widetilde{F}$ as $\Lan_F K$.
\end{proof}

\begin{prop}\label{DENSITY_COMONAD_COCONT_PROP}
 Let $\ca{A}$ be a small $\ca{V}$-category, and let $\omega\colon \ca{A} \rightarrow \Prs{B}$ be a $\ca{V}$-functor whose image is contained in $\overline{\ca{B}}$. Then the right adjoint $\widetilde{\omega} \colon \Prs{B} \rightarrow \Prs{A}$ of $L_\omega$ (see Section~\ref{FREE_COCOMPLETION_SECTION}) is cocontinuous, and the density comonad $L(\omega)$ of $\omega$ (see Section~\ref{DENSITY_COMONAD_SECTION}) is cocontinuous.
\end{prop}

\begin{proof}
 Since colimits in $\Prs{A}=[\ca{A}^{\op},\ca{V}]$ are computed pointwise (cf.\ \cite[Section~3.3]{KELLY_BASIC}) it suffices to check that for every $A\in \ca{A}$,
\[
\widetilde\omega(-)(A)=[\ca{B}^{\op},\ca{V}]\bigl(\omega(A),-\bigr)\colon [\ca{B}^{\op},
\ca {V}] \longrightarrow \ca{V}
\]
 preserves weighted colimits. This follows immediately from the fact that the objects of $\overline{\ca{B}}$ are the small projectives in $[\ca{B}^{\op},\ca{V}]$ (cf.\ Section~\ref{CAUCHY_COMPLETION_SECTION} and \cite[Section~5.5]{KELLY_BASIC}). The density comonad $L(\omega)$ of $\omega$ is by Proposition~\ref{KAN_EXTENSION_PROP} equal to the composite $L_\omega \widetilde{\omega}$, hence it is cocontinuous as composite of cocontinuous $\ca{V}$-functors.
\end{proof}

\begin{dfn}[Tannakian adjunction]\label{TANNAKIAN_ADJUNCTION_DFN}
 Let $\ca{B}$ be a $\ca{V}$-category whose Cauchy completion $\overline{\ca{B}}$ is small. The \emph{Tannakian adjunction} is the restriction of the adjunction
\[
 \xymatrix{ *!<33pt,0pt>+{\Comon([\Prs{B},\Prs{B}])} \rtwocell<7>^{L(i\circ -)}_{*!<0pt,2pt>+{\Prs{B}^c_{(-)}}}{`\perp} &*!<-10pt,0pt>+{\Vcat\slash \overline{\ca{B}}} }
\]
 from Section~\ref{PARTIAL_ADJUNCTION_SECTION} to the full subcategory $\CC{B}$ of cocontinuous comonads (cf.\ Proposition~\ref{DENSITY_COMONAD_COCONT_PROP}).

 We use the fully faithful functor $i\circ -$ to identify $\Vcat\slash \overline{\ca{B}}$ with the full subcategory of $\VCAT\slash \Prs{B}$ consisting of $\ca{V}$-functors $\omega \colon \ca{A} \rightarrow \Prs{B}$ with small domain whose image is contained in $\overline{\ca{B}}$. In other words, we omit the inclusion $i \colon \overline{\ca{B}} \rightarrow \Prs{B}$ from the notation. We write
\[
\xymatrix{ \ca{A} \ar[rd]_{\omega} \ar[rr]^-N && \Prs{B}^c_{L(\omega)} \ar[ld]^{V^c_{L(\omega)}} \\
& \Prs{B}} 
\]
for the $(\ca{A},\omega)$-component of the unit of the Tannakian adjunction. For a cocontinuous comonad $T \colon \Prs{B} \rightarrow \Prs{B}$, we denote the $T$-component of the counit by $\nu \colon L(V^c_T) \rightarrow T$.
\end{dfn}

\subsection[TRIVIAL_OBSERVATION_SECTION]{}\label{TRIVIAL_OBSERVATION_SECTION}
 Let $\omega \colon \ca{A} \rightarrow \Prs{B}$ be a functor whose image is contained in $\overline{\ca{B}}$, and let $T \in \CC{B}$. If $F \colon \ca{A} \rightarrow \Prs{B}_T$ is a $\ca{V}$-functor making the diagram
\[
 \xymatrix{ \ca{A} \ar[rd]_{\omega} \ar[rr]^-F && \Prs{B}_T \ar[ld]^{V_T}\\ &\Prs{B}}
\]
 commutative, then the image of $F$ is necessarily contained in $\Prs{B}^c_T$. We therefore get a natural bijection
\[
 \VCAT\slash\Prs{B}\bigl((\ca{A},\omega),(\Prs{B}_T, V_T)\bigr) \cong
\Vcat\slash\overline{\ca{B}}\bigl((\ca{A},\omega),(\Prs{B}^c_T, V^c_T)\bigr)
\]
 whose inverse is given by composing with the inclusion $\Prs{B}_T^c \rightarrow \Prs{B}_T$.

\begin{prop}\label{TANNAKIAN_ADJUNCTION_PROP}
Let $\omega \colon \ca{A} \rightarrow \Prs{B}$ be a $\ca{V}$-functor with small domain whose image is contained in $\overline{\ca{B}}$. With the notation of Section~\ref{KAN_EXTENSION_SECTION}, the density comonad of $\omega$ is given by $L(\omega)=L_\omega \widetilde{\omega}$. Let $T \colon \Prs{B} \rightarrow \Prs{B}$ be a cocontinuous comonad. The natural bijection
\[
 \CC{B}\bigl(L(\omega), T\bigr)\cong \Vcat \slash \overline{\ca{B}}\bigl((\ca{A},\omega),(\Prs{B}^c_T, V^c_T)\bigr)
\]
which exhibits $L(-)$ as left adjoint of the comodule functor sends a morphism of comonads $\varphi \colon L(\omega) \Rightarrow T$ to the factorization $\overline{\omega} \colon \ca{A} \rightarrow \Prs{B}^c_T$ of $\overline{\omega} \colon \ca{A} \rightarrow \Prs{B}_T$ through the inclusion $\Prs{B}^c_T \rightarrow \Prs{B}_T$, where $\overline{\omega}$ is the lift of $\omega$ corresponding to the $T$-coaction
\[
\vcenter{
\vspace{-25pt}
\xymatrix@C=30pt@R=10pt{ \ar[r]^(0.7){Y}
\xuppertwocell[0,4]{}<12>^{\omega} {<-4>*!<-4pt,0pt>{\alpha_\omega}}
\xtwocell[0,4]{}\omit {\omit *!<7.5pt,-21pt>{\cong}}
\xlowertwocell[1,2]{}<-6>_{\omega}{<1>*!<15pt,-4pt>{\alpha_\omega^{-1}}}
\xtwocell[1,2]{}\omit{\omit *!<-6pt,7pt>{\cong}} &
\rruppertwocell<4>^{*!<23pt,2pt>{\id}} {<0.5>*!<-3pt,0pt>{\eta^\omega}} 
\ar[rd]^{L_\omega} & &
\ar[r]^(0.4){L_\omega} & \\
& & \ar[ru]_{\widetilde{\omega}}
\xlowertwocell[-1,2]{}<-6>_{T}{<1>\varphi}}}
\]
on $\omega$ (cf.\ Proposition~\ref{COACTION_LIFT_PROP}).
\end{prop}

\begin{proof}
 We have $L(\omega)=L_\omega \widetilde{\omega}$ by Proposition~\ref{KAN_EXTENSION_PROP}. The described bijection is the composite of the the bijection from Section~\ref{TRIVIAL_OBSERVATION_SECTION} and the bijection
\[
 \CC{B}(L_\omega\widetilde{\omega}, T)\cong \VCAT\slash\Prs{B}\bigl((\ca{A},\omega),(\Prs{B}_T, V_T)\bigr)
\]
from Theorem~\ref{SEMANTICS_STRUCTURE_THM}, where the unit $\pi \colon \omega \Rightarrow L(\omega)$ of $L(\omega)$ is as described in Proposition~\ref{KAN_EXTENSION_PROP}, i.e., $\pi = L_\omega \widetilde{\omega} \alpha_\omega^{-1} \circ L_\omega \eta^\omega Y \circ \alpha_\omega$.
\end{proof}

\section{Reconstruction}\label{RECONSTRUCTION_SECTION}

\subsection[NOTATION_SECTION]{}\label{NOTATION_SECTION}
 The goal of this section is to give a necessary and sufficient condition for the counit of the Tannakian adjunction (see Definition~\ref{TANNAKIAN_ADJUNCTION_DFN}) to be an isomorphism. In order to do this we have to find a suitable description of this counit, hence we fix some notation first. We fix a $\ca{V}$-category $\ca{B}$ with small Cauchy completion $\overline{\ca{B}}$, and a cocontinuous comonad $T \colon \Prs{B} \rightarrow \Prs{B}$. For the sake of brevity we let $\ca{A}=\Prs{B}^c_T$ be the category of Cauchy comodules, we write $K \colon \ca{A} \rightarrow \Prs{B}_T$ for the inclusion functor (see Section~\ref{CAUCHY_COMODULE_SECTION}) and we let $\omega$ be the restriction $\omega = V^c_T \colon \ca{A} \rightarrow \overline{\ca{B}}$ of $V_T$. The functors $K$ and $\omega$ induce adjunctions $\eta^K, \varepsilon^K \colon L_K \colon \Prs{A} \rightleftarrows \Prs{B}_T \colon \widetilde{K}$ and $\eta^\omega, \varepsilon^\omega \colon L_\omega \colon \Prs{A} \rightleftarrows \Prs{B} \colon \widetilde{\omega}$ (see Section~\ref{FREE_COCOMPLETION_SECTION}). The diagram
\[
\xymatrix@C=30pt@R=30pt{
\ca{A} \ar[rd]_-{\omega} \ar[r]^-K & \Prs{B}_T \ar@<-0.5ex>[r]_-{\widetilde{K}} \ar@<0.5ex>[d]^-{V_T} & \Prs{A} \ar@<-0.5ex>[l]_-{L_K} \ar@<1ex>[ld]^{L_\omega}\\
& \Prs{B} \ar@<0.5ex>[u]^{W_T} \ar[ru]^(0.45){\widetilde{\omega}}& 
}
\]
summarizes the situation.

\begin{prop}\label{COUNIT_CHARACTERIZATION_PROP}
 Let $T\in \CC{B}$. With the notation introduced in Section~\ref{NOTATION_SECTION}, there is an isomorphism $\sigma \colon L(\omega) \cong V_T L_K \widetilde{K} W_T$ such that
\[
\nu = \vcenter{
\vspace{-15pt}
\xymatrix@C=30pt@R=10pt{
 \ar[rd]_{W_T}%
\xuppertwocell[0,4]{}<9>^{L(\omega)}{<-1.75>\sigma}%
\xtwocell[0,4]{}\omit{\omit*!<7.5pt,-8pt>{\cong}}%
& & \ar[rd]^{L_K} & &\\%
 & \ar[ru]^{\widetilde{K} }%
\rrlowertwocell<-3>_{\id}{<-0.5>*!<-3pt,0pt>{\varepsilon^K}} & & \ar[ru]_{V_T} & 
}}
\]
is the $T$-component $\nu \colon L(\omega) \Rightarrow T$ of the counit of the Tannakian adjunction.
\end{prop}

\begin{proof}
 By definition of the Tannakian adjunction (see Definition~\ref{TANNAKIAN_ADJUNCTION_DFN}), in order to compute its counit, we have to compute the pasted composite
\[
\vcenter{\xymatrix@C=3pt@R=10pt{%
&& \Vcat\slash \overline{\ca{B}} \ar[rd]^{i\circ -}\\
& \VCAT\slash\Prs{B} \ar[ru]^-{(-)^{c}}\rrlowertwocell<-3>_{\id}{<-1>*!<-3pt,0pt>{\varepsilon_1}}
&&\ar[rd]^-{L(-)} \VCAT\slash\Prs{B} \\
\ar[ru]^-{\Prs{B}_{(-)}} \xlowertwocell[0,4]{}<-7>_{\id}{<0.5>*!<-4pt,0pt>{\varepsilon_2}} %
\CC{B} &&&& \CC{B}}
\vspace{-10pt}}
\]
where we write $\varepsilon_i$ for the counits of the respective adjunctions. The $T$-component of $\varepsilon_2$ is the identity (see Theorem~\ref{SEMANTICS_STRUCTURE_THM}). The $\Prs{B}_T$-component of $\varepsilon_1$ is given by the inclusion $\Prs{B}^c_T \rightarrow \Prs{B}_T$, i.e., by $K \colon \ca{A} \rightarrow \Prs{B}$. It follows that the above pasting composite is given by $L(K)$. Since we only work up to unspecified isomorphism $\sigma$, it suffices to show that there is a natural transformation $\pi \colon \omega \Rightarrow V_T L_K \widetilde{K} W_T \omega $ which exhibits $V_T L_K \widetilde{K} W_T$ as left extension of $\omega$ along itself, and that $L(K)=V_T \varepsilon^T W_T$ if we use this natural transformation to define $L(\omega)=V_T L_K \widetilde{K} W_T$. By definition of $L(K)$ via the universal property of left extensions, this just means that we have to show that $V_T \varepsilon^K W_T \omega \circ\pi = V_T \eta^T K$. Using one of the triangular identities for $\eta^K$ and $\varepsilon^K$ it is not hard to see that
\[
\pi \defl \vcenter{
\xymatrix@C=20pt@!R=0.5pt{
\ca{A} \ar@/^1pc/[0,4]^{\omega} \xuppertwocell[1,3]{}<3>^{K}{^<2>*!<-15pt,5pt>{\alpha_K^{-1}}}
\xtwocell[1,3]{}\omit{\omit *!<10pt,7pt>{\cong}} \ar[2,2]_{Y}
\xlowertwocell[5,3]{}<-7>_{K}{^<0.5>*!<5pt,-2pt>{\alpha_K}}
\xtwocell[5,3]{}\omit{\omit *!<-4pt,8pt>{\cong}}
\ar@/_3pc/[8,4]_{\omega}
 &&&& \Prs{B} \ar[2,0]^{W_T} \\
 &&& \Prs{B}_T \ar[-1,1]^{V_T} \xlowertwocell[1,1]{}<-1>_{\id}{^<-3>*!<-15pt,4pt>{\eta^T}} \\
&& \Prs{A} \ar[-1,1]_{L_K} \ar[3,1]_{L_K} \xlowertwocell[2,2]{}<-2>_{\id}{^<-4>{\eta^K}} && \Prs{B}_T \ar[2,0]^{\widetilde{K}}\\
\\
&&&&\Prs{A}\ar[2,0]^{L_K} \\
&&&\Prs{B}_T \ar[3,1]_{V_T}\\
&&&& \Prs{B}_T \ar[2,0]^{V_T}\\
\\
&&&&\Prs{B}
}}
\]
does satisfy the desired equation. The left adjoint $V_T$ preserves left Kan extensions (see \cite[Proposition~4.14]{KELLY_BASIC}), which implies that $V_T\alpha_K$ exhibits $V_T L_K$ as $\Lan_Y \omega$. The considerations in Section~\ref{KAN_EXTENSION_II_SECTION} therefore show that $\pi$ exhibits $V_T L_K \widetilde{K} W_T$ as $\Lan_\omega \omega$, which concludes the proof.
\end{proof}

\subsection[COEND_SECTION]{}\label{COEND_SECTION}
 For $M\in \Prs{B}_T$, we have
\[
 L_K \widetilde{K}(M)=\Prs{B}_T(K-,M)\star K= \int^{A\in\ca{A}} \Prs{B}_T(KA,M)\odot KA
\]
 (see Section~\ref{FREE_COCOMPLETION_SECTION} and Section~\ref{WEIGHTED_COLIMIT_SECTION}), and one can check that $\varepsilon^K_M \colon L_K \widetilde{K}(M) \rightarrow M$ is given by the comparison map induced by the canonical maps
\[
\Prs{B}_T(KA,M) \odot KA \rightarrow M\rlap{.}
\]
 This observation allows us to rephrase the condition that $\varepsilon^K_M$ is an isomorphism in terms of a universal property. Since this is purely formal and we do not have a specific application in mind we do not elaborate on this point. The essence of the following theorem is that we can reduce the reconstruction problem to showing that $\varepsilon^K_M$ is an isomorphism for certain $M$. We will then use this in Section~\ref{DENSITY_SECTION} to give a more tractable characterization for a special class of cosmoi (including cosmoi of modules for a commutative ring).

\begin{thm}\label{RECONSTRUCTION_THM}
 With the notation introduced in Section~\ref{NOTATION_SECTION}, the $T$-component $\nu \colon L(V^c_T) \rightarrow T$ of the counit of the Tannakian adjunction is an isomorphism if and only if for each $B\in \ca{B}$, the morphism $\varepsilon^K_{(TB,\delta_{B})} \colon L_K \widetilde{K} \bigl((TB,\delta_{B})\bigr) \rightarrow {(TB,\delta_{B})}$ is an isomorphism. This is the case if and only if the canonical maps
\[
\Prs{B}_T\bigl(KA,(TB,\delta_{B})\bigr) \odot KA \rightarrow (TB,\delta_{B})
\]
 exhibit $(TB,\delta_{B})$ as coend of
\[
 \Prs{B}_T\bigl(K(-),(TB,\delta_{B})\bigr) \odot K(-) \colon \ca{A}^{\op} \otimes \ca{A} \rightarrow \Prs{B}_T \rlap{.}
\]
\end{thm}

\begin{proof}
 We use the notation of Proposition~\ref{COUNIT_CHARACTERIZATION_PROP}. We have to show that $\nu \colon L(\omega) \rightarrow T$ is an isomorphism if and only if $\varepsilon^K_{(TB,\delta_{B})} \colon L_K \widetilde{K}(TB,\delta_{B}) \rightarrow {(TB,\delta_{B})}$ is an isomorphism. Density of the Yoneda embedding (see Section~\ref{DENSE_FUNCTORS_SECTION}) and the fact that both $L(\omega)$ and $T$ are cocontinuous (see Proposition~\ref{DENSITY_COMONAD_COCONT_PROP}) imply that $\nu$ is an isomorphism if and only if $\nu Y$ is. Since $\sigma$ is an isomorphism, $\nu Y$ is an isomorphism if and only if $V_T \varepsilon^K W_TY$ is, i.e., if and only if $V_T(\varepsilon^K_{W_T(B)})$ is an isomorphism for every $B \in \ca{B}$. The first claim thus follows from the fact that $V_T$ reflects isomorphisms (see Theorem~\ref{COMONADICITY_THM}). The second assertion follows immediately from the remarks in Section~\ref{COEND_SECTION}.
\end{proof}

\subsection[RECONSTRUCTION_REMARKS_SECTION]{}\label{RECONSTRUCTION_REMARKS_SECTION}
 If we take $\ca{V}=\CGTop$, the cosmos of compactly generated topological spaces, with monoidal structure given by the cartesian product, and $\ca{B}=\ast$, then every object $X$ of $\Prs{B}\simeq \CGTop$ has a unique comonoid structure. Moreover, giving an $X$-coaction on $Y$ is the same as giving a map $Y \rightarrow X$. Since $\CGTop$ is cartesian, the terminal object is the only object which has a dual. The category of Cauchy comodules of $X$ is therefore the full subcategory of $\CGTop/X$ generated by the objects $\ast \rightarrow X$. Clearly this category is independent of the topology on $X$, so it is impossible to reconstruct $X$ from its category of Cauchy comodules. In other words, the counit of the Tannakian adjunction need not be an isomorphism in general. In Theorem~\ref{DENSE_RECONSTRUCTION_THM} we will give a better characterization of when the maps $\varepsilon^K_{(TB,\delta_B)}$ are isomorphisms. We will discuss the relationship between the above result and the classical case of Tannaka duality where $\ca{V}$ is the cosmos of vector spaces over some field $k$ in Section~\ref{THE_CASE_OF_VECT_SECTION}

\section{Recognition}\label{RECOGNITION_SECTION}

\subsection[RECOGNITION_OVERVIEW_SECTION]{}
 In this section we study the unit of the Tannakian adjunction. In the entire section we use the notation introduced in Proposition~\ref{TANNAKIAN_ADJUNCTION_PROP}. As we will see in Proposition~\ref{UNIT_CHARACTERIZATION_PROP}, the unit of the Tannakian adjunction is closely related to the unit of the semantics-structure adjunction (Theorem~\ref{SEMANTICS_STRUCTURE_THM}). It is well-known when the latter is an equivalence (Beck's monadicity theorem, see Theorem~\ref{COMONADICITY_THM}), and this is the starting point for proving our recognition result.

\begin{prop}\label{UNIT_CHARACTERIZATION_PROP}
 Let $\omega \colon \ca{A} \rightarrow \Prs{B}$ be as in Proposition~\ref{TANNAKIAN_ADJUNCTION_PROP}. The composite of the unit $N \colon \ca{A} \rightarrow \Prs{B}^c_{L(\omega)}$ of the Tannakian adjunction with the inclusion functor $\Prs{B}^c_{L(\omega)} \rightarrow \Prs{B}_{L(\omega)}$ is naturally isomorphic to the composite
\[
\xymatrix{ \ca{A} \ar[r]^-Y & \Prs{A} \ar[r]^-{J} & \Prs{B}_{L_(\omega)}}
\]
 of the Yoneda embedding and the $(\Prs{A}, L_\omega)$-component $J \colon \Prs{A} \rightarrow \Prs{B}_{L_(\omega)}$ of the unit of the semantics-structure adjunction (see Theorem~\ref{SEMANTICS_STRUCTURE_THM}).
\end{prop}

\begin{proof}
 By Proposition~\ref{COACTION_LIFT_PROP}, we only have to check that both functors correspond to coactions which are related by a suitable natural isomorphism. By Proposition~\ref{TANNAKIAN_ADJUNCTION_PROP}, the composite of $N$ with the inclusion $\Prs{B}^c_{L(\omega)} \rightarrow \Prs{B}_{L(\omega)}$ corresponds to the $L(\omega)$-coaction
\[
\varrho=
\vcenter{
\vspace{-25pt}
\xymatrix@C=30pt@R=10pt{ \ar[r]^(0.7){Y} 
\xuppertwocell[0,4]{}<12>^{\omega} {<-4>*!<-4pt,0pt>{\alpha_\omega}}
\xtwocell[0,4]{}\omit{\omit *!<7.5pt,-21pt>{\cong}}
\xlowertwocell[1,2]{}<-6>_{\omega}{<1>*!<15pt,-4pt>{\alpha_\omega^{-1}}}
\xtwocell[1,2]{}\omit{\omit *!<-6pt,7pt>{\cong}} &
\rruppertwocell<4>^{*!<23pt,2pt>{\id}} {<0.5>*!<-3pt,0pt>{\eta^\omega}} 
\ar[rd]^{L_\omega} & &
\ar[r]^(0.4){L_\omega} & \\
& & \ar[ru]_{\widetilde{\omega}}}
}
\]
 and $J$ corresponds to the $L(\omega)$-coaction
\[
\varrho^\prime=
\vcenter{
\xymatrix@C=30pt@R=10pt{
\rruppertwocell<4>^{\id} {<0.5>*!<-3pt,0pt>{\eta^\omega}}
\ar[rd]_{L_\omega} & & \ar[r]^{L_\omega} &  \\
& \ar[ru]_{\widetilde{\omega}}
}}
\]
(see Theorem~\ref{SEMANTICS_STRUCTURE_THM} and Section~\ref{KAN_EXTENSION_SECTION}). The uniqueness part of Proposition~\ref{COACTION_LIFT_PROP} implies that $JY\colon \ca{A} \rightarrow \Prs{B}_{L_\omega \widetilde{\omega}}$ corresponds to the coaction $\varrho^\prime Y$. Clearly $L_\omega \widetilde{\omega} \alpha_\omega \circ \varrho=\varrho^\prime Y \circ \alpha_\omega$, so by Proposition~\ref{COACTION_LIFT_PROP}, $\alpha_\omega$ lifts to the desired $\ca{V}$-natural isomorphism.
\end{proof}

\subsection[OUTLINE_SECTION]{}\label{OUTLINE_SECTION}
The above proposition implies that the $(\ca{A},\omega)$-component of the unit of the Tannakian adjunction is fully faithful if $J$ is an equivalence, i.e., if $L_\omega \colon \Prs{A} \rightarrow \Prs{B}$ is comonadic (see Theorem~\ref{COMONADICITY_THM}). However, this is not to be expected in general, because the category $\Prs{B}_T$ is usually not equivalent to a category of enriched presheaves. Hence we first have to analyze the situation where $(\ca{A},\omega)$ is equal to $(\Prs{B}^c_T,V_T)$ more carefully. Writing $K\colon \ca{A} \rightarrow \Prs{B}_T$ for the inclusion functor, this is summarized by the diagram
\[
\xymatrix@C=30pt@R=30pt{
\ca{A} \ar[rd]_-{\omega} \ar[r]^-K & \Prs{B}_T \ar@<-0.5ex>[r]_-{\widetilde{K}} \ar@<0.5ex>[d]^-{V_T} & \Prs{A} \ar@<-0.5ex>[l]_-{L_K} \ar@<1ex>[ld]^{L_\omega}\\
& \Prs{B} \ar@<0.5ex>[u]^{W_T} \ar[ru]^(0.45){\widetilde{\omega}}& \\
}
\]
(see Section~\ref{NOTATION_SECTION}). Under the additional assumption that $K$ is $\ca{V}$-dense, we find that $\Prs{B}_T$ is a reflective subcategory of $\Prs{A}$ and that the comonads $L_\omega \widetilde{\omega}$ and $V_T W_T$ are isomorphic. Moreover, we have $\widetilde{K}\circ K=Y$ (see Section~\ref{FREE_COCOMPLETION_SECTION}), which suggests that we adopt the following strategy:
\begin{enumerate}
\item[a)]
 Find conditions for a $\ca{V}$-functor $\omega \colon \ca{A} \rightarrow \Prs{B}$ which imply the existence of a reflective subcategory $\ca{C}$ of $\Prs{A}$, together with a comonadic adjunction $V \colon \ca{C} \rightleftarrows \Prs{B} \colon W$ such that the density comonad $VW$ of $V$ is equal to $L(\omega)$,
\item[b)]
 Find conditions for a $\ca{V}$-functor $\omega$ which imply that the Yoneda embedding $Y\colon \ca{A}\rightarrow \Prs{A}$ factors through the embedding $\ca{C}\rightarrow \Prs{A}$ from a).
\end{enumerate}
 Under these conditions, it is not hard to see that we get a commutative diagram
\[
 \xymatrix@!=15pt{ &\ca{A} \ar@{-->}[ld] \ar[rd]^{Y} \\
 \ca{C} \ar[rd]_{J^\prime} \ar[rr] && \Prs{A} \ar[ld]^{J} \\
 & \Prs{B}_{L(\omega)}}
\]
 where $J^\prime$ and $J$ are the respective components of the unit of the semantics-structure adjunction. Thus, if we succeed with showing a) and b) above, we know that $J^\prime$ is an equivalence, and therefore that $JY$ is fully faithful. By Proposition~\ref{UNIT_CHARACTERIZATION_PROP}, it follows that the unit $N \colon \ca{A} \rightarrow \Prs{B}^c_{L(\omega)}$ of the Tannakian adjunction is fully faithful. This is the hard part of our proof; showing that $N$ is essentially surjective is straightforward.

 Since comonadic functors reflect isomorphisms, a reasonable candidate for $\ca{C}$ is the full subcategory of objects which are orthogonal to the class $\Sigma$ of morphisms in $\Prs{A}$ which get sent to isomorphisms by $L_\omega \colon \Prs{A} \rightarrow \Prs{B}$ (an object $X \in \Prs{A}$ is \emph{orthogonal} to a class $\Sigma$ of morphisms if, for all $f \in \Sigma$, the morphism $\Prs{A}(f,X)$ is an isomorphism; see \cite[Section~6.2]{KELLY_BASIC}). If this subcategory is reflective, then it is the localization of $\Prs{A}$ at $\Sigma$, so it can be thought of as the category obtained by formally inverting the morphisms in $\Sigma$. We first show that a reflection functor exists if the monoidal category $\ca{V}$ is locally presentable.

 Recall that we denote the underlying ordinary category of a $\ca{V}$-category $\ca{C}$ by $\ca{C}_0$, and that we denote the underlying ordinary functor of a $\ca{V}$-functor $F \colon \ca{C} \rightarrow \ca{C}^{\prime}$ by $F_0 \colon \ca{C}_0 \rightarrow \ca{C}^{\prime}_0$. For any small $\ca{V}$-category $\ca{A}$, the underlying ordinary category $\Prs{A}_0$ of the category $\Prs{A}$ of enriched presheaves on $\ca{A}$ has objects the $\ca{V}$-functors $\ca{A}^{\op} \rightarrow \ca{V}$, and morphisms the $\ca{V}$-natural transformations between them. This is not to be confused with $\mathcal{P}(\ca{A}_0)$, the category of ordinary $\Set$-valued presheaves on the ordinary category $\ca{A}_0$.

\begin{lemma}\label{PRESHEAVES_LOCALLY_PRESENTABLE_LEMMA}
Let $\ca{A}$ be a small $\ca{V}$-category. If $\ca{V}_0$ is locally presentable, then $\Prs{A}_0$ is locally presentable.
\end{lemma}

\begin{proof}
Since $\Prs{A}$ is cocomplete as a $\ca{V}$-category, its underlying ordinary category $\Prs{A}_0$ is cocomplete, so it suffices to show that there is a regular cardinal $\lambda$ and a set of $\lambda$-small objects which forms a strong $\Set$-generator of $\Prs{A}_0$. By assumption we can choose a regular  cardinal $\lambda$ such that $\ca{V}_0$ is locally $\lambda$-presentable. Thus there exists a set $\ca{G}$ of $\lambda$-small objects of $\ca{V}_0$ which forms a strong $\Set$-generator of $\ca{V}_0$. For $A \in \ca{A}$ and $G \in \ca{G}$, the object $\ca{A}(-,A)\odot G$ is $\lambda$-small: by definition of tensor products (see Section~\ref{WEIGHTED_COLIMIT_SECTION}) and by Yoneda, we have for each $F \in \Prs{A}_0$
\[
 \Prs{A}_0\bigl(\ca{A}(-,A)\odot G,F\bigr)\cong [G, \Prs{A}\bigl(\ca{A}(-,A),F\bigr)]_0 \cong [G,FA]_0\cong \ca{V}_0(G,FA)\rlap{,}
\]
and both $\ca{V}(G,-)$ and $F \mapsto FA$ preserve $\lambda$-filtered colimits. In fact, the latter preserves all colimits, because colimits in $\Prs{A}$ are computed pointwise.

 Let $\alpha \colon F \Rightarrow F^{\prime}$ be a $\ca{V}$-natural transformation, i.e., a morphism in $\Prs{A}_0$. Since the above isomorphism is natural, we find that $\Prs{A}_0\bigl(\ca{A}(-,A)\odot G,\alpha\bigr)$ is an isomorphism if and only if $\ca{V}(G,\alpha_A)$ is. This implies that $\{\ca{A}(-,A)\odot G \vert G\in \ca{G}, A \in \ca{A}\}$ forms a strong $\Set$-generator of $\Prs{A}_0$. Thus $\Prs{A}_0$ is locally presentable.
\end{proof}

\begin{lemma}\label{ORTHOGONAL_SUBCATEGORY_LEMMA}
 Let $\omega \colon \ca{A} \rightarrow \Prs{B}$ and $L_\omega \colon \Prs{A} \rightleftarrows \Prs{B} \colon \widetilde{\omega}$ be as in Proposition~\ref{TANNAKIAN_ADJUNCTION_PROP} and Section~\ref{FREE_COCOMPLETION_SECTION} respectively, and let $\Sigma$ be the class of morphisms $f$ in $\Prs{A}_0$ for which $L_\omega f$ is an isomorphism. If $\ca{V}_0$ is locally presentable, then the full subcategory $\ca{C}=\Sigma^\perp$ of $\Prs{A}$ of objects which are orthogonal to $\Sigma$ is reflective, i.e., the inclusion functor $I \colon \ca{C} \rightarrow \Prs{A}$ has a left adjoint $R \colon \Prs{A} \rightarrow \ca{C}$.
\end{lemma}

\begin{proof}
 The functor $L_\omega$ preserves tensor products, so for any morphism $f\colon X \rightarrow X^{\prime}$ in $\Sigma$ and any object $V \in \ca{V}$, the morphism $V\odot f \colon V\odot X \rightarrow V \odot X^{\prime}$ is in $\Sigma$. Hence an object $F \in \Prs{A}_0$ is orthogonal to $\Sigma$ in the enriched sense if and only if it is orthogonal in the classical sense (cf. \cite[Section~6.2]{KELLY_BASIC}).

 By Lemma~\ref{PRESHEAVES_LOCALLY_PRESENTABLE_LEMMA}, the categories $\Prs{A}_0$ and $\Prs{B}_0$ are locally presentable, and the functor $(L_\omega)_0 \colon \Prs{A}_0 \rightarrow \Prs{B}_0$ preserves all colimits because it is a left adjoint, so it is accessible. It follows from \cite[Remark~2.50 and \S~2.60]{ADAMEK_ROSICKY} that the full subcategory $\Sigma$ of $\Mor(\Prs{A}_0)$ is an accessible, accessibly embedded subcategory. Thus there is a regular cardinal $\lambda$ and a subset $\Sigma_0 \subseteq \Sigma$ such that every element of $\Sigma$ is a $\lambda$-filtered colimit of elements of $\Sigma_0$, and this colimit is computed as in $\Mor(\Prs{A}_0)$. It follows immediately that $\Sigma^\perp=\Sigma_0^\perp$, and \cite[Theorem~1.39]{ADAMEK_ROSICKY} implies that the inclusion $I_0 \colon \ca{C}_0 \rightarrow \Prs{A}_0$ has a left adjoint $R_0 \colon \Prs{A}_0 \rightarrow \ca{C}_0$. The claim now follows from \cite[Theorem~4.85]{KELLY_BASIC}.
\end{proof}

\begin{lemma}\label{ADJUNCTION_FACTORIZATION_LEMMA}
 Let $\ca{V}_0$ be locally presentable. With the notation from Lemma~\ref{ORTHOGONAL_SUBCATEGORY_LEMMA}, $\widetilde{\omega} \colon \Prs{B} \rightarrow \Prs{A}$ factors through the inclusion $I \colon \ca{C} \rightarrow \Prs{A}$, i.e., there is a functor $W \colon \Prs{B} \rightarrow \ca{C}$ such that $\widetilde{\omega}=IW$. The restriction $V=L_\omega I$ of $L_\omega$ to $\ca{C}$ is left adjoint to $W$, and $V$ reflects isomorphisms.
\end{lemma}

\begin{proof}
 We have to check that for any $Y \in \Prs{B}$, the object $\widetilde{\omega}(Y)$ lies in $\ca{C}$, i.e., that for any morphism $f \colon X \rightarrow X^{\prime}$ in $\Sigma$, the map 
\[
\Prs{A}\bigl(f,\widetilde{\omega}(Y)\bigr) \colon \Prs{A}\bigl(X^{\prime},\widetilde{\omega}(Y)\bigr) \rightarrow
\Prs{A}\bigl(X,\widetilde{\omega}(Y)\bigr)
\]
 is an isomorphism. By adjunction, this is is equivalent to $\Prs{A}(L_\omega f,Y)$ being an isomorphism, which is evident: the morphism $L_\omega f$ is an isomorphism by definition of $\Sigma$ (see Lemma~\ref{ORTHOGONAL_SUBCATEGORY_LEMMA}). It is clear that the restriction of a left adjoint still gives a left adjoint, so it remains to check that $V$ reflects isomorphisms.

 Let $\eta$ and $\varepsilon$ be the unit and counit of the adjunction $R \colon \Prs{A} \rightleftarrows \ca{C} \colon I$ from Lemma~\ref{ORTHOGONAL_SUBCATEGORY_LEMMA}, and let $\eta_0$ and $\varepsilon_0$ be the unit and counit of $V \colon \ca{C} \rightleftarrows \Prs{B} \colon W$. We clearly have $I\eta_0=\eta^{\omega} I$ and $\varepsilon_0=\varepsilon^{\omega}$, so it might seem obfuscating to use different names, but doing so clarifies the pasting diagram
\[
\vcenter{ \xymatrix@C=15pt@R=10pt{
&& \ar[rd]|{R} \xuppertwocell[0,5]{}^{\id}{\omit} \xuppertwocell[1,2]{}<5>^{\id}{<-0.5>\eta} & \xuppertwocell[0,4]{}\omit{<0.5>*!<-4pt,0pt>{\eta^{\omega}}}&&&&\\
& \ar[ru]^{I} \rrlowertwocell<-3>_{\id}{<-0.5>*!<12pt,0pt>{\varepsilon}}& & \ar[r]_{I} & \ar[rd]^{L_\omega}\\
\ar[ru]^{W} \xlowertwocell[0,5]{}_{\id}{<0.75>*!<-2pt,0pt>{\varepsilon_0}} &&&&& \ar[rruu]_{\widetilde{\omega}}
}}=\vcenter{\xymatrix@C=20pt@R=20pt{
& \rruppertwocell^{\id}<3>{<1>*!<-4pt,0pt>{\eta^{\omega}}} \ar[rd]|{L_\omega} & & \\
\ar[ru]^{IW} \rrlowertwocell<-3>_{\id}{<-1>*!<14pt,0pt>{\varepsilon^{\omega}}} & &  \ar[ru]_{\widetilde{\omega}}
}}=\xymatrix@1{\rrtwocell^{IW}_{\widetilde{\omega}}{*!<-2pt,0pt>{\id}} & & }
\]
 which shows that $L_\omega \eta$ is the `mate' of the identity natural transformation $IW=\widetilde{\omega}$ under the adjunctions $L_\omega \dashv \widetilde{\omega}$ and $VR \dashv IW$, in the sense of \cite{KELLY_STREET}. Since the mate of a natural isomorphism is an isomorphism, it follows that $L_\omega \eta$ is an isomorphism, i.e., that the components of $\eta$ lie in $\Sigma$. If we apply $L_\omega$ to the commutative diagram
\[
 \xymatrix{X \ar[r]^-{\eta_X} \ar[d]_{f} & IR(X) \ar[d]^{IR(f)} \\
X^{\prime} \ar[r]_-{\eta_{X^{\prime}}} & IR(X^{\prime})}
\]
 we find that a morphism $f$ of $\Prs{A}_0$ lies in $\Sigma$ if and only if $R(f)=IR(f)$ lies in $\Sigma$. Since $R(f)$ is a morphism whose domain and codomain lie in $\ca{C}=\Sigma^\perp$, this is the case if and only if $R(f)$ is an isomorphism. Now let $f$ be a morphism in $\ca{C}$ such that $Vf=L_\omega I(f)$ is an isomorphism. Then $I(f)$ lies in $\Sigma$, and the above reasoning shows that $RI(f)$ is an isomorphism. But $\varepsilon \colon RI \Rightarrow \id$ is a natural isomorphism, which implies that $f$ is an isomorphism. Thus $V$ reflects isomorphisms, as claimed.
\end{proof}

\begin{thm}\label{RECOGNITION_THM}
 Let $\ca{V}$ be a cosmos with $\ca{V}_0$ locally presentable. With the notation of Definition~\ref{TANNAKIAN_ADJUNCTION_DFN}, let
\[
 \xymatrix{\ca{A} \ar[rr]^-{N} \ar[rd]_{\omega}& & \Prs{B}_{L(\omega)}^c \ar[ld]^{V^c_{L(\omega)}}\\
& \Prs{B}}
\]
 be the unit of the Tannakian adjunction, and write $\Sigma$ for the class of morphisms $f$ of $\Prs{A}_0$ for which $L_\omega(f)$ is an isomorphism. The $\ca{V}$-functor $N \colon \ca{A} \rightarrow \Prs{B}_{L(\omega)}^c$ is an equivalence of $\ca{V}$-categories if the following hold:
\begin{enumerate}
\item[i)]
 The functor $\omega \colon \ca{A} \rightarrow \Prs{B}$ reflects isomorphisms,
\item[ii)]
 The left adjoint $L_\omega \colon \Prs{A} \rightarrow \Prs{B}$ preserves equalizers of $L_\omega$-cosplit pairs in $\Sigma^\perp$ (see Section~\ref{COSPLIT_EQUALIZER_SECTION}), and
\item[iii)]
 Whenever $J \colon \ca{D}^{\op} \rightarrow \ca{V}$ is a weight and $G \colon \ca{D} \rightarrow \ca{A}$ is a $\ca{V}$-functor such that the weighted colimit $J \star \omega G$ lies in the subcategory $\overline{\ca{B}} \subseteq \Prs{B}$, then $J\star G$ exists and is preserved by $\omega \colon \ca{A} \rightarrow \Prs{B}$.
\end{enumerate}
 Moreover, if $\mathbf{\Phi}$ is a class of weights with the property that for each $X \in \Prs{A}$ there is a weight $J \colon \ca{D}^{\op} \rightarrow \ca{V}$ in $\mathbf{\Phi}$ and a $\ca{V}$-functor $G \colon \ca{D} \rightarrow \ca{A}$ such that $X \cong J \star YG$, then $N \colon \ca{A} \rightarrow \Prs{B}_{L_\omega \widetilde{\omega}}^c$ is an equivalence if i) and ii) hold, and iii) holds for all weights in the class $\mathbf{\Phi}$.
\end{thm}

\begin{proof}
 Note that it suffices to prove the statement involving a class $\mathbf{\Phi}$ of weights. We can consider the Yoneda embedding $Y \colon \ca{A} \rightarrow \Prs{A}$ as a diagram on $\ca{A}$, and any $X \in \Prs{A}$ as a weight $\ca{A}^{\op} \rightarrow \ca{V}$. In \cite[Formula~(3.17)]{KELLY_BASIC} it was shown that the weighted colimit $X \star Y$ is isomorphic to $X$. Thus we can always let $\mathbf{\Phi}$ be the class of all weights $\ca{A}^{\op} \rightarrow \ca{V}$.

 By Proposition~\ref{UNIT_CHARACTERIZATION_PROP}, the composite of the unit $N \colon \ca{A} \rightarrow \Prs{B}_{L(\omega)}^c$ with the inclusion $\Prs{B}_{L(\omega)}^c \rightarrow \Prs{B}_{L(\omega)}$ is isomorphic to the composite
\[
\xymatrix{ \ca{A} \ar[r]^Y & \Prs{A} \ar[r]^-{J} & \Prs{B}_{L(\omega)}}
\]
 of the Yoneda embedding and the $(\Prs{A},L_\omega)$-component of the unit of the semantics structure adjunction. Let $\ca{C}$ be the full subcategory $\Sigma^\perp$ of $\Prs{A}$ consisting of objects which are orthogonal to $\Sigma$. Let $R \colon \ca{C} \rightleftarrows \Prs{A} \colon I$ be as in Lemma~\ref{ORTHOGONAL_SUBCATEGORY_LEMMA}, and let $V \colon \ca{C} \rightleftarrows \Prs{B} \colon W$ be as in Lemma~\ref{ADJUNCTION_FACTORIZATION_LEMMA}. We write $J^{\prime} \colon \ca{C} \rightarrow \Prs{A}_{VW}=\Prs{A}_{L_\omega \widetilde{\omega}}$ for the $(\ca{C},V)$ component of the unit of the semantics structure adjunction (see Theorem~\ref{SEMANTICS_STRUCTURE_THM}). Since the unit $\eta_0$ of $V \dashv W$ satisfies $I \eta_0=\eta^\omega I$ by construction, Proposition~\ref{COACTION_LIFT_PROP} shows that the diagram
\[
\xymatrix{ \ca{C} \ar[rr]^-I \ar[rd]_-{J^{\prime}} && \Prs{A} \ar[ld]^-{J}\\
& \Prs{B}_{L(\omega)}^c}
\]
 is commutative (cf.\ Theorem~\ref{SEMANTICS_STRUCTURE_THM} and Section~\ref{KAN_EXTENSION_SECTION}). By ii), the $\ca{V}$-functor $V=L_\omega I \colon \ca{C} \rightarrow \Prs{B}$ preserves equalizers of $V$-cosplit pairs, and by Lemma~\ref{ADJUNCTION_FACTORIZATION_LEMMA} it reflects isomorphisms. Thus $J^{\prime}$ is an equivalence (see Theorem~\ref{COMONADICITY_THM}). We will now show that the Yoneda embedding factors through the full subcategory $\ca{C}=\Sigma^{\perp}$, and therefore that $JY=J^{\prime}Y$ is fully faithful. If this holds, then condition iii) implies that all objects $X$ of $\ca{C}$ with $V(X) \in \overline{\ca{B}}$ lie in the essential image of the Yoneda embedding, and since the equivalence $J^{\prime} \colon \ca{C} \rightarrow \Prs{B}_{L(\omega)}$ satisfies $V_{L(\omega)} J^{\prime} =V$, it follows that $J^\prime Y \colon \ca{A} \rightarrow \Prs{B}^c_{L(\omega)}$ is an equivalence.

 In other words, it remains to show that for each $A\in \ca{A}$, the presheaf $\ca{A}(-,A)$ lies in $\Sigma^{\perp}$. Since $L_\omega \colon \Prs{A} \rightarrow \Prs{B}$ is left adjoint, $\Sigma$ is closed under tensor products, so it suffices to check that for any $f_0 \colon X \rightarrow X^{\prime}$ in $\Sigma$, and any $g \colon X \rightarrow \ca{A}(-,A)$, the unique dotted arrow
\[
 \xymatrix{X \ar[d]_{f_0} \ar[r]^-{g} & \ca{A}(-,A) \\ X^{\prime} \ar@{-->}[ru] }
\]
 exists. Since $L_\omega$ preserves colimits, the pushout $f \colon \ca{A}(-,A) \rightarrow X^{\prime\prime}$ of $f_0$ along $g$ lies in $\Sigma$, too. It follows that the dotted arrow in the above diagram exists and is unique if and only if there is a unique dotted arrow making
\[
\xymatrix{\ca{A}(-,A) \ar[d]_{f} \ar@{=}[r] & \ca{A}(-,A) \\ X^{\prime\prime} \ar@{-->}[ru] &}
\]
 commutative. By assumption, there is a weight $J \colon \ca{D}^{\op} \rightarrow \ca{V}$ in the class $\mathbf{\Phi}$ and a functor $G \colon \ca{D} \rightarrow \ca{A}$ such that $X^{\prime\prime} \cong J\star YG$, so we might as well assume that $f$ is a morphism from $\ca{A}(-,A)$ to $J\star YG$.
 
 \textbf{Step 1: Existence.} By the remarks in Section~\ref{WEIGHTED_COLIMIT_SECTION} we have a chain of isomorphisms
\[
 \xymatrix@C=25pt{\omega(A) \ar[r]^-{\alpha_\omega} & L_\omega YA \ar[r]^-{L_\omega f} & L_\omega (J\star YG) \ar[r]^-{\widehat{L_\omega}^{-1}} & J\star L_\omega YG \ar[r]^-{J\star \alpha_\omega G}  & J\star \omega G},
\]
 which shows that $J\star \omega G$ lies in $\overline{\ca{B}}$. By iii) it follows that the colimit $J \star G \in \ca{A}$ exists, and that $\widehat{\omega} \colon J\star \omega G \rightarrow \omega(J\star G)$ is an isomorphism. The fact that $\alpha_\omega \colon \omega \Rightarrow L_\omega Y$ is an isomorphism implies that $\widehat{L_\omega Y} \colon J\star L_\omega G \rightarrow L_\omega Y(J\star G)$ is an isomorphism. An easy application of the Yoneda Lemma and the definition of the comparison morphism (see Section~\ref{WEIGHTED_COLIMIT_SECTION}) show that $\widehat{L_\omega Y}=L_\omega (\widehat{Y}) \circ \widehat{L_\omega}$, and it follows that $L_\omega (\widehat{Y})$ is an isomorphism. Since the Yoneda embedding is full, the composite $\widehat{Y} \circ f \colon \ca{A}(-,A) \rightarrow \ca{A}(-,J\star G)$ is of the form $Yh=\ca{A}(-,h)$ for a unique morphism $h \colon A \rightarrow J\star G$. Note that $L_\omega Yh=L_\omega(\widehat{Y}) \circ L_\omega (f)$ is an isomorphism, by the above argument and by assumption on $f$. Since $\omega \cong L_\omega Y$ it follows that $\omega(h)$ is an isomorphism, and condition i) implies that $h$ itself is an isomorphism. Thus
\[
 \xymatrix{\ca{A}(-,A) \ar@{=}[r] \ar[d]_{f} & \ca{A}(-,A) \\
J\ast YG \ar[d]_{\widehat Y} \\
\ca{A}(-,J\star G) \ar[ruu]_{\ca{A}(-,h^{-1})}}
\]
 gives the desired lift.

 \textbf{Step 2: Uniqueness.} We use the notation introduced above. We have to show that for any $k\colon J\star YG \rightarrow \ca{A}(-,A)$ with $kf=\id$, the equality $k=Yh^{-1} \circ \widehat{Y}$ holds. This follows from the fact that the target $\ca{A}(-,A)$ of $k$ lies in the image of the Yoneda embedding, and within that image, $Y(J\star G)$ \emph{is} the colimit of $YG$ weighted by $J$, which means that $\ca{A}(-,A)$ cannot detect whether or not  $\widehat{Y} \colon J\star YG \rightarrow Y(J\star G)$ is an isomorphism. We explain this more carefully. From Section~\ref{WEIGHTED_COLIMIT_SECTION} we know that there is a unique morphism $k^\prime \colon J\star G \rightarrow A$ such that part (2) of the diagram
\[
 \xymatrix{& J \ar@{}[rddd]|(0.4){(2)} \ar@{}[lddd]|(0.4){(1)} \ar[ld] \ar[rd] \ar[dd]\\ 
 \Prs{A}(YG-,J\star YG)\ar[dd]_{\Prs{A}(YG-,\widehat{Y})} && \Prs{A}(YG-,J\star YG) \ar[d]^{\Prs{A}(YG-,k)} \\
 &\ca{A}(G-,J\star G) \ar@{}[dd]|{(3)} \ar[ld]_-{Y} \ar[rd]|{\ca{A}(G-, k^{\prime})}& \Prs{A}\bigl(YG-,\ca{A}(-,A)\bigr) \ar[d]^{Y^{-1}}\\
 \Prs{A}\bigl(YG-,Y(J\star G)\bigr) \ar[rd]|{\Prs{A}(YG-,Yk^{\prime})} & & \ca{A}(G-,A) \ar[ld]^-{Y}\\
 & \Prs{A}(YG-,YA)}
\]
 is commutative, where the unlabeled arrows denote the respective units. Part (1) is commutative by definition of $\widehat{Y}$ (see Section~\ref{WEIGHTED_COLIMIT_SECTION}), and part (3) is commutative by functoriality of the Yoneda embedding. It follows that $Yk^{\prime} \circ \widehat{Y}=k$, and therefore that $Yk^{\prime} \circ Yh= Yk^{\prime} \circ \widehat{Y} \circ f =kf=\id$. But $h$ is an isomorphism, hence we must have $k^{\prime}=h^{-1}$ and $k=Yk^\prime \circ \widehat{Y}=Y(h^{-1})\circ \widehat{Y}$, which shows that the lift is unique.
\end{proof}

\subsection[RECOGNITION_REMARKS_SECTION]{}\label{RECOGNITION_REMARKS_SECTION}
 We now take a closer look at the different conditions in Theorem~\ref{RECOGNITION_THM}. From Proposition~\ref{COLIMIT_CREATION_PROP} it follows that iii) is a necessary condition, and from Theorem~\ref{COMONADICITY_THM} we know that i) is a necessary condition too. On the other hand, in the proof of Theorem~\ref{RECOGNITION_THM} we have seen that under the given conditions, the inclusion functor $K \colon \ca{A} \rightarrow \Prs{B}_{L(\omega)}$ of the category of Cauchy comodules of $L(\omega)$ in the category of all comodules is a dense functor (cf.\ \cite[Proposition~5.16]{KELLY_BASIC}). But this implies that the right adjoint $\widetilde{K} \colon \Prs{B} \rightarrow \Prs{A}$ is fully faithful. Therefore the counit $\varepsilon^K$ is an isomorphism, so by Theorem~\ref{RECONSTRUCTION_THM} it follows that the $L(\omega)$-component of the counit of the Tannakian adjunction is an isomorphism. Thus Theorem~\ref{RECOGNITION_THM} does not provide a complete characterization of categories of comodules if the reconstruction problem can't be solved, for example if $\ca{V}=\CGTop$ (see Section~\ref{RECONSTRUCTION_REMARKS_SECTION}). Also, this condition is probably the most difficult to check. In some cases it might be easier to check whether $L_\omega$ preserves \emph{all} weighted finite limits. If this is the case, $\omega$ is called a \emph{flat} $\ca{V}$-functor (cf.\ \cite[\S~6]{KELLY_FINLIM}). For $\ca{V}=\Vect$, the notion of a flat functor defined in the introduction (see Section~\ref{DESCRIPTION_OF_RESULTS_SECTION}) agrees with the notion from \cite{KELLY_FINLIM}, but for general cosmoi the notion of flatness is a bit more subtle. In the next section, we study the recognition question for cosmoi with dense autonomous generator. In this case we can avoid the more general definition of flat functors.

\section{Cosmoi with dense autonomous generator}\label{DENSITY_SECTION}

\subsection[CONICAL_COLIMIT_SECTION]{}\label{CONCIAL_COLIMIT_SECTION}
 When working with additive $R$-linear categories for some commutative ring $R$, the notion of weighted colimits is generally not needed. For example, any $R$-linear presheaf is a conical colimit of representable functors. Any additive $R$-linear category has in particular tensor products with finitely generated free $R$-modules: we have $R^n\odot A \cong \oplus_{i=1}^n A$. The full subcategory of finitely generated free $R$-modules is a $\Set$-dense subcategory which is closed under the tensor product and under the formation of duals.

\begin{dfn}\label{DAG_DFN}
 An essentially small full subcategory $\ca{X} \subseteq \ca{V}$ of a cosmos $\ca{V}$ is called a \emph{dense autonomous generator} if $\ca{X}$ consists of objects with duals, is closed under the tensor product and under the formation of duals, and is $\Set$-dense in $\ca{V}$.
\end{dfn}

 In Appendix~\ref{DENSITY_IN_COSMOI_WITH_DAG_APPENDIX}, we show that cosmoi $\ca{V}$ with dense autonomous generator $\ca{X}$ and $\ca{X}$-tensored $\ca{V}$-categories share some of the nice properties of additive $R$-linear categories. In the $R$-linear context it is not hard to prove the following result directly.

\begin{cit}[Corollary~\ref{CONICAL_COLIMIT_COR}]
 Let $\ca{V}$ be a cosmos which has a dense autonomous generator $\ca{X}$. Let $\ca{A}$ be a small $\ca{X}$-tensored $\ca{V}$-category. Fix a presheaf $F \in \Prs{A}$, and let $\ca{A}\slash F$ be the category of representable functors over $F$. Then $F$ is the conical colimit of the domain functor $D \colon \ca{A} \slash F \rightarrow \Prs{A}$.
\end{cit}

\begin{thm}\label{DENSE_RECONSTRUCTION_THM}
 Let $\ca{V}$ be a cosmos which has a dense autonomous generator $\ca{X}$. Let $\ca{B}$ be a $\ca{V}$-category with small Cauchy completion. Then the counit $\nu \colon L(V^c_T) \rightarrow T$ of the Tannakian adjunction is an isomorphism if and only if for each $B \in \ca{B}$, the tautological cocone on the domain functor $D \colon \Prs{B}^c_T \slash (TB,\delta_B) \rightarrow \Prs{B}_T$ exhibits $(TB,\delta_B)$ as the conical colimit of $D$.
\end{thm}

\begin{proof}
 We use the notation from Theorem~\ref{RECONSTRUCTION_THM}. By Proposition~\ref{COLIMIT_CREATION_PROP} it follows that the category $\ca{A}=\Prs{B}_T^c$ is $\ca{X}$-tensored (note that for $B \in \overline{\ca{B}}$, $X \in \ca{X}$, we have $X \odot B \in \overline{\ca{B}}$ because $\ca{X}$ consists of objects with duals). By Corollary~\ref{CONICAL_COLIMIT_COR}, the presheaf $F = \Prs{B}_T\bigl(K-,(TB,\delta_B)\bigr)$ on $\ca{A}$ is the colimit of the diagram of representable functors over $F$. We have $L_K(\ca{A}(-,A)) \cong KA$ and $L_K$ preserves colimits, so we find that $L_K \widetilde{K}(TB,\delta_B)$ is the colimit of $D \colon \ca{A}\slash (TB,\delta_B) \rightarrow \Prs{B}_T$ of Cauchy comodules over $(TB,\delta_B)$. It is not hard to check that the comparison morphism induced by the tautological cocone is precisely the $(TB,\delta_B)$-component of $\varepsilon^K$, so the conclusion follows from Theorem~\ref{RECONSTRUCTION_THM}.
\end{proof}

\subsection[THE_CASE_OF_VECT_SECTION]{}\label{THE_CASE_OF_VECT_SECTION}
 In Section~\ref{GENERAL_OVERVIEW_SECTION} we have seen that the reconstruction theorem of classical Tannakian duality can be restated: the reconstruction theorem says that the counit of the Tannakian adjunction is an isomorphism for $\ca{V}$ the cosmos $\Vect$ of vector spaces over a field $k$, and $\ca{B}=k$, the unit $k$-linear category. Using the above theorem it is now quite easy to show this. Let $(C,\delta,\varepsilon)$ be a coalgebra in $\Vect$. The fact that the $C$-comodule $(C,\delta)$ is a union of finite dimensional subcomodules implies that the comparison morphism $\alpha \colon \colim D \rightarrow (C,\delta)$ is surjective. Any element of $\colim D$ lies in the image of one of its structure maps $\kappa_\varphi$ for a finite dimensional comodule $\varphi\colon M\rightarrow (C,\delta)$ over $(C,\delta)$. Now let $x_1,x_2 \in \colim D$ with $\alpha(x_1)=\alpha(x_2)$. For $i=1,2$ let $\varphi_i \colon M_i \rightarrow (C,\delta)$ be finite dimensional comodules over $(C,\delta)$, with elements $x_i^\prime \in M$ satisfying $\kappa_{\varphi_i}(x_i^\prime)=x_i$. By definition of the comparison map we have $\alpha \kappa_{\varphi_i}=\varphi_i$, hence $\varphi_1(x_1^\prime)=\varphi_2(x_2^\prime)$. It follows that there is an element $y$ in the pullback $N$ of $\varphi_1$ and $\varphi_2$ which gets sent to $x_1^\prime$ and $x_2^\prime$ under the respective projection maps $N \rightarrow M_i$. These projection maps are morphisms in the category of finite dimensional comodules over $(C,\delta)$, so the structure map $N \rightarrow \colim D$ corresponding to $N \rightarrow (C,\delta)$ sends $y$ to $x_i$ for $i=1,2$. Thus $x_1=x_2$, and we have shown that $\alpha$ is an isomorphism.

 Note that we have used the fact that, in $\Vect$, a pullback of two objects with duals is again an object with a dual. This argument can be generalized to flat coalgebras over Dedekind rings (cf. \cite{WEDHORN}), precisely because a submodule of a finitely generated projective module over a Dedekind ring is again finitely generated and projective. We can generalize this result further to flat coalgebras $C$ over an arbitrary commutative ring $R$ which have \emph{enough Cauchy comodules}, meaning that for every $C$-comodule $M$ and every $m\in M$ there is a Cauchy comodule $M^\prime$ and a morphism $M^\prime \rightarrow M$ of $C$-comodules whose image contains $m$.

\begin{cor}
 Let $R$ be a commutative ring and let $C$ be a flat coalgebra which has enough Cauchy comodules (see above). Then the $C$-component of the counit of the Tannakian adjunction is an isomorphism.
\end{cor}

\begin{proof}
 We proceed as in Section~\ref{THE_CASE_OF_VECT_SECTION}. Since $C$ has enough comodules it follows immediately that the comparison morphism $\alpha \colon \colim D \rightarrow (C,\delta)$ is surjective. To show that it is also injective we form the pullback comodule $N$ as in Section~\ref{THE_CASE_OF_VECT_SECTION}. Flatness of $C$ implies that this pullback is computed as in the category of $R$-modules. We thus have the desired element $y\in N$, but $N$ itself need no longer be a Cauchy comodule. But $C$ has enough Cauchy comodules, so there is a Cauchy comodule $M$ with an element $z\in M$ and a morphism of comodules $\varphi \colon M \rightarrow N$ with $\varphi(z)=y$. The structure map corresponding to the composite $M \rightarrow N \rightarrow (C,\delta)$ sends $z$ to $x_i$, $i=1,2$, which shows that the comparison map is an isomorphism.
\end{proof}

 We now turn to a discussion of the recognition principle in cosmoi with dense autonomous generator.

\begin{dfn}\label{OMEGA_RIGID_EXTENDED_DFN}
 Let $\ca{A}$, $\ca{B}$ be small $\ca{V}$-categories, and let $\omega \colon \ca{A} \rightarrow \Prs{B}$ be a $\ca{V}$-functor, with underlying ordinary functor $\omega_0 \colon \ca{A}_0 \rightarrow \Prs{B}_0$. Let $\ca{D}$ be a small ordinary category. A diagram $D \colon \ca{D} \rightarrow \ca{A}_0$ is called \emph{$\omega$-rigid} if the colimit of $\omega_0 D \colon \ca{D} \rightarrow \Prs{B}_0$ lies in the Cauchy completion $\overline{\ca{B}}$ of $\ca{B}$. This extends Definition~\ref{OMEGA_RIGID_DFN} of $\omega$-rigid diagrams for $\Prs{B}=\ca{V}$.
\end{dfn}

\begin{thm}\label{DENSE_RECOGNITION_THM}
 Let $\ca{V}$ be a cosmos which has a dense autonomous generator $\ca{X}$, and assume that the functor $\ca{V}_0(I,-) \colon \ca{V}_0 \rightarrow \Set$ preserves filtered colimits. Let $\ca{A}$ be a small $\ca{X}$-tensored $\ca{V}$-category, let $\ca{B}$ be a small $\ca{V}$-category, and let $\omega \colon \ca{A} \rightarrow \Prs{B}$ be a $\ca{V}$-functor whose image is contained in $\overline{\ca{B}}$. Then the $(\ca{A},\omega)$-component of the unit of the Tannakian adjunction is an equivalence if
\begin{enumerate}
 \item[i)]
 The functor $\omega$ reflects isomorphisms,
 \item[ii)]
 For each object $B\in\ca{B}$, the category $\el \bigl(V (\ev_B \omega)_0\bigr)$ of elements of the functor $V(\ev_B\omega)_0 \colon \ca{A}_0 \rightarrow \Set$ is cofiltered, where $\ev_B \colon \Prs{B} \rightarrow \ca{V}$ denotes the evaluation functor, and
 \item[iii)]
 The category $\ca{A}_0$ has colimits of $\omega$-rigid diagrams (see Definition~\ref{OMEGA_RIGID_EXTENDED_DFN}), and $\omega_0$ preserves them.
\end{enumerate}
If these conditions are satisfied, then the comonad $L(\omega) \colon \Prs{B} \rightarrow \Prs{B}$ preserves finite limits.
\end{thm}

\begin{proof}
 We use the notation from Proposition~\ref{TANNAKIAN_ADJUNCTION_PROP} and Theorem~\ref{RECOGNITION_THM}. First, note that $\ca{V}_0$ is locally finitely presentable: for $X\in \ca{X}$ we have
\[
 \ca{V}_0(X,-)\cong \ca{V}_0(I,[X,-])\cong\ca{V}_0(I,X^\vee \otimes -),
\]
 which preserves filtered colimits by our assumption on the unit object $I$. But $\ca{X}_0$ is $\Set$-dense, so it is in particular a strong generator, which shows that $\ca{V}_0$ is indeed locally finitely presentable. It follows that the Cauchy completion $\overline{\ca{B}}$ of $\ca{B}$ is small (see \cite{JOHNSEN}), hence it makes sense to speak about the Tannakian adjunction in this context. Moreover, Theorem~\ref{RECOGNITION_THM} can be applied. Condition i) coincides with condition i) in Theorem~\ref{RECOGNITION_THM}. To see that iii) holds, we first note that $\ca{A}$ has cotensor products\footnote{Cotensor products are the dual notion of tensor products.} with objects in $\ca{X}$: since any $X\in \ca{X}$ has a dual $X^\vee$, any $\ca{V}$-functor preserves tensor products with $X$ (see \cite{STREET_ABSOLUTE}). In particular we have natural isomorphisms
\[
 \ca{A}(B,X^\vee \odot A) \cong X^\vee \otimes \ca{A}(B,A) \cong[X,\ca{A}(B,A)],
\]
 which shows that $\ca{A}$ is $\ca{X}$-cotensored. Since $\ca{X}\subseteq \ca{V}$ is $\Set$-dense, this implies that the notion of conical colimit in $\ca{A}$ coincides with the notion of ordinary colimit in $\ca{A}_0$ (cf.\ \cite[\S~3.8]{KELLY_BASIC}). If we let $\mathbf{\Phi}$ be the class of conical weights, the above observation, Corollary~\ref{CONICAL_COLIMIT_COR} and iii) imply that condition iii) of Theorem~\ref{RECOGNITION_THM} is satisfied. Therefore it only remains to check that ii) implies that $L_\omega \colon \Prs{A} \rightarrow \Prs{B}$ preserves the necessary coequalizers. We will in fact show that $L_\omega$ preserves all \emph{finite} weighted limits (see \cite[\S~4]{KELLY_FINLIM} for a definition of finite weighted limits).

 To see this, we first note that the functor
\[
(\ev_B)_{B\in \ca{B}} \colon \Prs{B} \rightarrow \prod_{B\in \ca{B}} \ca{V} 
\]
 preserves and reflects all limits and all colimits because both limits and colimits in $\Prs{B}$ are computed pointwise. Hence, it suffices to check that for a fixed $B\in \ca{B}$, the functor $\ev_B L_\omega \colon \Prs{A} \rightarrow \ca{V}$ preserves finite limits, i.e., that $\ev_B L_\omega$ is \emph{left exact} in the terminology of \cite{KELLY_FINLIM}. Since the functor $(\ev_B)_{B\in \ca{B}}$ preserves all colimits, it also preserves left Kan extensions (see \cite[Proposition~4.14]{KELLY_BASIC}). Thus, $\ev_B L_\omega$ is naturally isomorphic to the left Kan extension of $\ev_B \omega$ along the Yoneda embedding $Y \colon \ca{A} \rightarrow \Prs{A}$. This reduces the problem to showing that $\ev_B\omega$ is \emph{flat}, and by \cite[\S~6.3]{KELLY_FINLIM} it suffices to check that $\ev_B \omega \in [\ca{A},\ca{V}]$ is a filtered (conical) colimit of representable functors.

 We have seen that $\ca{A}$ is $\ca{X}$-cotensored, i.e., that $\ca{A}^{\op}$ is $\ca{X}$-tensored. Corollary~\ref{CONICAL_COLIMIT_COR}, applied to the $\ca{V}$-category $\ca{A}^{\op}$ and the contravariant Yoneda embedding $Y^\prime \colon \ca{A}^{\op} \rightarrow [\ca{A}, \ca{V}]$, shows that $\ev_B \omega$ is isomorphic to the conical colimit of the domain functor $\ca{A}^{\op} \slash \ev_B \omega \rightarrow [\ca{A},\ca{V}]$. Using the weak Yoneda lemma (see \cite[\S~1.9]{KELLY_BASIC}) one can show that the category $\ca{A}^{\op} \slash \ev_B \omega$ is isomorphic to the opposite of the category of elements of $V(\ev_B \omega)_0$. The latter is filtered by assumption ii), hence the remarks in \cite[\S~6.3]{KELLY_FINLIM} show that $\ev_B \omega$ is indeed flat. The second statement follows immediately: the functor $\widetilde{\omega}$ is right adjoint, so it preserves all limits, and we have just shown that $L_\omega$ preserves finite limits.
\end{proof}

\begin{cor}\label{COMMUTATIVE_RING_COR}
 Let $R$ be a commutative ring, $B$ an $R$-algebra, and let $\ca{V}$ be the category $\Mod_R$ of $R$-modules. Let $\ca{A}$ be a small additive $R$-linear category, equipped with an $R$-linear functor $\omega \colon \ca{A} \rightarrow \Mod_B$ into the category of left $B$-modules such that $\omega(A)$ is finitely generated and projective for all $A \in \ca{A}$. If 
\begin{enumerate}
 \item[i)] The functor $\omega$ reflects isomorphisms,
 \item[ii)] The category $\el(\omega)$ of elements of $\omega$ is cofiltered,
 \item[iii)] The functor $\omega$ detects and preserves those (ordinary) colimits which are finitely generated projective,
\end{enumerate}
 then the unit $N \colon \ca{A} \rightarrow (\Mod_B)^c_{L(\omega)}$ of the Tannakian adjunction is an equivalence of categories. Moreover, the comonad $L(\omega) \colon \Mod_{B} \rightarrow \Mod_{B}$ preserves finite limits, i.e., the corresponding $B$-$B$-bimodule is flat as a right $B$-module.
\end{cor}

\begin{proof}
 We let $\ca{X}$ be the full monoidal subcategory of $\Mod_R$ consisting of finitely generated free modules, which is clearly a dense autonomous generator. Since $\ca{A}$ is additive, it is $\ca{X}$-tensored: the tensor product of $A\in \ca{A}$ with $R^n$ is simply the $n$-fold direct sum $\oplus_{i=1}^n A$. The Cauchy completion of $B^{\op}$, considered as a one-object $\ca{V}$-category, is the full subcategory of $\Mod_B$ consisting of finitely generated projective modules. Thus the image of $\omega(A)$ is an object of the Cauchy completion of $B^{\op}$, and it makes sense to speak of the $(\ca{A},\omega)$-component of the unit of the Tannakian adjunction (cf.\ Definition~\ref{TANNAKIAN_ADJUNCTION_DFN}). The remaining conditions are precisely the conditions in Theorem~\ref{DENSE_RECOGNITION_THM}.
\end{proof}

\subsection[DAG_REMARKS_SECTION]{}\label{DAG_REMARKS_SECTION}
 We conclude with a few remarks about the necessity of the conditions in Theorem~\ref{DENSE_RECOGNITION_THM}. More precisely, we are interested in the following question: under what assumptions on the cosmos $\ca{V}$, the small category $\ca{B}$ and the comonad $T \colon \Prs{B}\rightarrow \Prs{B}$ does the functor $V_T \colon \Prs{B}^c_T \rightarrow \Prs{B}$ satisfy the conditions in Theorem~\ref{DENSE_RECONSTRUCTION_THM}? Since this theorem is a consequence of Theorem~\ref{RECOGNITION_THM}, the general remarks from Section~\ref{RECOGNITION_REMARKS_SECTION} apply: conditions i) and iii) hold without any further assumptions. If $\overline{\ca{B}} \subseteq \Prs{B}$ is closed under finite limits and if $T \in \CC{B}$ preserves finite limits, then the functor $V^c_T \colon \Prs{B}_T^c \rightarrow \Prs{B}$ satisfies condition ii) of Theorem~\ref{DENSE_RECONSTRUCTION_THM}. Indeed, it is well known that $\Prs{B}_T$ has finite limits under these conditions, and that $V_T$ preserves them. Since $\overline{\ca{B}}$ is closed under finite limits, so is $\Prs{B}_T^c$, and it follows that the category of elements of $V(\ev_B V^c_T)_0$ is cofiltered for every $B \in \ca{B}$. In the special case of Corollary~\ref{COMMUTATIVE_RING_COR}, where $\ca{V}=\Mod_R$ for some commutative ring $R$ and $\ca{B}$ is an $R$-algebra $B$, seen as a one-object $\ca{V}$-category, the category $\overline{\ca{B}}$ is simply the category of finitely generated projective $B$-modules. It is closed under finite limits if $B$ is a \emph{hereditary} Noetherian algebra, i.e., if submodules of finitely generated projective $B$-modules are again finitely generated projective.

\section{The category of filtered modules}\label{F_MODULES_SECTION}
\subsection[F_MODULES_OVERVIEW_SECTION]{}\label{F_MODULES_OVERVIEW_SECTION}
 In this section, we give an application of the generalized theory of Tannakian duality developed so far to the category of filtered modules introduced by Fontaine and Laffaille in \cite{FONTAINE_LAFFAILLE}. We fix a perfect field $k$ of characteristic $p>0$, and we let $W$ be the ring of Witt vectors with coefficients in $k$. For our purposes it suffices to know that $W$ is a discrete valuation ring with residue field $k$ which contains the ring of $p$-adic integers $\mathbb{Z}_p$, and that $p \in \mathbb{Z}_p$ is a uniformizer of $W$. A construction of the ring can be found in \cite[\S~II.6]{SERRE}. There is an automorphism $\sigma \colon W \rightarrow W$ of $\mathbb{Z}_p$-algebras which lifts the Frobenius automorphism on the residue field $k$ of $W$ (see \cite[Theorem~II.8 and Proposition~II.10]{SERRE}). This automorphism $\sigma$ is again called the Frobenius automorphism. For a $W$-module $M$, we write $M_\sigma$ for the $W$-module obtained by base change along $\sigma$. More concretely, we have $M_\sigma = M$, and the $W$-action on $M_\sigma$ is given by $w\cdot m=\sigma(w)m$. For a morphism $g\colon M \rightarrow M^\prime$ of $W$-modules, we have $g_{\sigma}(m)=g(m)$ for all $m\in M$. In the following definition we use the same notation and nomenclature which was introduced in \cite{WINTENBERGER}. We write $W_n$ for the quotient ring $W/p^n W$.

\begin{dfn}
 A \emph{filtered $F$-Module}\footnote{A filtered $F$-module is a filtered $W$-module with additional structure; the $F$ is part of the name and does not stand for a ring.} consists of
\begin{itemize}
\item
 a $W$-module $M$ with a decreasing filtration $(\Fil^i M)_{i \in \mathbb{Z}}$ of submodules $\Fil^i M \subseteq M$. The filtration is \emph{exhaustive}, $\bigcup_{i\in \mathbb{Z}} \Fil^i M = M$, and \emph{separated}, $\bigcap_{i\in \mathbb{Z}} \Fil^i M =0$;
\item
 for each $i\in \mathbb{Z}$, a morphism $\varphi^i \colon \Fil^i M \rightarrow M_\sigma$ of $W$-modules such that the restriction of $\varphi^i$ to $\Fil^{i+1} M$ is $p\varphi^{i+1}$.
\end{itemize}
 A morphism of filtered $F$-modules $M \rightarrow M^\prime$ is a morphism $g\colon M\rightarrow M^\prime$ of $W$-modules such that for all $i\in \mathbb{Z}$, $g(\Fil^i M ) \subseteq \Fil^i M^\prime$ and $\varphi^i_{M^{\prime}} \circ g = g_\sigma \circ \varphi^i_M$. We denote the category of filtered $F$-modules by $\MF$, and we write $\MF_{fl}$ for the full subcategory of objects $M$ which satisfy
\begin{itemize}
 \item
 the $W$ module $M$ has finite length;
 \item
 the images of the $\varphi^i$ span $M$, i.e., $\sum_{i\in \mathbb{Z}} \varphi^i(\Fil^i M) =M$.
\end{itemize}
The category $\MF^n_{fl}$ is the full subcategory of $\MF_{fl}$ consisting of those objects whose underlying $W$-module $M$ is annihilated by $p^n$, i.e., for which $M$ is a $W_n$-module. We write $\MF^n_{\proj}$ for the full subcategory of $\MF^n_{fl}$ consisting of objects whose underlying module is a finitely generated projective $W_n$-module.
\end{dfn}

\subsection[M_BAR_SECTION]{}\label{M_BAR_SECTION}
 J.-M.\ Fontaine and G.\ Laffaille have shown that $\MF_{fl}$ is an abelian $\mathbb{Z}_p$-linear category, and that the forgetful functor $\omega \colon \MF_{fl} \rightarrow \Mod_W$ is an exact $\mathbb{Z}_p$-linear functor \cite{FONTAINE_LAFFAILLE}. It follows immediately that $\MF^n_{fl}$ is an abelian $\mathbb{Z}/p^n\mathbb{Z}$-linear category, and that $\omega$ restricts to a $\mathbb{Z}/p^n \mathbb{Z}$-linear functor $\omega \colon \MF^n \rightarrow \Mod_{W_n}$. Thus, $\MF^n_{\proj}$ is $\mathbb{Z}/p^n \mathbb{Z}$-linear, and we can further restrict $\omega$ to a functor on $\MF^n_{\proj}$ whose image is contained in the category of finitely generated projective $W_n$-modules. The goal of this section is to show that the $(\MF^n_{\proj},\omega)$-component of the unit of the Tannakian adjunction is an equivalence of categories. Thus the category $\MF^n_{\proj}$ is equivalent to the category of comodules of a certain comonoid in the category of $W_n$-$W_n$-bimodules. Such a comonoid is precisely a $W_n$-$W_n$-coalgebra (cf.\ Section~\ref{F_MODULES_INTRO_SECTION}). Since we will use Corollary~\ref{COMMUTATIVE_RING_COR}, it will also follow that the $W_n$-$W_n$-coalgebra in question is flat as a right $W_n$-module. In order to check the conditions in Corollary~\ref{COMMUTATIVE_RING_COR}, we will use an alternative description of the category $\MF_{fl}$.

 We write $\ca{F}_{W}$ for the category of $W$-modules with an exhaustive $\mathbb{Z}$-filtration. Note that we do not assume that the filtration of an object in $\ca{F}_W$ is separated. Given a filtered $W$-module $M$ we write $\overline {M}$ for the colimit of the diagram
\[
 \xymatrix@H=17pt@!R=40pt@!C=15pt{&  \cdots \ar@{_{(}->}[ld] \ar[rd]^{p\cdot} & & \Fil^i M \ar@{_{(}->}[ld] \ar[rd]^{p\cdot}  & & \Fil^{i+1} M \ar@{_{(}->}[ld] \ar[rd]^{p\cdot} & & \cdots \ar@{_{(}->}[ld] \ar[rd]^{p\cdot} \\
\cdots && \Fil^{i-1}M && \Fil^i M && \Fil^{i+1} M && \cdots
}
\]
 where $i$ runs over all the integers. A morphism of filtered $W$-modules induces a natural transformation between the corresponding diagrams, so the assignment $M \mapsto \overline{M}$ extends to a functor $\overline{(-)} \colon \ca{F}_W \rightarrow \Mod_W$.

\begin{lemma}\label{FILTRATION_LEMMA}
 Let $M$ be a filtered $W$-module whose filtration is exhaustive and separated. If $\len(M)$ is finite, then $\len(\overline{M})$ is finite and $\len(\overline{M})=\len(M)$.
\end{lemma}

\begin{proof}
 This is \cite[Lemma~{1.7}]{FONTAINE_LAFFAILLE}.
\end{proof}

 It is evident from the definition of $\MF_{fl}$ that giving an object $M$ of $\MF_{fl}$ is equivalent to giving a $W$-module $M$ of finite length together with an exhaustive separated filtration $(\Fil^i M)_{i\in \mathbb{Z}}$ and a surjective morphism $\varphi \colon \overline{M} \rightarrow M_\sigma$. By the above lemma, such a morphism is surjective if and only if it is an isomorphism. Giving a morphism $g \colon M \rightarrow M^\prime$ in the category $\MF_{fl}$ is equivalent to giving a homomorphism $g$ of $W$-modules such that the diagram
\[
\xymatrix{\overline{M} \ar[r]^{\overline{g}} \ar[d]_{\varphi} & \overline{M^\prime} \ar[d]^{\varphi^\prime}\\
M_\sigma \ar[r]_{g_\sigma} & M^\prime_\sigma}
\]
is commutative.

\begin{lemma}\label{F_MODULE_COLIMIT_LEMMA}
 The forgetful functor $\omega \colon \MF^n_{fl} \rightarrow \Mod_{W_n}$ detects and preserves colimits which have finite length, i.e., if $G \colon \ca{D} \rightarrow \MF^n_{fl}$ is a diagram such that the colimit of $\omega G$ exists and has finite length, then the colimit of $G$ exists and is preserved by $\omega$.
\end{lemma}

\begin{proof}
 The category $\MF^n_{fl}$ already has finite colimits. Because any colimit can be written as a filtered colimit of its finitely generated partial colimits (see \cite[Proposition~2.13.7]{BORCEUX}), it suffices to show the claim for the case where $\ca{D}$ is a filtered category. Let $M$ be the colimit of $\omega G$, and let $M^i$ be the colimit of the diagram $\Fil^i G \colon \ca{D} \rightarrow  \Mod_{W_n}$. From exactness of filtered colimits in $\Mod_{W_n}$ it follows that the inclusions $\Fil^i Gd \subseteq Gd$ induce inclusions of the colimits, i.e., we have $M^i \subseteq M$. We define a filtration on $M$ by $\Fil^i M \defl M^i$. This filtration is exhaustive: since $M$ has finite length, it is in particular finitely generated. The category $\ca{D}$ is filtered, so one of the structure maps $Gd \rightarrow M$ of the colimit must be surjective. From the definition of $\MF^n_{fl}$ it follows that there is an integer $j$ such that $\Fil^j Gd=Gd$, and it follows that $M^j=M$. The module $M$ endowed with this filtration is clearly the colimit of the diagram $\ca{D} \rightarrow \ca{F}_W$ of underlying filtered $W$-modules. Since colimits in diagram categories are computed pointwise it is also immediate that $\overline M$ is the colimit of $d \mapsto \overline{Gd}$. The isomorphisms $\overline{Gd} \rightarrow (Gd)_\sigma$ therefore induce an isomorphism of $W_n$-modules $\varphi \colon \overline{M} \rightarrow M_\sigma$. It only remains to check that the filtration of $M$ is separated.

 Let $M^{\infty}=\cap_{i\in \mathbb{Z}} M^i$. Since $\len(M)< \infty$, we must have $M^j =M^\infty$ for $j$ sufficiently large. The module $M/M^{\infty}$ has a natural exhaustive $\mathbb{Z}$-filtration given by $\Fil^i M/M^{\infty}=M^{i}/M^{\infty}$, and since $M^j=M^{\infty}$ for $j$ large enough, this filtration is separated. From Lemma~\ref{FILTRATION_LEMMA} we get that $\len(\overline{M/M^\infty})=\len(M/M^{\infty})$. For $i \in \mathbb{Z}$ arbitrary we let $\varphi^i$ be the composite
\[
 \xymatrix@1{M^i \ar[r] & \overline{M} \ar[r]^{\varphi}& M_\sigma}
\]
 of $\varphi$ with the structure map of the colimit $\overline{M}$. Thus $\varphi^{i}(x)=p\varphi^{i+1}(x)$ for each $x\in M^{i+1}$. Now let $x \in M^{\infty}$. Since $M^\infty \subseteq M^{i+n}$ we have $\varphi^i(x)=p^n \varphi^{i+n}(x)$. But $p^n$ annihilates $M_{\sigma}$, so we find that $\varphi^i(x)=0$ for each $i\in \mathbb{Z}$ and each $x \in M^\infty$. We therefore get induced homomorphisms $\psi^i \colon M^i/M^{\infty} \rightarrow M_\sigma$, which in turn induce a factorization $\overline{M} \rightarrow \overline{M/M^{\infty}} \rightarrow M_\sigma$ of $\varphi \colon \overline{M} \rightarrow M_{\sigma}$. Since $\varphi$ is an isomorphism, we have a surjective homomorphism $\overline{M/M^\infty} \rightarrow M_\sigma$, which implies that $\len(M/M^{\infty})=\len(\overline{M/M^{\infty}}) \geq \len(M_\sigma)=\len(M)$. It follows that $\len(M^\infty)=0$, i.e., that $M^\infty=0$. Thus the filtration of $M$ is indeed separated, which shows that the module $M$ with filtration $\Fil^i M=M^i$ and structure morphisms $\varphi^i \colon M^i \rightarrow M_\sigma$ gives an object in $\MF^n_{fl}$. It is clear from the construction that $M$ is the colimit of $G$ in the category $\MF^n_{fl}$.
\end{proof}

 In order to prove that $\MF^n_{\proj}$ is the category of comodules for some comonoid, we need to introduce one more auxiliary category. We write $\MF_{fg}$ for the full subcategory of $\MF$ of filtered $F$-modules $M$ which satisfy
\begin{itemize}
 \item
 the $W$-module $M$ is finitely generated;
 \item
 the modules $\Fil^i M$ are direct summands of $M$;
 \item
 the images of the $\varphi^i$ span $M$, i.e., $\sum_{i\in \mathbb{Z}} \varphi^i(\Fil^i M) =M$.
\end{itemize}
 The following proposition was proved by J.-P.\ Wintenberger in \cite{WINTENBERGER}. It shows in particular that we have a sequence $\MF^n_{\proj} \subseteq \MF^n_{fl} \subseteq \MF_{fl} \subseteq MF_{fg}$ of full subcategories.

\begin{prop} \label{F_MODULE_COFILTERED_PROP}
 The category of filtered $F$-modules has the following properties.
\begin{enumerate}
 \item[i)]
 The category $\MF_{fg}$ is abelian, and the forgetful functor $\omega \colon \MF_{fg} \rightarrow \Mod_W$ is an exact $\mathbb{Z}_p$-linear functor.
 \item [ii)]
 For any object $M$ of $\MF_{fl}$, the filtration by submodules consists of direct summands. Thus $\MF_{fl}$ can be identified with the full subcategory of $\MF_{fg}$ consisting of objects which are annihilated by some power of $p$.
 \item[iii)]
 For any object $M$ of $\MF_{fg}$ there exists an object $M^\prime$ of $\MF_{fg}$ and an epimorphism $g \colon M^{\prime} \rightarrow M$ such that the underlying $W$-module of $M^\prime$ is free.
\end{enumerate}
\end{prop}

\begin{proof}
 Parts i) and ii) are proved in \cite[Proposition~1.4.1]{WINTENBERGER}, and part iii) is \cite[Proposition~1.6.3]{WINTENBERGER}.
\end{proof}

\begin{thm} \label{FILTERED_MODULE_THEOREM}
 The $(\MF^n_{\proj},\omega)$-component of the unit of the Tannakian adjunction is an equivalence. Thus $\MF^n_{\proj}$ is equivalent to the category of left Cauchy comodules of the $W_n$-$W_n$-coalgebra
\[
 L = \int^{\MF^n_{\proj}} \omega(M) \otimes_{\mathbb{Z}/p^n \mathbb{Z}} \omega(M)^\vee\rlap{,}
\]
 where the right action on $L$ is induced by the $W_n$-actions on $\omega(M)^\vee$, and the left action is induced by the $W_n$-actions on $\omega(M)$. The $W_n$-$W_n$-coalgebra $L$ is flat as a right $W_n$-module.
\end{thm}

\begin{proof}
 We will first show that the conditions of Corollary~\ref{COMMUTATIVE_RING_COR} are satisfied for $R=\mathbb{Z}/p^n \mathbb{Z}$ and $B=W_n$, i.e., that the functor $\omega \colon \MF^n_{\proj} \rightarrow \Mod_{W_n}$ reflects isomorphisms, that $\omega$ detects and preserves those colimits in $\Mod_{W_n}$ which are finitely generated and projective, and that the category $\el(\omega)$ is cofiltered. The first fact is immediate from the alternative description of $\MF_{fl}$ given in Section~\ref{M_BAR_SECTION}, and detection and preservation of the colimits in question follows from Lemma~\ref{F_MODULE_COLIMIT_LEMMA}. It remains to to check that $\el(\omega)$ is cofiltered. Since $\MF^n_{\proj}$ has direct sums, it suffices to check that for any pair of morphisms $f,g \colon (M,x) \rightarrow (M^\prime,x^\prime)$ in $\el(\omega)$, there is an object $(N,y)$ in $\el(\omega)$ and a morphism $h \colon (N,y) \rightarrow (M,x)$ such that $fh=gh$. Let $K$ be the equalizer of $f,g$ in $\MF^n_{fl}$. We have $f(x)=x^\prime=g(x)$, so $x \in K$. From Proposition~\ref{F_MODULE_COFILTERED_PROP}, part iii) we know that there is an object $L$ of $\MF_{fg}$ with an epimorphism $k \colon L \rightarrow K$ such that the underlying $W$-module of $L$ is free. Multiplication with $p^n$ defines an endomorphism of $L$ in $\MF_{fg}$. The cokernel $N$ of this endomorphism is a free $W_n$-module. Since $K$ is an object of $\MF^n_{fl}$, it is annihilated by $p^n$, so we get a morphism $h\colon N \rightarrow K$ in $\MF_{fg}$ making the diagram
\[
 \xymatrix{L \ar[r]^{p^n} \ar[rd]_{0} & L \ar[r] \ar[d]^{k} & N \ar[ld]^{h} \\
& K }
\]
 commutative. Surjectivity of the morphism $k$ implies that $h$ is surjective. In particular, there is an element $y \in N$ with $h(y)=x$. Since $N$ is a finitely generated free $W_n$-module, this gives the desired morphism $h \colon (N,y) \rightarrow (M,x)$ in the category of elements of $\omega \colon \MF^n_{\proj} \rightarrow \Mod_{W_n}$. By Corollary~\ref{COMMUTATIVE_RING_COR} it follows that the $(\MF^n_{\proj},\omega)$-component of the unit of the Tannakian adjunction is an equivalence, and that $L(\omega) \colon \Mod_{W_n} \rightarrow \Mod_{W_n}$ is left exact. By \cite[Formula~4.25]{KELLY_BASIC}, we have
\begin{align*}
 L(\omega)(X) & \cong \int^{\MF^n_{\proj}} \omega(M) \otimes_{\mathbb{Z}/p^n\mathbb{Z}} \Mod_{W_n} \bigl(\omega (M),X\bigr)\\
 &\cong \int^{\MF^n_{\proj}} \omega(M) \otimes_{\mathbb{Z}/p^n \mathbb{Z}} \omega(M)^{\vee} \otimes_{W_n} X \rlap{,}
\end{align*}
 (cf.\ Proposition~\ref{KAN_EXTENSION_PROP} and Section~\ref{WEIGHTED_COLIMIT_SECTION}) which shows that $L_\omega \widetilde{\omega}$ is naturally isomorphic to $L\otimes_{W_n}(-)$.
\end{proof}

\section{The Tannakian biadjunction}\label{TANNAKIAN_BIADJUNCTION_SECTION}

\subsection[BI_ADJUNCTION_OVERVIEW_SECTION]{}\label{BI_ADJUNCTION_OVERVIEW_SECTION}
 In this section we assume some more 2-categorical background. As already mentioned in the introduction, the reconstruction and recognition problems for comonoids are only part of classical Tannakian duality. We give an outline on how one could approach the reconstruction and recognition problems for comonoids with additional structure. More precisely, we show how to approach them from a categorical point of view.

 A bialgebra or bimonoid is an object of $\ca{V}$ which is both a monoid and a comonoid in a compatible way. There is a monoidal structure on the category of comonoids in $\ca{V}$ such that bimonoids in $\ca{V}$ are precisely the monoids in the category of comonoids. We would like to introduce a monoidal structure on the category $\Vcat\slash \ca{V}^c$ such that monoids in $\Vcat\slash \ca{V}^c$ are small monoidal $\ca{V}$-categories equipped with a strong monoidal $\ca{V}$-functor to $\ca{V}^c$. The goal is to show that the Tannakian adjunction becomes a monoidal adjunction, and consequently lifts to an adjunction between monoids in the respective categories. However, a small monoidal $\ca{V}$-category is only a monoid `up to natural isomorphism', so we first need to show that the Tannakian adjunction is suitably compatible with natural transformations.

\subsection[2_CATEGORIES_SECTION]{}\label{2_CATEGORIES_SECTION}
 The target category $\Vcat\slash \overline{\ca{B}}$ of the comodule functor is a 2-category in a natural way. A 2-cell from $F\colon (\ca{A},\omega) \rightarrow (\ca{A}^{\prime},\omega^{\prime})$ to $G \colon (\ca{A},\omega) \rightarrow (\ca{A}^{\prime},\omega^{\prime})$ is a $\ca{V}$-natural transformation $\alpha \colon F \Rightarrow F^{\prime}$ such that $\omega^{\prime} \alpha=\id$, i.e., such that the equality
\[
\vcenter{
\xymatrix{
\ca{A} \rrtwocell^{F}_{F^{\prime}}{\alpha} \ar[rd]_{\omega} && \ca{A}^{\prime} \ar[ld]^{\omega^{\prime}} \\
& \overline{\ca{B}}}}=
\vcenter{
\xymatrix{
\ca{A} \ar[rr]^{F} \ar[rd]_{\omega} && \ca{A}^{\prime} \ar[ld]^{\omega^{\prime}} \\
& \overline{\ca{B}}}}
\]
 of pasting composites holds. We think of $\CC{B}$ as a 2-category with no non-identity 2-cells.

\begin{prop}\label{TANNAKIAN_2_ADJUNCTION}
 With 2-cells defined as in Section~\ref{2_CATEGORIES_SECTION}, the Tannakian adjunction becomes a 2-adjunction. In other words, the assignment from Proposition~\ref{TANNAKIAN_ADJUNCTION_PROP} gives a 2-natural isomorphism of categories
\[
\CC{B}\bigl(L(\omega), T\bigr)\cong
\Vcat \slash \overline{\ca{B}}\bigl((\ca{A},\omega),(\Prs{B}^c_T, V^c_T)\bigr)
\]
(see Definition~\ref{TANNAKIAN_ADJUNCTION_DFN} for the notation).
\end{prop}

\begin{proof}
 It suffices to check that there are no non-identity two-cells between any two 1-cells in $\Vcat \slash \overline{\ca{B}}\bigl((\ca{A},\omega),(\Prs{B}^c_T, V^c_T)\bigr)$. Thus, let $F,F^\prime \colon (\ca{A},\omega) \rightarrow (\Prs{B}^c_T, V_T^c)$ be two such 1-cells. This implies that $V^c_T \circ F=V^c_T \circ F^{\prime}=\omega$. Let $\overline{\alpha} \colon F \Rightarrow F^\prime$ be a 2-cell. By Proposition~\ref{COACTION_LIFT_PROP}, the $\ca{V}$-natural transformation $\overline{\alpha}$ is uniquely determined by $\alpha \defl V^c_T \circ \overline{\alpha}$. But $\alpha$ is equal to the identity natural transformation $\omega \Rightarrow \omega$ by definition of the 2-category structure (see Section~\ref{2_CATEGORIES_SECTION}). Again by Proposition~\ref{COACTION_LIFT_PROP}, it follows that the coactions corresponding to $F$ and $F^\prime$ are equal, and that $\overline{\alpha}$ is the identity. The fact that this isomorphism is 2-natural is now trivial. We have to show that an equality between 2-cells holds, and since there are no non-identity 2-cells, any two 2-cells with the same domain and codomain are equal.
\end{proof}

 There is an obvious candidate for a monoidal structure on $\Vcat\slash \ca{V}^c$: given $\ca{V}$-functors $\omega \colon \ca{A} \rightarrow \ca{V}^c$ and $\omega^\prime \colon \ca{A}^{\prime} \rightarrow \ca{V}^c$ we can form a new $\ca{V}$-functor
\[
 \xymatrix{\ca{A}\otimes \ca{A}^{\prime} \ar[r]^-{\omega \otimes \omega^{\prime}} & \ca{V}^c \otimes \ca{V}^c \ar[r]^-{\otimes} & \ca{V}^c \rlap{.}}
\]
 However, there is in general no way to define associativity constraints for this tensor product in $\Vcat\slash\ca{V}^c$. To remedy this we have to introduce a new 2-category as target of the Tannakian adjunction. In Section~\ref{ISOFIBRATION_SECTION} we show that the comodule functor lands in a full sub-2-category $\ca{K}$ of $\Vcat\slash \overline{\ca{B}}$. In Proposition~\ref{BIEQUIVALENCE_PROP} we will show that the 2-category $\ca{K}$ is biequivalent to another 2-category $\Vcat\sslash \overline{\ca{B}}$, which has a nice monoidal structure in the case where $\overline{\ca{B}}=\ca{V}^c$. The resulting composite will only be a biadjunction, which is a weakening of the notion of a 2-adjunction. The definitions of biequivalences and biadjunctions can be found in \cite[Section~1]{STREET_FIBRATIONS}.

\subsection[ISOFIBRATION_SECTION]{}\label{ISOFIBRATION_SECTION}
 Let $(\ca{A},\omega)$ be an object of $\Vcat\slash \overline{\ca{B}}$. We call $\omega$ an \emph{isofibration} if it has the following property: for any two $\ca{V}$-functors $F\colon \ca{X} \rightarrow \ca{A}$, $G \colon \ca{X} \rightarrow \overline{\ca{B}}$ together with a $\ca{V}$-natural isomorphism $\beta \colon G \Rightarrow \omega F$, there exists a $\ca{V}$-functor $\overline{G} \colon \ca{X} \rightarrow \ca{A}$ and a $\ca{V}$-natural isomorphism $\overline{\beta} \colon \overline{G} \Rightarrow F$ with $\omega \overline{G}=G$ and $\omega \overline{\beta}=\beta$. We write $\ca{K}$ for the full sub-2-category of $\Vcat\slash \overline{\ca{B}}$ whose objects are the isofibrations.

 For $T \in \CC{B}$, the functor $V^c_T \colon \Prs{B}^c_T \rightarrow \overline{\ca{B}}$ is an isofibration. Indeed, by Proposition~\ref{COACTION_LIFT_PROP}, a $\ca{V}$-functor $F \colon \ca{X} \rightarrow \Prs{B}^c_T$ corresponds to a coaction $\varrho \colon V_T F \Rightarrow T V_T F$, and $\varrho^\prime \defl T\beta \circ \varrho \circ \beta^{-1} \colon G \Rightarrow TG$ defines a coaction on $G \colon \ca{X} \rightarrow \overline{\ca{B}}$ which satisfies $T\beta \circ \varrho = \varrho^\prime \circ \beta$. Applying Proposition~\ref{COACTION_LIFT_PROP} again we get the desired lifts. In other words, the image of the comodule 2-functor is contained in the full sub-2-category $\ca{K}$.

\subsection[WEAK_SLICE_SECTION]{}\label{WEAK_SLICE_SECTION}
 The 2-category $\Vcat\sslash\overline{\ca{B}}$ has the same objects as $\Vcat \slash \overline{\ca{B}}$, but the morphisms are given by triangles
\[
 \xymatrix{\ca{A} \ar[rd]_-{\omega} \ar[rr]^{F} & \dtwocell\omit{*!<0pt,-13pt>{\sigma}} & \ca{A}^{\prime} \ar[ld]^-{\omega^{\prime}}\\
 &  \overline{\ca{B}}}
\]
 which commute up to natural isomorphism $\sigma \colon \omega^\prime F \Rightarrow \omega$. The 2-cells between $(F,\sigma)$ and $(F^{\prime},\sigma^{\prime})$ are the $\ca{V}$-natural transformations $\alpha \colon F \Rightarrow F^{\prime}$ for which the equality
\[
\vcenter{
\xymatrix@!R=25pt{
\ca{A} \rrtwocell^{F}_{F^{\prime}}{\alpha} \ar[rd]_{\omega} & \dtwocell\omit{\sigma^{\prime}} & \ca{A}^{\prime} \ar[ld]^{\omega^{\prime}} \\
& \overline{\ca{B}}}}=
\vcenter{
 \xymatrix{\ca{A} \ar[rd]_-{\omega} \ar[rr]^{F} & \dtwocell\omit{*!<0pt,-13pt>{\sigma}} & \ca{A}^{\prime} \ar[ld]^-{\omega^{\prime}}\\
 &  \overline{\ca{B}}}
}
\]
 holds.

 We write $I\colon \ca{K} \rightarrow \Vcat \sslash \overline{\ca{B}}$ for the natural inclusion 2-functor. In order to show that $I$ is a biequivalence we need a characterization of equivalences in the 2-category $\Vcat\sslash \overline{\ca{B}}$. Recall that a 1-cell $F \colon \ca{A} \rightarrow \ca{A}^{\prime}$ in a 2-category is an equivalence if there is a 1-cell $G \colon \ca{A}^{\prime} \rightarrow \ca{A}$ with invertible 2-cells $\eta \colon \id \Rightarrow GF$ and $\varepsilon \colon FG \Rightarrow \id$. These 2-cells can be chosen such that the triangular identities hold, in which case they form an \emph{adjoint} equivalence.

\begin{lemma}\label{EQUIVALENCE_LEMMA}
 A 1-cell $(F,\sigma) \colon (\ca{A},\omega) \rightarrow (\ca{A}^{\prime},\omega^{\prime})$ in $\Vcat \sslash \overline{\ca{B}}$ is an equivalence if and only if $F \colon \ca{A} \rightarrow \ca{A}^{\prime}$ is an equivalence in $\Vcat$.
\end{lemma}

\begin{proof}
 The `only if' part is trivial. Thus, let $F \colon \ca{A} \rightarrow \ca{A}^{\prime}$ be an equivalence, and choose an inverse equivalence $G \colon \ca{A}^{\prime} \rightarrow \ca{A}$ and $\ca{V}$-natural isomorphisms $\eta \colon \id \Rightarrow GF$ and $FG \Rightarrow \id$ which satisfy the triangular identities. We claim that for
\[
\tau \defl \vcenter{\xymatrix@!C=15pt{
 \ca{A}^{\prime} \ar[rd]_{\id} \ar[rr]^G & \dtwocell\omit{*!<0pt,-12pt>{\varepsilon}} & \ca{A} \ar[ld]|F \ar[rd]^{\omega} \\
& \ca{A}^{\prime} \ar[rr]_{\omega^{\prime}} &\utwocell\omit{^\sigma^{-1}}& \ca{A}\rlap{,}
}}
\]
 the 1-cell $(G,\tau)$ is an inverse equivalence of $(F,\sigma)$. It follows directly from the definition that the equality
\[
\vcenter{\xymatrix{
\ca{A}^{\prime} \rrtwocell^{FG}_{\id}{\varepsilon} \ar[rd]_{\omega^{\prime}} && \ca{A}^{\prime} \ar[ld]^{\omega^{\prime}} \\
& \overline{\ca{B}}}}=
\vcenter{\xymatrix@!C=25pt{
\ca{A}^{\prime} \ar[rd]_{\omega^{\prime}} \ar[r]^G & \ca{A} \ar[d]|{\omega} \ar[r]^F \dtwocell\omit{<3>*!<3pt,-13pt>{\tau}} \dtwocell\omit{<-3>*!<-3pt,-13pt>{\sigma}} & \ca{A}^{\prime} \ar[ld]^{\omega^{\prime}}\\
 & \overline{\ca{B}}
}}
\]
 holds, i.e., that $\varepsilon$ is a 2-cell in the 2-category $\Vcat \sslash \overline{\ca{B}}$. Since inverses are unique, one of the triangular identities implies that $\varepsilon F=F \eta^{-1}$. Using pasting composites we find that the equality
\[
\vcenter{\xymatrix@!C=25pt{
\ca{A} \ar[rd]_{\omega} \ar[r]^F & \ca{A}^{\prime} \ar[d]|{\omega^{\prime}} \ar[r]^G \dtwocell\omit{<3>*!<3pt,-13pt>{\sigma}} \dtwocell\omit{<-3>*!<-3pt,-13pt>{\tau}} & \ca{A}  \ar[ld]^{\omega} \\
& \overline{\ca{B}}
}}=
\vcenter{\xymatrix@!C=15pt{
\ca{A} \ar[rd]_{\omega} \ar@/_8pt/[rr]_{\id} \rrcompositemap<3>_{F}^{G}{<-0.5>*!<-6pt,0pt>{\eta^{-1}}} && \ca{A} \ar[ld]^{\omega}\\
& \overline{\ca{B}}
}}
\]
 holds, which implies that $\eta$ is a 2-cell $\id \Rightarrow (G,\tau) \circ (F,\sigma)$ in the 2-category $\Vcat \sslash \overline{\ca{B}}$. Thus $(G,\tau)$ is an inverse equivalence of $(F,\sigma)$, as claimed.
\end{proof}

\begin{prop}\label{BIEQUIVALENCE_PROP}
 The natural inclusion $I \colon \ca{K} \rightarrow \Vcat\sslash \overline{\ca{B}}$ (see Section~\ref{WEAK_SLICE_SECTION}) is a biequivalence.
\end{prop}

\begin{proof}
 We have to show that $I$ induces an equivalence on hom-categories and that every object in the codomain is equivalent to an object in the image of $I$. It is clear from the definition of 2-cells that $I$ is fully faithful on hom-categories (cf.\ Section~\ref{2_CATEGORIES_SECTION} and Section~\ref{WEAK_SLICE_SECTION}), and the definition of an isofibration implies that $I$ is essentially surjective on hom-categories. It remains to show that any object $(\ca{A},\omega)$ in $\Vcat\sslash\overline{\ca{B}}$ is equivalent to an isofibration. This follows from the fact that $\omega$ can be written as $\omega=\omega^{\prime} \circ u$ where $\omega^{\prime}$ is an isofibration and $u$ is an equivalence (see \cite[Section~3]{LACK}). From Lemma~\ref{EQUIVALENCE_LEMMA} we know that $u$ is an equivalence in the 2-category ${\Vcat \sslash \overline{\ca{B}}}$.
\end{proof}

\subsection[ADDITIONAL_STRUCTURE_SECTION]{}\label{ADDITIONAL_STRUCTURE_SECTION}
 In Section~\ref{ISOFIBRATION_SECTION} we have seen that the Tannakian 2-adjunction can be restricted to the full sub-2-category $\ca{K}$ whose objects are isofibrations. It follows from Proposition~\ref{BIEQUIVALENCE_PROP} that the comodule functor with codomain the 2-category $\Vcat\sslash\overline{\ca{B}}$ (i.e., the composite $I\circ \Prs{B}^c_{(-)}$) has a left bi-adjoint. Moreover, any property that is stable under bi-equivalence holds for this biadjunction if and only if it holds for the Tannakian 2-adjunction. Thus, we can use our construction in Section~\ref{TANNAKIAN_ADJUNCTION_SECTION} to deduce facts about the Tannakian biadjunction. For example, the $T$-component of the counit of the biadjunction is an equivalence if and only if the $T$-component of the counit of the Tannakian 2-adjunction is, and a 1-cell in $\CC{B}$ is an equivalence if and only if it is an isomorphism (recall that $\CC{B}$ has no non-identity 2-cells, so all equivalences are isomorphisms). Thus, the reconstruction results in Section~\ref{RECONSTRUCTION_SECTION} apply to both the Tannakian 2-adjunction and the Tannakian biadjunction. 

 The 2-category $\Vcat\sslash \ca{V}^c$ is a monoidal 2-category. The tensor product of $(\ca{A},\omega)$ and $(\ca{A}^{\prime},\omega^{\prime})$ is given by the composite
\[
 \xymatrix{\ca{A}\otimes \ca{A}^{\prime} \ar[r]^-{\omega \otimes \omega^{\prime}} & \ca{V}^c \otimes \ca{V}^c \ar[r]^-{\otimes} & \ca{V}^c \rlap{.}}
\]
 Note that the problem in defining the associativity isomorphisms in $\Vcat\slash \ca{V}^c$ no longer exists, since we are looking at triangles which commute only up to natural isomorphism. Proposition~6.4 in \cite{MCCRUDDEN_MASCHKE} and the computations in \cite[\S~16]{STREET_QUANTUM_GROUPS} indicate that the left adjoint preserves the tensor product up to isomorphism, and thus lifts to a functor between pseudomonoids on on both sides. A pseudomonoid in the category $\Comon(\ca{V})$ is just a bimonoid, and a pseudomonoid in $\Vcat\sslash\ca{V}^c$ is precisely a monoidal $\ca{V}$-category $\ca{A}$, equipped with a strong monoidal $\ca{V}$-functor $\omega \colon \ca{A} \rightarrow \ca{V}^c$. By making this precise one could therefore extend our reconstruction results in Section~\ref{RECONSTRUCTION_SECTION} from comonoids to bimonoids. This is of course just a starting point for the reconstruction of additional structures on a comonoid.

\appendix

\section{Pasting composites}\label{PASTING_COMPOSITES_APPENDIX}
 Pasting composites were introduced in \cite{KELLY_STREET}, where they were defined in the more general context of arbitrary strict 2-categories. They provide a convenient way to streamline computations with $\ca{V}$-natural transformations. For $\ca{V}$-functors and $\ca{V}$-natural transformations as in
\[
\xymatrix{{\ca{A}} \ar[r]^L & {\ca{B}} \rtwocell^F_G{\alpha} & {\ca{C}} \ar[r]^K & {\ca{D}}}
\quad\text{and}\quad
\xymatrix{
\ca{A} \ruppertwocell^F{\beta} \ar[r]^(0.3)G \rlowertwocell_H{\gamma} & \ca{B}\smash{\rlap{,}}\\
}
\]
 we write $K \alpha L$ for the $\ca{V}$-natural transformation from $KFL$ to $KGL$ whose component at the object $A$ is given by $K\alpha_{LA}$. We say that $K\alpha L$ is obtained by \emph{whiskering} $\alpha$ with $K$ and $L$. We call the $\ca{V}$-natural transformation from $F$ to $H$ with $A$-component given by the composite $\gamma_A \beta_A$ the \emph{vertical composite} of $\gamma$ and $\beta$, and we denote it by $\gamma \beta$. Whiskering and vertical composition are called the \emph{basic pasting operations}.

 Given a diagram of $\ca{V}$-functors and $\ca{V}$-natural transformations such as
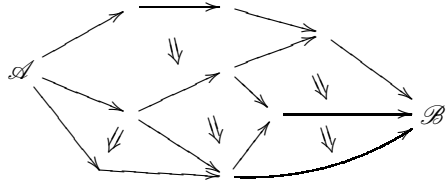
\[
\xymatrix@R=4pt@C=12pt{
      & & \ar[rr] & & \ar[rrd] &            &          & &        \\
      & &         & &          &            & \ar[2,2] & &        \\
\ca{A}
\ar[-2,2]
\ar[1,2]
\xtwocell[-2,6]{}\omit{<1>}
\xcompositemap[3,4]{}<-4>{<1.5>}
      & &         & & \ar[rru] \ar[rd]
                      \xtwocell[-2,3]{}\omit{<4>}
                               &            &          & &        \\
      & & \ar[rru] \ar[2,2]
          \xtwocell[-1,3]{}\omit{<3>}
                  & &          & \ar[rrr]   &          & & \ca{B} \\
      & &         & &          &            &          & &        \\
      & &         & & \ar[ruu]
                      \xlowertwocell[-2,4]{}<-3>{<-0.5>}
                                            &          & &        \\
}
\]
we can use the basic pasting operations to get a $\ca{V}$-natural transformation going from the composite of the $\ca{V}$-functors on the top of the diagram to the composite of those on the bottom as follows: first, we choose any $\ca{V}$-natural transformation whose domain is contained in the top chain of the diagram. Then we `split' the diagram along the codomain of this $\ca{V}$-natural transformation:
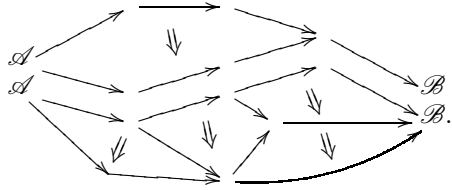
\[
\xymatrix@R=5pt@C=12pt{
      & & \ar[rr] & & \ar[rrd] &            &          & &        \\
      & &         & &          &            & \ar[2,2] & &        \\
{\smash{\ca{A}}}
\ar[-2,2]
\ar[1,2]
\xtwocell[-2,6]{}\omit{<0.5>}
      & &         & & \ar[rru] &            & \ar[2,2] & &        \\
{\smash{\ca{A}}}
\ar[1,2]
\xcompositemap[3,4]{}<-3>{<1>}
      & & \ar[rru] & & \ar[rru] \ar[rd]
                      \xtwocell[-2,3]{}\omit{<3>}
                               &            &          & & {\smash{\ca{B}}} \\
      & & \ar[rru] \ar[2,2]
          \xtwocell[-1,3]{}\omit{<2>}
                  & &          & \ar[rrr]   &          & &
                  {\smash{\ca{B}}\smash{\rlap{.}}} \\
      & &         & &          &            &          & &        \\
      & &         & & \ar[ruu]
                      \xlowertwocell[-2,4]{}<-3>{<-0.5>}
                                            &          & &        \\
}
\]
We proceed by whiskering the diagram on the top to obtain a $\ca{V}$-natural transformation between $\ca{V}$-functors with domain $\ca{A}$ and codomain $\ca{B}$, and then iterate this whole process with the rest of the diagram. The \emph{pasted composite} of the diagram is the vertical composite of the resulting collection of $\ca{V}$-natural transformations. Usually this process involves choices, namely whenever there are several $\ca{V}$-natural transformations whose domains are contained in the top chain of the diagram. The resulting $\ca{V}$-natural transformation is independent of these choices (see \cite{POWER}).

Note that a diagram of $\ca{V}$-functors is commutative if and only if one can place the identity natural transformation in the diagram. We therefore introduce the following useful convention from \cite{KELLY_STREET}: when we compute the pasted composite of a diagram of $\ca{V}$-categories, $\ca{V}$-functors and $\ca{V}$-natural transformations between them, if the diagram has parts containing no natural transformation, these parts must be commutative and they are treated as identity natural transformations.

\section{Density in cosmoi with dense autonomous generator} \label{DENSITY_IN_COSMOI_WITH_DAG_APPENDIX}
\subsection[ACTEGORY_SECTION]{}\label{ACTEGROY_SECTION}
 Let $\ca{V}$ be a cosmos with dense autonomous generator $\ca{X}$ (see Definition~\ref{DAG_DFN}). To each $\ca{X}$-tensored $\ca{V}$-category $\ca{A}$ we can associate an ordinary category endowed with an action of $\ca{X}_0$. Such a category is called a \emph{$\ca{X}_0$-actegory} (in \cite{MCCRUDDEN_REPR_COALGEBROIDS}) or \emph{$\ca{X}_0$-representation} (in \cite{GORDON_POWER}). An $\ca{X}_0$-representation is an ordinary category $\ca{L}$, together with a functor $-\odot - \colon \ca{X}_0 \times \ca{L} \rightarrow \ca{L}$ and natural isomorphisms $l\colon L \rightarrow I \odot L$ and $a \colon X \odot (X^\prime \odot L) \rightarrow (X\otimes X^\prime) \odot L$ for all $L\in \ca{L}$, subject to certain coherence conditions (details can be found in \cite[Section~2]{GORDON_POWER} or \cite[Section~3]{MCCRUDDEN_REPR_COALGEBROIDS}). Since we assume that $\ca{X}_0$ is $\Set$-dense, the assignment which sends an $\ca{X}$-tensored $\ca{V}$-category $\ca{A}$ to the $\ca{X}_0$-representation $\ca{A}_0$, with action given by the tensor functor $-\odot- \colon \ca{X}_0 \times \ca{A}_0 \rightarrow \ca{A}_0$ is in fact a fully faithful 2-functor (see \cite[Theorem~3.4]{GORDON_POWER}). This means that giving a $\ca{V}$-functor $F \colon \ca{A} \rightarrow \ca{A}^\prime$ between $\ca{X}$-tensored $\ca{V}$-categories is the same as giving an ordinary functor $F_0 \colon \ca{A}_0 \rightarrow \ca{A}^\prime_0$, together with morphisms $\widehat{F} \colon X\odot F_0 A \rightarrow F_0(X\odot A)$ making the diagrams
\[
\vcenter{\xymatrix{
F_0 A \ar[rd]_{F_0 l} \ar[r]^{l^\prime} & I \odot F_0 A \ar[d]^{\widehat{F}}\\
& F_0(I\odot A)
}}
\quad\mathrm{and}\quad
\vcenter{\xymatrix{ X \odot (X^\prime \odot F_0 A) \ar[r]^{a^\prime}  \ar[d]_{X \odot \widehat{F}} & (X\otimes X^\prime) \odot F_0 A \ar[dd]^{\widehat{F}} \\
X\odot F_0(X^\prime\odot A) \ar[d]_{\widehat{F}} &\\
F_0\bigl(X\odot(X^\prime \odot A)\bigr) \ar[r]_{F_0 a} & F_0\bigl((X\otimes X^\prime) \odot A\bigr)
}}
\]
commutative. The arrows $l$ and $a$ correspond to the canonical isomorphisms $\id \cong [I,-]$ and $[X,[X^\prime,-]] \cong [X\otimes X^\prime,-]$ under the $\ca{V}$-natural isomorphisms which define the respective tensor products, and $\widehat{F}$ is given by the map of the same name introduced in Section~\ref{WEIGHTED_COLIMIT_SECTION}. Moreover, we know that tensor products with objects in $\ca{X}$ are absolute colimits (see \cite{STREET_ABSOLUTE}), so the morphisms $\widehat{F} \colon X\odot FA \rightarrow F(X\odot A)$ are isomorphisms. Still under the assumption that $\ca{X}_0$ is $\Set$-dense and that $\ca{A}$, $\ca{A}^\prime$ are $\ca{X}$-tensored, giving a $\ca{V}$-natural transformation $\alpha \colon F \Rightarrow F^\prime \colon \ca{A} \rightarrow \ca{A}^\prime$ is the same as giving an ordinary natural transformation $\alpha \colon F_0 \Rightarrow F^\prime_0$ such that
\[
 \xymatrix{X\odot F_0A \ar[r]^{X\odot \alpha_A} \ar[d]_{\widehat{F}} & X \odot F^\prime_0 A \ar[d]^{\widehat{F^\prime}}\\
F_0(X\odot A) \ar[r]_{\alpha_{X\odot A}} & F^\prime_0 (X\odot A)}
\]
is commutative.

\begin{prop}\label{DENSITY_PROP}
 Let $\ca{V}$ be a cosmos which has a dense autonomous generator $\ca{X}$. Let $\ca{A}$ be an $\ca{X}$-tensored $\ca{V}$-category, and let $\ca{C}$ be a $\ca{V}$-category which is cotensored. A $\ca{V}$-functor $K \colon \ca{A} \rightarrow \ca{C}$ is $\ca{V}$-dense if and only if the underlying ordinary functor $K_0 \colon \ca{A}_0 \rightarrow \ca{C}_0$ is $\Set$-dense.
\end{prop}

\begin{proof}
 The assumption that $\ca{C}$ is cotensored implies that $K$ is $\ca{V}$-dense if and only if the map $\ca{C}_0(C,D) \rightarrow \VNat\bigl(\ca{C}(K-,C),\ca{C}(K-,D)\bigr)$ which sends $g \colon C \rightarrow D$ to the $\ca{V}$-natural transformation $\ca{C}(K-,g)$ is a bijection of sets (see \cite[Section~5.1]{KELLY_BASIC}). For $C\in \ca{C}$, let $K\slash C$ be the category with objects the morphisms $\varphi \colon KA \rightarrow C$, $A \in \ca{A}$, and morphisms $\varphi \rightarrow \varphi^\prime$ the morphisms in $\ca{A}_0$ which make the evident triangle commutative. From \cite[Formula~5.4]{KELLY_BASIC} we know that $K_0$ is $\Set$-dense if and only if each object $C$ is the colimit of the tautological cocone on the functor $V_C \colon K\slash C \rightarrow \ca{C}$ which sends $\varphi$ to its domain. We write $S_D$ for the set of cocones on $V_C$ with vertex $D$, and we let $V=\ca{V}_0(I,-) \colon \ca{V}_0 \rightarrow \Set$ be the canonical forgetful functor. Let $\chi \colon \VNat\bigl(\ca{C}(K-,C), \ca{C}(K-,D)\bigr) \rightarrow S_D$ be the map which sends $\alpha$ to the cocone $\chi(\alpha)\defl\bigl(V\alpha_A (\varphi)\bigr)_{\varphi \in K \slash C}$. The composite
\[
 \xymatrix{ \ca{C}_0(C,D) \ar[r] & \VNat \bigl(\ca{C}(K-,C),\ca{C}(K-,D)\bigr)  \ar[r]^-{\chi} & S_D}
\]
 sends $g$ to the cocone $(g \varphi)_{\varphi \in K\slash C}$. This composite is a bijection if and only if $C$ is the colimit of the tautological cocone, i.e., if and only if $K_0$ is $\Set$-dense. If we can show that $\chi$ is a bijection, then $K_0$ is $\Set$-dense if and only if $K$ is $\ca{V}$-dense, as claimed.

 We now construct an inverse for $\chi$, as follows. Given a cocone $\gamma=(\gamma_\varphi)_{\varphi \in K\slash C}$, we let
\[
\beta_A \colon \ca{C}_0(KA,C) \rightarrow \ca{C}_0(KA,D)
\]
 be the map with $\beta_A(\varphi)=\gamma_\varphi$. We write $F, G \colon \ca{A} \rightarrow \ca{V}^{\op}$ for the functors $\ca{C}(K-,C)$ and $\ca{C}(K-,D)$ respectively. Note that we have $VF_0=\ca{C}_0(KA,C)$, and $\beta$ is a natural transformation between the $\Set$-valued functors $VF_0$ and $VG_0$. We first use the density assumption to lift this to a natural transformation $\xi(\gamma) \colon F_0 \rightarrow G_0$ between the underlying ordinary $\ca{V}_0$-valued functors of $F$ and $G$, and we then show that $\xi(\gamma)$ is in fact $\ca{V}$-natural. The tensor product of $B$ by $X$ in $\ca{V}^{\op}$ is given by $[X,B]$. Since all $\ca{V}$-functors preserve tensor products with objects which have duals (see \cite{STREET_ABSOLUTE}), we get isomorphisms
\[
\xymatrix{F(X\odot A) \ar[r]^-{\widehat{F}} & [X,FA]}\quad\mathrm{and}\quad \xymatrix{G(X\odot A) \ar[r]^-{\widehat{G}} & [X,GA] },
\]
 and the composite $V\widehat{G}_0 \circ \beta_{X\odot A} \circ V \widehat{F}_0^{-1} \colon V[X,\ca{C}(KA,C)]_0 \rightarrow V[X,\ca{C}(KA,D)]_0$ is natural in $X$. Since $\ca{V}_0(X,-)$ is naturally isomorphic to $V[X,-]_0$, it follows by $\Set$-density of $\ca{X}$ in $\ca{V}$ that there is a unique morphism $\xi(\gamma)_A\colon \ca{C}(KA,C) \rightarrow \ca{C}(KA,D)$ in $\ca{V}$ such that
\[
 \xymatrix@C=45pt{\ca{V}_0\bigl(X,\ca{C}(KA,C)\bigr) \ar[d]_{\cong}  \ar[r]^-{\ca{V}_0(X,\xi(\gamma)_A)} & \ca{V}_0\bigl(X,\ca{C}(KA,D)\bigr) \ar[d]^{\cong} \\
 V[X,\ca{C}(KA,C)]_0 \ar[d]_{V\widehat{F}_0^{-1}} \ar[r]^-{V[X,\xi(\gamma)_A]_0} & V[X,\ca{C}(KA,D)]_0 \ar[d]^{V\widehat{G}_0^{-1}}\\
 \ca{C}_0\bigl(K(X\odot A),C\bigr) \ar[r]_-{\beta_{X\odot A}} & \ca{C}_0\bigl(K(X\odot A),D\bigr) }
\]
 is commutative for every $X \in \ca{X}$. Hence part (1) and (3) of the diagram
\[
 \xymatrix@!C=45pt{V[X,[X^{\prime},FA]]_0 \ar[rrrr]^{V[X,[X^\prime,\xi(\gamma)_A]]_0} && \ar@{}[d]|(0.4){(0)} && V[X,[X^{\prime},FA]]_0 \\
 & V[X,F(X\odot A)]_0 \ar[lu]_(0.4){V[X,\widehat{F}]_0} \ar[rr]^{V[X,\xi(\gamma)_{X^\prime \odot A}]_0} &\ar@{}[d]|(0.4){(1)}& V[X,G(X\odot A)]_0 \ar[ru]^(0.4){V[X,\widehat{G}]_0}\\
 & VF_0\bigr(X\odot(X^\prime\odot A)\bigl) \ar[u]^{V\widehat{F}_0} \ar[rr]^{\beta_{X\odot(X^\prime \odot A)}} &\ar@{}[d]|{(2)}& VG_0\bigr(X\odot(X^\prime\odot A)\bigl) \ar[u]_{V\widehat{G}_0}\\
 & VF_0\bigr((X\otimes X^\prime)\odot A)\bigl) \ar[u]^{VF_0 a} \ar[ld]^(0.4){V\widehat{F}_0} \ar[rr]_{\beta_{(X\otimes X^\prime)\odot A}} && VG_0\bigr((X\otimes X^\prime)\odot A\bigl) \ar[u]_{VG_0 a} \ar[rd]_(0.4){V\widehat{G}_0}\\
V[X\otimes X^{\prime},FA]_0 \ar[uuuu]^{V a^\prime} \ar[rrrr]_-{V[X\otimes X^\prime,\xi(\gamma)_A]_0}&&\ar@{}[u]|(0.4){(3)}&& V[X\otimes X^{\prime},GA]_0 \ar[uuuu]_{V a^\prime} }
\]
 are commutative. Part (2) is commutative since $\beta$ is natural, and the two pentagons are instances of the coherence diagrams in Section~\ref{ACTEGROY_SECTION}. The outer diagram is commutative because $a^\prime$ is natural, hence it follows that part (0) is commutative. Since $\ca{X}$ is $\Set$-dense we find that $[X^\prime,\xi(\gamma)_A] \circ \widehat{F}=\widehat{G} \circ \xi(\gamma)_{X^\prime \odot A}$. The considerations in Section~\ref{ACTEGROY_SECTION} therefore imply that $\xi(\gamma)$ is a $\ca{V}$-natural transformation. The second coherence diagram of Section~\ref{ACTEGROY_SECTION} implies that $V\xi(\gamma)_A=\beta_A$ for all objects $A \in \ca{A}$, and it follows that $\chi\bigl(\xi(\gamma)\bigr) = \gamma$. Moreover, if we start with a $\ca{V}$-natural transformation $\alpha \colon F \Rightarrow G$ and construct the $\beta_A$ associated to the cocone $\chi(\alpha)$, we clearly get $\beta_{X\odot A}=V\alpha_{X \odot A}$, i.e., $V\xi\bigl(\chi(\alpha)\bigr)_{X\odot A}=V\alpha_{X \odot A}$. Both $\alpha$ and $\xi(\gamma)$ are $\ca{V}$-natural, hence we must have $V[X,\xi\bigl(\chi(\alpha)\bigr)_A]_0=V[X,\alpha_A]_0$, and by density of $\ca{X}$ it follows that $\xi\bigl(\chi(\alpha)\bigr)=\alpha$. In other words, the assignment which sends a cocone $\gamma$ to the $\ca{V}$-natural transformation $\xi(\gamma)$ constructed above gives the desired inverse to $\chi$.
\end{proof}

\begin{cor}\label{CONICAL_COLIMIT_COR}
 Let $\ca{V}$ be a cosmos which has a dense autonomous generator $\ca{X}$. Let $\ca{A}$ be a small $\ca{X}$-tensored $\ca{V}$-category. Fix a presheaf $F \in \Prs{A}$, and let $\ca{A}\slash F$ be the category of representable functors over $F$. Then $F$ is the conical colimit of the domain functor $D \colon \ca{A} \slash F \rightarrow \Prs{A}$.
\end{cor}

\begin{proof}
 The fact that the Yoneda embedding is always $\ca{V}$-dense and completeness of $\Prs{A}$ imply that the conditions of Proposition~\ref{DENSITY_PROP} are satisfied. It follows that $F$ is the ordinary colimit of the tautological cocone on $D \colon \ca{A}\slash F \rightarrow \Prs{A}_0$. But $\Prs{A}$ is cotensored, hence the notion of conical colimit and ordinary colimit coincide.
\end{proof}

\bibliographystyle{amsalpha}
\bibliography{tannaka}

\end{document}